\documentclass[pdflatex,sn-mathphys-num]{sn-jnl}

\usepackage{graphicx}
\usepackage{multirow}
\usepackage{amsmath,amssymb,amsfonts}
\usepackage{amsthm}
\usepackage{mathrsfs}\usepackage[title]{appendix}\usepackage{xcolor}\usepackage{textcomp}\usepackage{manyfoot}\usepackage{booktabs}\usepackage{algorithm}\usepackage{algorithmicx}\usepackage{algpseudocode}\usepackage{listings}

\theoremstyle{thmstyleone}\newtheorem{theorem}{Theorem}\newtheorem{proposition}{Proposition}\newtheorem{lemma}{Lemma}

\theoremstyle{thmstyletwo}\newtheorem{example}{Example}\newtheorem{remark}{Remark}

\theoremstyle{thmstylethree}\newtheorem{definition}{Definition}\newtheorem{assumption}{Assumption}

\usepackage{amssymb,amsmath,amsfonts,latexsym}
\usepackage{pstricks,pst-plot,pst-node,pst-text,pst-3d}
\usepackage{enumerate}
\usepackage[english]{babel}
\usepackage{geometry,booktabs,fancyhdr}
\usepackage{amsmath,latexsym,mathrsfs,textcomp,fontenc,yhmath}
\usepackage{amssymb,amsfonts,amscd,amssymb,amsthm}
\usepackage{graphics,verbatim,url}
\usepackage{epic,curves}
\usepackage{index}
\usepackage[T1]{fontenc}
\usepackage{tabularx,ragged2e,booktabs,caption}
\usepackage{pgf,tikz}
\usepackage{mathrsfs,mathtools}
\usetikzlibrary{arrows}
\usepackage{pgflibraryarrows}
\usepackage{pgflibrarysnakes}
\usetikzlibrary{calc}

\makeatletter
\newsavebox{\@brx}
\newcommand{\llangle}[1][]{\savebox{\@brx}{\(\m@th{#1\langle}\)}\mathopen{\copy\@brx\mkern2mu\kern-0.9\wd\@brx\usebox{\@brx}}}
\newcommand{\rrangle}[1][]{\savebox{\@brx}{\(\m@th{#1\rangle}\)}\mathclose{\copy\@brx\mkern2mu\kern-0.9\wd\@brx\usebox{\@brx}}}
\makeatother

\DeclareMathOperator{\divv}{{div \/}}

\DeclareMathOperator*{\esssup}{ess\,sup} \newcommand{\ac}[1]{\textcolor{black}{#1}}

\begin{document}

\title[Gevrey class regularity for {parametric} Navier-Stokes equations]{Gevrey class regularity for steady-state incompressible Navier-Stokes equations {in parametric domains and related models} }

\author*[1]{\fnm{Alexey} \sur{Chernov}}\email{alexey.chernov@uni-oldenburg.de}

\author[1]{\fnm{T{\`u}ng} \sur{L{\^e}}}\email{tung.le@uni-oldenburg.de}

\affil[1]{\orgdiv{Institut f{\"u}r Mathematik}, \orgname{Carl von Ossietzky Universit{\"a}t Oldenburg}, \orgaddress{\street{Ammerl{\"a}nder Heerstra{\ss}e 114-118}, \city{Oldenburg}, \postcode{26129}, \country{Germany}}}

\abstract{We investigate {parametric} Navier-Stokes equations for {a viscous, incompressible flow in  bounded domains}. The coefficients {of the equations}
are perturbed by high-dimensional random parameters,
{this fits in particular for modelling flows in domains with uncertain perturbations.}
Our focus is on deriving bounds for arbitrary high-order derivatives of the pressure and the velocity fields with respect to the random parameters in the context of incompressible Navier-Stokes equation under a small-data assumption. To achieve this, we analyze mixed and saddle-point problems and employ the \emph{alternative-to-factorial technique} to establish generalized Gevrey-class regularity for the solution pair. Thereby the analytic regularity follows as a special case. In the numerical experiments, we validate and illustrate our theoretical findings using Gauss-Legendre quadrature and Quasi-Monte Carlo methods.}



\maketitle

\section{Introduction} \label{sec: intro}
The Navier-Stokes equations are fundamental partial differential equations that describe the motion of incompressible fluid flows, encompassing both liquids and gases. They are central to a wide range of applications, including pipe flow, laminar and turbulent flows \cite{Sperone23,John2016}, airflow over an airfoil \cite{PlotnikovSokolowski13}, ocean currents \cite{ConstantinJohnson23,ConstantinJohnson19}, and cloud simulation \cite{Chertock2019,Chertock2023}. Numerical simulation and solution of the Navier-Stokes equations have become increasingly important for reducing development costs in design processes. However, real-world circumstances often unpredictably diverge from ideal computational environments, leading to random deviations in the actual solutions of the equations compared to simulation results. To further optimize manufacturing processes and minimize costs, it is crucial to quantify these uncertainties and their impact on the outcomes.

In order to introduce the parametric regularity approach, let $\boldsymbol{y}=(y_1,y_2,\dots)\in U$ be a random variable, where $\boldsymbol{y}$ has either finitely many or countably many components, {and $U$ is a bounded parameter domain.} We are interested in uncertainty quantification
for the steady state incompressible Navier-Stokes equations in mixed variational problem: Find $(\boldsymbol{u},p)\in \mathcal{H}\times\mathcal{L}$ such that
\begin{align}\label{NS gen equation}
	a(\boldsymbol{u},\boldsymbol{v};A(\boldsymbol{y}))+b(
	\boldsymbol{v},p;B(\boldsymbol{y}))+m( \boldsymbol{u},\boldsymbol{u},\boldsymbol{v};M(\boldsymbol{y}))
	&=\left\langle{\boldsymbol{f}(\boldsymbol{y})},{\boldsymbol{v}}
	\right\rangle
	,\quad \forall\boldsymbol{v}\in {\mathcal H}\notag\\
	b(
	\boldsymbol{u},q;B(\boldsymbol{y}))
	&=\llangle {g(\boldsymbol{y})},{q}
	\rrangle
	,\quad \forall q\in \mathcal L 
\end{align}
where ${\mathcal H}$ and $\mathcal L$ are Hilbert spaces, $a(\cdot,\cdot;A(\boldsymbol{y})):{\mathcal H}\times {\mathcal H} \rightarrow \mathbb R$ is {a} continuous, ${\mathcal H}$-elliptic, and symmetric bilinear form, and $b(\cdot,\cdot,B(\boldsymbol{y})): {\mathcal H}\times\mathcal L \rightarrow \mathbb R$ is continuous and {satisfies} the inf-sup condition
\begin{align}
	{\inf_{q\in \mathcal L}}\, \sup_{\boldsymbol{v}\in {\mathcal H}}\, \frac{b( \boldsymbol{v},q;B(\boldsymbol{y}))}{\left\|{\boldsymbol{v}}\right\|_{{\mathcal H}}\left\|{q}\right\|_{\mathcal L}}
	\geq \beta.
\end{align}
{The} trilinear form $m( \cdot,\cdot,\cdot;M(\boldsymbol{y}))$ is continuous on ${\mathcal H}\times{\mathcal H}\times{\mathcal H}\rightarrow\mathbb R$.  The dependence of the bilinear forms $a, b$ and the trilinear form $m$ on matrices $A,B$  and $M$  is assumed to be linear. We denote by  ${\mathcal H}^*$ and  $\mathcal L^*$ the dual spaces of ${\mathcal H}$ and $\mathcal L$, respectively. {The single bracket $\left\langle{\cdot},{\cdot}
	\right\rangle$ represents the {duality product} between ${\mathcal H}^*$ and  ${\mathcal H}$ while the double bracket $\llangle {\cdot},{\cdot}
	\rrangle$ represents the {duality product} between $\mathcal L^*$ and  $\mathcal L$}.  The {prescribed forcing term} is denoted by  $ \boldsymbol{f}(\boldsymbol{y})\in {\mathcal H}^*$, whereas the source mass term is represented by $ g(\boldsymbol{y}) \in \mathcal L^*$. The model is motivated by bilinear forms originating from the weak formulation of partial differential equations, the standard examples for our considerations are as follows:

\begin{example}[Stochastic Steady-State Navier-Stokes Flow]
	Consider the Navier-Stokes equations on the domain $D = [-1,1] \times [0,2\pi]$ with Dirichlet boundary conditions on $\Gamma_{\text{Dirichlet}} = [-1,1] \times \{0,2\pi\}$ and periodic boundary conditions on $\Gamma_{\text{periodic}} = \{-1,1\} \times [0,2\pi]$:
	\begin{align*}
		-\nu(\boldsymbol{y})\Delta \boldsymbol{u} + (\boldsymbol{u}\cdot \nabla)\boldsymbol{u} + \nabla p &= \boldsymbol{f}(\boldsymbol{x}, \boldsymbol{y}) \quad \text{in } D,\\
		\nabla \cdot \boldsymbol{u} &= 0  \quad \text{in } D,\\
		\boldsymbol{u} &= \mathbf{0} \quad \text{on } \Gamma_{\text{Dirichlet}},\\
		\boldsymbol{u} &\text{ is periodic on } \Gamma_{\text{periodic}}.
	\end{align*}
	Here, the stochastic kinematic viscosity $\nu(\boldsymbol{y})$ and force term $\boldsymbol{f}(\boldsymbol{x}, \boldsymbol{y})$ depend on the parameter $\boldsymbol{y} \in U$. We consider the functional spaces as follows 
	\begin{align*}
		{\mathcal H}
		&=
		\left\{\boldsymbol{u}\in H^1(D)^d : \boldsymbol{u}=0 \text{ on } \Gamma_{\text{Dirichlet}}
			\text{ and }
			\boldsymbol{u} \text{ is periodic on } \Gamma_{\text{periodic}} \right\}, \\
		\mathcal L&=L^2({D})\setminus \mathbb R=\left\{q\in L^2(D):\displaystyle\int_D q =0\right\}	  .
	\end{align*}
	The corresponding bilinear forms are given by
	\begin{align*}
		a(\boldsymbol{u},\boldsymbol{v};A(\boldsymbol{y})) &= \nu(\boldsymbol{y}) \int_{D} \text{tr} \big(\nabla {\boldsymbol{u}}^{\top}  \nabla{\boldsymbol{v}} \big),\\
		b(\boldsymbol{v},p;B(\boldsymbol{y})) &= -\int_{D} \divv {\boldsymbol{v}} \,\, {p} ,\\
		m(\boldsymbol{u},\boldsymbol{v},\boldsymbol{w};M(\boldsymbol{y})) &= \int_{D} {\boldsymbol{u}}^{\top} \nabla {\boldsymbol{v}}^{\top}  {\boldsymbol{w}},\\
		\langle \boldsymbol{f}(\boldsymbol{y}), \boldsymbol{v} \rangle&= \int_{D} {\boldsymbol{f}(\boldsymbol{x}, \boldsymbol{y})} \cdot {\boldsymbol{v}}.
	\end{align*}
{	Here, we denote the gradient of a vector field by the Jacobian matrix $\nabla \boldsymbol{u} := \left( \frac{\partial u_i}{\partial x_j} \right)_{i,j=1}^d$, and the associated coefficient matrices are set to $A(\boldsymbol{y})=\nu(\boldsymbol{y}) I$ and $B(\boldsymbol{y})=M(\boldsymbol{y})=I$.}
\end{example}

\begin{example}[Parametric Domain {Perturbations}] \label{Ex: dom trans}
		We consider the Navier-Stokes equations for a viscous, incompressible fluid in a parameter-dependent bounded Lipschitz domain $D_{\boldsymbol{y}}$. The system is subject to homogeneous Dirichlet boundary conditions on $\partial D_{\boldsymbol{y}}$:
		\begin{align*}
				-\Delta \boldsymbol{u} + (\boldsymbol{u}\cdot \nabla)\boldsymbol{u} + \nabla p &= \boldsymbol{f} \quad \text{in } D_{\boldsymbol{y}},\\
				\nabla \cdot \boldsymbol{u} &= g  \quad \text{in } D_{\boldsymbol{y}},\\
				\boldsymbol{u} &= \mathbf{0} \quad \text{on } \partial D_{\boldsymbol{y}}.
			\end{align*}
		Here, $D_{\boldsymbol{y}}$ is obtained as the image of a fixed nominal Lipschitz domain $\widehat{D} \subset \mathbb{R}^d$, where $d \in \{2,3\}$, under a bijective transformation $\boldsymbol{T}_{\boldsymbol{y}}: \widehat{D} \to D_{\boldsymbol{y}}$. We consider the function spaces ${\mathcal H} = H^1_0(\widehat{D})^d$ and $\mathcal L = L^2(\widehat{D})\setminus \mathbb R $. Noting the plain pullback transformation $\widehat{\varphi}(\widehat{\boldsymbol{x}})=\varphi(\boldsymbol{x})$ for all $\boldsymbol{x}=\boldsymbol{T}_{\boldsymbol{y}}(\widehat{\boldsymbol{x}})$, the bilinear forms and given data are expressed as {integrals over the nominal domain}:
		\begin{align*}
					a(\widehat{\boldsymbol{u}},\widehat{\boldsymbol{v}};A(\boldsymbol{y})) &= \int_{\widehat{D}} \text{tr} \big(\widehat{\nabla} \widehat{\boldsymbol{u}}^{\top} A(\boldsymbol{y})^{\top} \widehat{\nabla}\widehat{\boldsymbol{v}} \big),\\
					b(\widehat{\boldsymbol{v}},\widehat{p};B(\boldsymbol{y})) &= -\int_{\widehat{D}} \text{tr} \big(\widehat{\nabla}\widehat{\boldsymbol{v}}  B(\boldsymbol{y})\widehat{p} \big),\\
					m(\widehat{\boldsymbol{u}},\widehat{\boldsymbol{v}},\widehat{\boldsymbol{w}};M(\boldsymbol{y})) &= \int_{\widehat{D}} (\widehat{\nabla} \widehat{\boldsymbol{v}}\, M(\boldsymbol{y}) \widehat{\boldsymbol{u}} )^{\top}  \widehat{\boldsymbol{w}},\\
					\left\langle \widetilde{\boldsymbol{f}}(\boldsymbol{y}), \widehat{\boldsymbol{v}} \right \rangle &= \int_{\widehat{D}} \widehat{\boldsymbol{f}} \cdot \widehat{\boldsymbol{v}} J_{\boldsymbol{y}},\\
					\llangle[\big] { \widetilde{g}(\boldsymbol{y})},{ \widehat{q}}
	\rrangle[\big] &= \int_{\widehat{D}} \widehat{g} \, \widehat{q} J_{\boldsymbol{y}}, 
				\end{align*}
		where $A(\boldsymbol{y}) = d\boldsymbol{T}_{\boldsymbol{y}}^{-1} d\boldsymbol{T}_{\boldsymbol{y}}^{-\top} J_{\boldsymbol{y}}$, $B(\boldsymbol{y}) = M(\boldsymbol{y}) = d\boldsymbol{T}_{\boldsymbol{y}}^{-1} J_{\boldsymbol{y}}$, and $J_{\boldsymbol{y}} = \det(d\boldsymbol{T}_{\boldsymbol{y}})$.
\end{example}

Note that in both cases mentioned above, the presence of the random parameter $\boldsymbol{y}$ renders the solution of the Navier–Stokes equations, the velocity field
$\boldsymbol{u}$ and the pressure $p$, also random.

In this paper, we establish the parametric regularity of solutions to the Navier-Stokes equations, where the coefficients and the physical domain depend on a parameter $\boldsymbol{y}$. The given data $A, B, M, f$, and $g$ are assumed to be infinitely differentiable functions of $\boldsymbol{y}$, belonging to the Gevrey class $G^{\delta}$ for some fixed $\delta \geq 1$. The scale of Gevrey classes is a nested scale with respect to $\delta$, serving as an intermediate space between class of analytic functions ${\mathcal A}$  and smooth functions $C^\infty$ :
\begin{equation} \label{Gdelta-nested}
	{\mathcal A} = G^{1} \subseteq G^{\delta} \subset G^{\delta'} \subset C^\infty, \qquad 1 \leq \delta < \delta'. 
\end{equation}
The class of analytic functions ${\mathcal A}$ is the simplest and most significant in this scale. This class has been previously studied in the context of shape holomorphy for the steady-state Navier-Stokes equations \cite{Cohen2018}. This work focuses on affine parametrization of the coefficients and establishes solution analyticity using elegant complex analysis techniques.

Beyond complex analysis, real-variable methods have also been employed as powerful tools for proving analytic regularity. For instance, \cite{KoLiYu2025} investigates unsteady 2D Navier-Stokes equations, providing bounds for mixed first-order derivatives, which are essential for rank-1 lattice rules in Quasi-Monte Carlo (QMC) methods. To obtain derivative bounds of arbitrary order without any loss of regularity, we introduce a modified approach in \cite{ChernovLe2024a,ChernovLe2024b}, referred to as the \emph{alternative-to-factorial technique}.

{The recent work by Harbrecht et al. \cite{Harbrecht2024} develops} an implicit function theorem in the Gevrey class, showing that the Gevrey regularity of parametric inputs is preserved in the solutions of a wide range of operator equations.
{In particular, \cite[Theorem 2.2]{Harbrecht2024} covers general Gevrey mappings in abstract Banach spaces, the authors also provide a specific example addressing scalar-valued semilinear elliptic PDEs.
 In the present paper, we utilize the \emph{alternative-to-factorial technique} to analyze parametric regularity of incompressible Navier-Stokes equations. Furthermore, we explicitly track all constants in the regularity estimates, which is necessary for computing the QMC generating vector.}

Recently, Djurdjevac et al. \cite{DjurdjevacKaarnioja2025} studied parametric regularity for the domain transform. In this paper, we not only consider the matrices $A(\boldsymbol{y}),B(\boldsymbol{y}), M(\boldsymbol{y})$ and the given data $\boldsymbol{f}(\boldsymbol{y}),g(\boldsymbol{y})$ in general form, but also we establish the parametric regularity for the plain pullback using falling factorial technique. Compared to  \cite{DjurdjevacKaarnioja2025}, our approach--presented in Section \ref{sec: Domain transform}--offers a {simple} proof with optimal constants, and introduces flexible tools that can be applied to a wider class of problems.

The paper is organized as follows. Section \ref{sec: Preliminary} introduces the essential notation, including multi-index, the falling factorial, and definitions of Gevrey classes and analytic functions. We also present results on the parametric regularity of matrix-valued mappings in Subsection \ref{sec: useful tool}, along with key properties of the multivariate Faà di Bruno formula in Subsection \ref{sec: faa di bruno}. In Section \ref{sec: well-posedness}, we define the functional spaces and norms used in the variational formulation of the Navier-Stokes problem, and summarize the properties of the associated operators in the mixed and saddle point formulations in \ref{sec: NS saddle point}. These properties are then used in Subsection \ref{sec: NS theory} to establish the well-posedness of the Navier-Stokes problem. Section \ref{sec: NS main result} contains our main theoretical contribution. Under general assumptions on the given data $A,B,M$, the force term $\boldsymbol{f}$ and source mass term $g$, we prove parametric regularity for the solution to the incompressible Navier-Stokes equations. In Section \ref{sec: Domain transform}, we show that, under mild assumptions, the standard pullback operator arising from domain transformations satisfies the framework developed in Section \ref{sec: NS main result}. In Section \ref{sec: NS num exp}, we present numerical experiments in the setting of domain transformations described in Example \ref{Ex: dom trans}. We provide two examples based on Gauss–Legendre quadrature and two using QMC methods to validate our theoretical results.

\section{Preliminaries} \label{sec: Preliminary}
\subsection{The falling factorial  and multi-index and notation} \label{sec: Multiindex}
The limitation of the real-variable inductive argument for nonlinear problems arises from the Leibniz product rule and the triangle inequality, {we refer to \cite[Section 2.1]{ChernovLe2024a} for a detailed discussion}. To address this {issue}, we use the \emph{alternative-to-factorial technique} introduced in \cite{ChernovLe2024a,ChernovLe2024b}.
We begin by recalling some elementary results on falling factorials.
For a non-negative integer $n \in \mathbb N_0$ {absolute value of the \emph{falling factorial} of $\tfrac{1}{2}$} is defined as
\begin{equation*} 
	\left[\tfrac{1}{2} \right]_{n} 
	:=
	\left\{
		\begin{array}{cl}
				1, & n=0, \\
				\tfrac{1}{2}, & n=1, \\
				\tfrac{1}{2}\cdot \tfrac{1}{2} \cdot \tfrac{3}{2} \dots \left|n-\tfrac{3}{2}\right|, & n \geq 2.
			\end{array}
				\right.
	.
\end{equation*}
Though non-standard, this notation is convenient for our analysis. A useful estimate is
\begin{align}\label{ff-estimates}
	 \left[\tfrac{1}{2} \right]_{n}  \leq n! \leq 2 \cdot 2^n \left[\tfrac{1}{2} \right]_{n} ,
	  \end{align}
which suffices for our purposes.

We adopt the standard {multi-index} notation commonly used in analysis, see e.g., \cite{BoutetDeMonvel1967,CohenDevoreSchwab2010}. The set of finitely supported sequences of nonnegative integers is denoted by
\begin{equation}\label{cF-def}
	\mathcal F := 
	\left\{
	\boldsymbol{\nu} =(\nu_1,\nu_2,\dots)~:~ \nu_j\in \mathbb N_0,
	\text{ and } \nu_j \neq 0
	\text{ for only a finite number of } j
	\right\} ,
\end{equation}
where the summation $\boldsymbol{\alpha} + \boldsymbol{\beta}$ and the partial order relations $\boldsymbol{\alpha} < \boldsymbol{\beta}$ and $\boldsymbol{\alpha} \leq \boldsymbol{\beta}$ of elements in $\boldsymbol{\alpha}, \boldsymbol{\beta} \in \mathcal F$ are understood component-wise. Similar to the approach in \cite{CONSTANTINE96}, we introduce a linear ordering on $\mathbb N^{\mathbb N}$. Specifically, we write $\boldsymbol{\ell}\prec\boldsymbol{\nu}$ if one of the following conditions is satisfied:
\begin{enumerate}[(i)]
	\item $\left|\boldsymbol{\ell}\right|<\left|\boldsymbol{\nu}\right|$
	\item $\left|\boldsymbol{\ell}\right|=\left|\boldsymbol{\nu}\right|$ and $\ell_1<\nu_1$; or
		\item $\left|\boldsymbol{\ell}\right|=\left|\boldsymbol{\nu}\right|$, $\ell_1=\nu_1,\dots,\ell_k=\nu_k$ and $\ell_{k+1}<\nu_{k+1}$ for some $k\geq 1$.
\end{enumerate}

In order to simplify expressions, {we also utilize the following notations}
\begin{align*}
	\boldsymbol{R}^{\boldsymbol{\nu}}:= \prod_{j\geq 1} R_j^{v_j},
		\qquad \qquad
		\boldsymbol{R}\cdot \boldsymbol{1}: = \sum_{j\geq 1} R_j,
			\qquad \qquad
	\left|\boldsymbol{\nu}\right|
	:=
	\sum_{j\geq 1} \nu_j, 
	\qquad \qquad
	\boldsymbol{\nu}!
	:=
	\prod_{j\geq 1} v_j!, 
\end{align*}
where $\boldsymbol{R}=\{R_j\}_{j\geq 1}$ is a sequence of positive real numbers and $\boldsymbol{\nu} = \left\{\nu_j\right\}_{j \geq 1}$ is a multi-index consisting of non-negative integers.

Although the expressions $\boldsymbol{R} \cdot \boldsymbol{1}$ and $\left|\boldsymbol{\nu}\right|$ both involve summation over a sequence, their meanings and contexts differ. The notation $\boldsymbol{R} \cdot \boldsymbol{1}$ is used for sequences of real numbers, while $\left|\boldsymbol{\nu}\right|$ is used for finitely supported multi-indices. In particular, $\left|\boldsymbol{\nu}\right|$ is finite if and only if $\boldsymbol{\nu} \in \mathcal F$.

 For every $\boldsymbol{\nu}\in \mathcal F$ supported in $\left\{1,2,\dots, n\right\}$, the corresponding partial derivatives with respect to the variables $\boldsymbol{y}$ are defined for any scalar function $u$ and a $d$-dimensional vector field $\boldsymbol{f}$ as follows
\begin{align*}
	\partial^{\boldsymbol{\nu}} u
	=
	\frac{\partial^{\left|\boldsymbol{\nu}\right|}u}
	{\partial y_1^{\nu_1} \partial y_2^{\nu_2} \dots \partial y_n^{\nu_n}},
	\qquad \qquad
	\partial^{\boldsymbol{\nu}} \boldsymbol{f}
	=
	\left(	\partial^{\boldsymbol{\nu}} f_i\right)_{i=1}^d.
\end{align*}
For two  {multi-indices} $\boldsymbol{\nu}, \boldsymbol{\eta} \in \mathcal F$, the binomial coefficient is given by
\begin{equation*}
\left(\boldsymbol{\nu} \atop \boldsymbol{\eta}\right)  = 
	\prod_{j\geq 1}  \left(\nu_j \atop \eta_j\right)
	=\frac{\boldsymbol{\nu}!}{\boldsymbol{\eta}!\,(\boldsymbol{\nu}-\boldsymbol{\eta})!
}.
\end{equation*}
This multi-index notation provides a convenient framework for handling multiparametric objects. The following technical lemma will be instrumental in our analysis, with its proof available in \cite[Section 7]{ChernovLe2024a} and \cite[Section 2]{ChernovLe2024b}.

\begin{lemma}
	Let $\delta \geq 1$,  $\boldsymbol{\nu}, \boldsymbol{\eta}\in \mathcal F$ be two multi-indices satisfying $\boldsymbol{\eta} \leq \boldsymbol{\nu}$ and $\boldsymbol{e} \mathcal F$ be a unit multi-index, i.e. $\left|\boldsymbol{e}\right|=1$, we have
	\begin{align}
		(|\boldsymbol{\nu}-\boldsymbol{\eta}|!)^{\delta-1} (|\boldsymbol{\eta}|!)^{\delta-1} 
		&\leq (|\boldsymbol{\nu}|!)^{\delta-1},\label{multiindex-est-1}
\\
		\sum_{ {0<\boldsymbol{\eta} < \boldsymbol{\nu}}}
		\left(\boldsymbol{\nu} \atop \boldsymbol{\eta}\right) 
		\left[\tfrac{1}{2} \right]_{|\boldsymbol{\eta}|} 
		\left[\tfrac{1}{2} \right]_{|\boldsymbol{\nu}-\boldsymbol{\eta}|} 
		&\leq
		 {2} [\tfrac{1}{2}]_{|\boldsymbol{\nu}|},\label{multiindex-est-7}
\\
		\sum_{ {0\leq \boldsymbol{\eta} \leq \boldsymbol{\nu}}}
		\left(\boldsymbol{\nu} \atop \boldsymbol{\eta}\right) 
		\left[\tfrac{1}{2} \right]_{|\boldsymbol{\eta}|} 
		\left[\tfrac{1}{2} \right]_{|\boldsymbol{\nu}-\boldsymbol{\eta}|} 
		&\leq
		{4} [\tfrac{1}{2}]_{|\boldsymbol{\nu}|},\label{multiindex-est-2}
\\
		 \sum_{0<\boldsymbol{\eta} \leq \boldsymbol{\nu}}
		\left(\boldsymbol{\nu} \atop \boldsymbol{\eta}\right) 
		\left[\tfrac{1}{2} \right]_{|\boldsymbol{\eta}|} 
		\left[\tfrac{1}{2} \right]_{|\boldsymbol{\nu}-\boldsymbol{\eta}|} 
		&\leq
		3 [\tfrac{1}{2}]_{|\boldsymbol{\nu}|},\label{multiindex-est-3}
\\
	\sum_{ {\boldsymbol{0}< \boldsymbol{\eta} \leq \boldsymbol{\nu}}}
	\left(\boldsymbol{\nu}\atop\boldsymbol{\eta}\right)
		\left[\tfrac{1}{2} \right]_{\left|\boldsymbol{\eta}\right|} 
	\left[\tfrac{1}{2} \right]_{\left|\boldsymbol{\nu}+\boldsymbol{e}-\boldsymbol{\eta}\right|} 
	&\leq {
	\left[\tfrac{1}{2} \right]_{\left|\boldsymbol{\nu}+\boldsymbol{e}\right|} },
	\label{multiindex-est-6}
\\
	\sum_{\boldsymbol{0}< \boldsymbol{\eta} \leq \boldsymbol{\nu}}
	\sum_{ {\boldsymbol{0}< \boldsymbol{\ell} < \boldsymbol{\eta}}}
	\left(\boldsymbol{\nu}\atop\boldsymbol{\eta}\right)
	\left(\boldsymbol{\eta}\atop \boldsymbol{\ell}\right)
	\left[\tfrac{1}{2} \right]_{\left|\boldsymbol{\eta}-\boldsymbol{\ell}\right|} 
	\left[\tfrac{1}{2} \right]_{\left|\boldsymbol{\ell}\right|} 
	\left[\tfrac{1}{2} \right]_{\left|\boldsymbol{\nu}+\boldsymbol{e}-\boldsymbol{\eta}\right|} 
	&\leq
	{2}\left[\tfrac{1}{2} \right]_{\left|\boldsymbol{\nu}+\boldsymbol{e}\right|} ,
	\label{multiindex-est-8}
	\\
	\sum_{\boldsymbol{0}< \boldsymbol{\eta} \leq \boldsymbol{\nu}}
	\sum_{ {\boldsymbol{0}< \boldsymbol{\ell} \leq \boldsymbol{\eta}}}
	\left(\boldsymbol{\nu}\atop\boldsymbol{\eta}\right)
	\left(\boldsymbol{\eta}\atop \boldsymbol{\ell}\right)
	\left[\tfrac{1}{2} \right]_{\left|\boldsymbol{\eta}-\boldsymbol{\ell}\right|} 
	\left[\tfrac{1}{2} \right]_{\left|\boldsymbol{\ell}\right|} 
	\left[\tfrac{1}{2} \right]_{\left|\boldsymbol{\nu}+\boldsymbol{e}-\boldsymbol{\eta}\right|} 
	&\leq
	 {3}\left[\tfrac{1}{2} \right]_{\left|\boldsymbol{\nu}+\boldsymbol{e}\right|} .
	\label{multiindex-est-5}
\end{align}
\end{lemma}
 
\subsection{Gevrey-class and analytic functions} \label{sec: Gevrey}
 The following definition of Gevrey-$\delta$ functions with countably many parameters will be used in our regularity analysis in Section \ref{sec: NS main result}.

\begin{definition} \label{def:G-delta-def}
	Let $\delta \geq 1$, $B$ be a Banach space, $I \subset \mathbb R^\mathbb N$ be an open domain and a function $f:I \to B$ be such that its $\boldsymbol{y}$-derivatives $\partial^{\boldsymbol{\nu}} f:I \to B$ are continuous for all $\boldsymbol{\nu} \in \mathcal F$. We say that the function $f$ is of class Gevrey-$\delta$ if for each $y_0 \in I$ there exist an open set $J\subseteq I$, and strictly positive constants $\boldsymbol{R} = (R_1,R_2,\dots) \subset \mathbb R_{>0}^{\mathbb N}$ and $C \in \mathbb R_{>0}$ that the derivatives of $f$ satisfy the bounds
	\begin{equation}\label{G-delta-def}
		\|\partial^{\boldsymbol{\nu}} f(\boldsymbol{y})\|_B \leq \frac{C}{\boldsymbol{R}^{\boldsymbol{\nu}}} (|\boldsymbol{\nu}|!)^\delta, \qquad \forall \boldsymbol{y} \in J, \quad \forall \boldsymbol{\nu} \in \mathcal F.
	\end{equation}
	In this case we write $f \in G^\delta(U,B)$.
\end{definition}
\begin{remark}
When $\boldsymbol{y} = (y_1,\dots,y_s)$ with $s<\infty$, we note that $|\boldsymbol{\nu}|! \leq s^{|\boldsymbol{\nu}|} \boldsymbol{\nu}!$. In this finite-dimensional setting, the Gevrey-$1$ class $G^1(\mathbb R^s, \mathbb{R})$ (or $G^1(\mathbb C^s, \mathbb{C})$) coincides with the space of analytic functions in $s$ real (or complex) variables; see e.g. \cite[Section 2.2]{Krantz2002}. Furthermore, the Gevrey scale $G^\delta$ increases monotonically with respect to $\delta$, as described in \eqref{Gdelta-nested}.
	\end{remark}
\begin{remark}\label{rem:def-G-delta-U}
	In case $U=\left[-\frac{1}{2}, \frac{1}{2}\right]^{\mathbb N}$, the topology in $U$ can be defined using the isomorphism mapping from $U$ to a compact set $U(\boldsymbol{\alpha})=\left\{\boldsymbol{w}\in \ell^{\infty}: \left|w_j\right|\leq \frac{1}{2}\left|\alpha_j\right|\right\}$ for some $\boldsymbol{\alpha}\in \ell^q$, where $q>1$ and $\left|\alpha_j\right|>0$ for all $j\in \mathbb N$ (see \cite[Section 2.2]{Gilbert2019}). For any $\boldsymbol{y}\in U$ and $\widetilde{\boldsymbol{y}}\in U(\boldsymbol{\alpha})$, the rescaling $\widetilde{\boldsymbol{y}}=\boldsymbol{\alpha} \odot \boldsymbol{y} = \left\{\alpha_j\, y_j\right\}_{j\in \mathbb N}$ and $\widetilde f(\widetilde{\boldsymbol{y}})= f\left(\boldsymbol{\alpha}^{-1}\odot(\boldsymbol{\alpha} \odot \boldsymbol{y})\right)=f(\boldsymbol{y})$ imply that $f \in G^{\delta}(U,B)$ if and only if $\widetilde f \in G^{\delta}(U(\boldsymbol{\alpha}),B)$.
\end{remark}

\subsection{Gevrey-$\delta$ regularity for matrix-valued mappings}\label{sec: useful tool}

In this subsection, we present essential tools for establishing parametric regularity in domain transformation and the solution of the Navier-Stokes equation.  We begin with the Frobenius norm, also known as the Hilbert-Schmidt norm, denoted by $\left|\cdot\right|_F$. For any matrix $A \in \mathbb{R}^{m \times n}$, the Frobenius norm is defined as
$$\left|A\right|_{F}:
=
\left(\sum_{i=1}^m \sum_{j=1}^n \left|A_{i,j}\right|^2\right)^{\frac{1}{2}}
=
\text{tr}(A^{\top}A)^{\frac{1}{2}}=\left(\sum_{i=1}^m \left|A_i\right|_2^2\right)^{\frac{1}{2}}=
\left(\sum_{i=1}^n \left|A_i^{\top}\right|_2^2\right)^{\frac{1}{2}}.$$
where $A_i$ and $A_i^{\top}$ is the $i$-th row and column of the matrix $A$, respectively. It's noteworthy that, in the special case of vectors, the Frobenius norm is identical to the standard Euclidean norm, denoted by $\left|\cdot\right|_2$.

 {The Proposition below outlines two well-known properties of the Frobenius norm following from  the Cauchy-Schwarz inequality. They will be extensively used in the upcoming proofs.}
\begin{proposition}\label{F norm prop}
	Suppose $A\in \mathbb R^{m\times n}$ and $B\in \mathbb R^{n\times k}$, there holds
	\begin{align}\label{Frobenius properties}
		\mathrm{tr}(A\,B)
		\leq
		\left|A\right|_F \left|B\right|_F \quad \text{and} \quad
		\left|A\,B\right|_F
		\leq 
		\left|A\right|_F \left|B\right|_F.
	\end{align}
\end{proposition}

\begin{remark}
	We state Proposition \ref{F norm prop} separately for $\text{tr}(AB)$ and $\left|AB\right|_F$ since {neither $\left|\text{tr}(AB)\right|$ nor $\left|AB\right|_F$ always dominates} the either. Indeed, for example, on the one hand, when $A=B=\left(\begin{smallmatrix}1 & 0 \\ 0 & 1\end{smallmatrix}\right)$, then we have $\left|\text{tr}(AB)\right|>\left|AB\right|_F$. On the other hand, when $A=\left(\begin{smallmatrix}1 & 0 \\ 0 & 1\end{smallmatrix}\right)$ and $B=\left(\begin{smallmatrix}0 & 1 \\ 1 & 0\end{smallmatrix}\right)$, there holds $\left|\text{tr}(AB)\right|<\left|AB\right|_F$.
\end{remark}

Next, we introduce the Sobolev space norms that will be used extensively in the subsequent sections. In particular, we define the $L^\infty$-norm for matrix-valued functions $\boldsymbol{H}$ over a domain $D$ as
\begin{align}\label{W k infty def}
	\left\|{\boldsymbol{H}}\right\|_{L^{\infty}({D})}
	&=
	\esssup_{\boldsymbol{x}\in{D}}
	\left|\boldsymbol{H}(\boldsymbol{x})\right|_{F}.
\end{align}
This norm allows us to measure the essential supremum of the Frobenius norm of $\boldsymbol{H}$ across the domain $D$. Utilizing this Sobolev norm, we proceed to establish the Gevrey-$\delta$ regularity of parametric matrix-valued mappings, as detailed below.
\begin{lemma}\label{matrix dev lem}
	Let $D$ be a non-empty domain with a nontrivial interior in $\mathbb R^d$, and {let $U$ be a compact set}. Let 	
	$M$ and $K$ be matrix-valued mappings $D\times U\mapsto \mathbb R^{d\times d}$ such that $M,K\in G^{\delta}(U,L^{\infty}(D))$ for some $\delta\geq 1$. Specially, for all $\boldsymbol{y}\in U$ and  for all multi-indices $\boldsymbol{\nu}\in \mathcal F$, there {exist} constants $C_M\geq 1$, $C_K>0$  and a sequence of real numbers $\boldsymbol{b} \in \mathbb R^s_+$ such that the following inequalities hold:
	\begin{align} \label{M property}
		\left\|{\partial^{\boldsymbol{\nu}} M}\right\|_{L^{\infty}(D)}
		\leq
		C_M \boldsymbol{b}^{\boldsymbol{\nu}} \left[\tfrac{1}{2} \right]_{\left|\boldsymbol{\nu}\right|}  (\left|\boldsymbol{\nu}\right|!)^{\delta-1},
			\end{align}
	and 
		\begin{align} \label{K property}
		\left\|{\partial^{\boldsymbol{\nu}} K}\right\|_{L^{\infty}(D)}
		\leq
		C_K \boldsymbol{b}^{\boldsymbol{\nu}} \left[\tfrac{1}{2} \right]_{\left|\boldsymbol{\nu}\right|}  (\left|\boldsymbol{\nu}\right|!)^{\delta-1}.
	\end{align}
	Additionally, we assume that $M$ has a uniformly bounded inverse for all $\boldsymbol{x}\in D$ and $\boldsymbol{y}\in U$, with the bound $\left\|{ M^{-1}}\right\|_{L^{\infty}(D)}\leq C_M$. Under these assumptions, the following estimates hold:
	\begin{enumerate}[(i)]
		\item \textbf{Product rule:} $	\left\|{\partial^{\boldsymbol{\nu}} (M\,K)}\right\|_{L^{\infty}(D)}	\leq
		4 C_M\, C_K \boldsymbol{b}^{\boldsymbol{\nu}} \left[\tfrac{1}{2} \right]_{\left|\boldsymbol{\nu}\right|}   (\left|\boldsymbol{\nu}\right|!)^{\delta-1}$.
		\item \textbf{Inverse Matrix Estimate: } $	\left\|{\partial^{\boldsymbol{\nu}} M^{-1}}\right\|_{L^{\infty}(D)}	\leq C_M^3
		(3{C_M^2})^{\left|\boldsymbol{\nu}\right|-1}\boldsymbol{b}^{\boldsymbol{\nu}} \left[\tfrac{1}{2} \right]_{\left|\boldsymbol{\nu}\right|}   (\left|\boldsymbol{\nu}\right|!)^{\delta-1}$ for all multi-indices $\boldsymbol{\nu}\neq \boldsymbol{0}$.
		\item \textbf{Determinant Estimate: }$	\left\|{\partial^{\boldsymbol{\nu}} \det M}\right\|_{L^{\infty}(D)}	\leq
		d!\, (4C_M)^d  \boldsymbol{b}^{\boldsymbol{\nu}} \left[\tfrac{1}{2} \right]_{\left|\boldsymbol{\nu}\right|}   (\left|\boldsymbol{\nu}\right|!)^{\delta-1}$.
	\end{enumerate}
\end{lemma}
\begin{proof}
	For the first statement, applying the Leibniz product rule, {the triangle inequality and \eqref{Frobenius properties}}, for every fixed $\boldsymbol{x}\in D$  we obtain 
	\begin{align*}
		\left|\partial^{\boldsymbol{\nu}}(M\,K)\right|_{F}
		=
		\left|
			\sum_{\boldsymbol{0} \leq  \boldsymbol{\eta}\leq \boldsymbol{\nu}}
			\left(\boldsymbol{\nu} \atop  \boldsymbol{\eta}\right) \partial^{\boldsymbol{\nu}- \boldsymbol{\eta}}M\partial^{ \boldsymbol{\eta}}K\right|_{F}\leq
		\sum_{\boldsymbol{0} \leq  \boldsymbol{\eta}\leq \boldsymbol{\nu}}
		\left(\boldsymbol{\nu} \atop  \boldsymbol{\eta}\right) 
		\left|\partial^{\boldsymbol{\nu}- \boldsymbol{\eta}}M	\right|_{F}	\left|\partial^{ \boldsymbol{\eta}}K\right|_{F}.
	\end{align*}
	Taking $L^{\infty}(D)$-norm both sides of the equation and using the giving bounds for $M$ and $K$, we obtain
	\begin{align*}
		\left\|{\partial^{\boldsymbol{\nu}}(M\,K)}\right\|_{L^{\infty}(D)}
		&\leq
		\sum_{\boldsymbol{0} \leq  \boldsymbol{\eta}\leq \boldsymbol{\nu}}
		\left(\boldsymbol{\nu} \atop  \boldsymbol{\eta}\right) 
		C_M \boldsymbol{b}^{\boldsymbol{\nu}- \boldsymbol{\eta}} \left[\tfrac{1}{2} \right]_{\left|\boldsymbol{\nu}- \boldsymbol{\eta}\right|}  (\left|\boldsymbol{\nu}- \boldsymbol{\eta}\right|!)^{\delta-1}
		C_K \boldsymbol{b}^{\boldsymbol{\eta}} \left[\tfrac{1}{2} \right]_{\left|\boldsymbol{\eta}\right|}  (\left|\boldsymbol{\eta}\right|!)^{\delta-1}\\
		&\leq
		C_M\,C_K  \boldsymbol{b}^{\boldsymbol{\nu}} (\left|\boldsymbol{\nu}\right|!)^{\delta-1}
		\sum_{\boldsymbol{0} \leq  \boldsymbol{\eta}\leq \boldsymbol{\nu}}
		\left(\boldsymbol{\nu} \atop  \boldsymbol{\eta}\right) 
		\left[\tfrac{1}{2} \right]_{\left|\boldsymbol{\nu}- \boldsymbol{\eta}\right|}  
		\left[\tfrac{1}{2} \right]_{\left|\boldsymbol{\eta}\right|} , 
\end{align*}
	where in the last step, we factorize terms and apply factorial inequalities \eqref{multiindex-est-1}. Employing \eqref{multiindex-est-2}, we arrive at
	\begin{align*}
		\left\|{\partial^{\boldsymbol{\nu}}(M K)}\right\|_{L^{\infty}(D)} \leq 4 C_M C_K  \boldsymbol{b}^{\boldsymbol{\nu}} \left[\tfrac{1}{2} \right]_{\left|\boldsymbol{\nu}\right|}  (\left|\boldsymbol{\nu}\right|!)^{\delta-1},
	\end{align*}
	which completes the proof of the first statement.
	
	For the second statement, we proceed by induction. For the base case $\boldsymbol{\nu}=\boldsymbol{e}$ and $\left|\boldsymbol{e}\right|=1$, by differentiating both sides of $MM^{-1}=\text{Id}$ with respect to $\boldsymbol{y}$, we obtain
	\begin{align*}
		{\partial^{\boldsymbol{e}} M^{-1}}(\boldsymbol{x})
		=
		{-M^{-1}(\boldsymbol{x})\, \partial^{\boldsymbol{e}}  M(\boldsymbol{x})\,  M^{-1}}(\boldsymbol{x}).
	\end{align*}
	Taking the $L^{\infty}(D)$-norm and using given bounds for $M$ and $M^{-1}$, we arrive
	\begin{align*}
		\left\|{\partial^{\boldsymbol{e}} M^{-1}}\right\|_{L^{\infty}(D)} 
		=
		\left\|{M^{-1}\, \partial^{\boldsymbol{e}}  M\,  M^{-1}}\right\|_{L^{\infty}(D)}
		\leq
		\left\|{M^{-1}}\right\|_{L^{\infty}(D)}^2\left\|{ \partial^{\boldsymbol{e}}  M}\right\|_{L^{\infty}(D)} \leq
		C_M^3 \boldsymbol{b}^{\boldsymbol{e}} \left[\tfrac{1}{2} \right]_{\boldsymbol{e}} .
	\end{align*}
	For the inductive step, we assume the statement holds for all multi-indices $\boldsymbol{\eta}$ such that $\boldsymbol{\eta}<\boldsymbol{\nu}$. We now aim to prove it for  $\boldsymbol{\nu}$ with $\left|\boldsymbol{\nu}\right|\geq 2$. Applying the triangle inequality and using \eqref{multiindex-est-3}, we have
	\begin{align*}
		\left|\partial^{\boldsymbol{\nu}} M^{-1}\right|_{F}
		&=
		\left| M^{-1}(\boldsymbol{x})
			\sum_{\boldsymbol{0} \leq \boldsymbol{\eta} < \boldsymbol{\nu}}
			\left(\boldsymbol{\nu} \atop \boldsymbol{\eta}\right)
			\partial^{\boldsymbol{\nu}-\boldsymbol{\eta}}M\,  \partial^{\boldsymbol{\eta}} M^{-1}\right|_{F}\\
		&\leq
		\left|M^{-1}\right|_{F}	 
		\left(
			\left|\partial^{\boldsymbol{\nu}}M\right|_{F}	 
			\left|M^{-1}\right|_{F}	 
			+
			\sum_{\boldsymbol{0} < \boldsymbol{\eta} < \boldsymbol{\nu}}
			\left(\boldsymbol{\nu} \atop \boldsymbol{\eta}\right)
			\left|\partial^{\boldsymbol{\nu}-\boldsymbol{\eta}}M\right|_{F}	 
			\left| \partial^{\boldsymbol{\eta}} M^{-1}\right|_{F}	\right) .
	\end{align*}
	Taking the $L^{\infty}(D)$-norm of both sides and applying inductive assumption, we obtain
	\begin{align*}		  
		&\left\|{\partial^{\boldsymbol{\nu}} M^{-1}}\right\|_{L^{\infty}(D)}
		\leq
		C_M
		\Biggl(
		C_M \boldsymbol{b}^{\boldsymbol{\nu}} \left[\tfrac{1}{2} \right]_{\left|\boldsymbol{\nu}\right|}  (\left|\boldsymbol{\nu}\right|!)^{\delta-1}\, C_M
		\\
		&\quad\quad\quad\quad
		\left.
		+
		\sum_{\boldsymbol{0} < \boldsymbol{\eta} < \boldsymbol{\nu}}
		\left(\boldsymbol{\nu} \atop \boldsymbol{\eta}\right)
		C_M \boldsymbol{b}^{\boldsymbol{\nu}-\boldsymbol{\eta}} \left[\tfrac{1}{2} \right]_{\left|\boldsymbol{\nu}-\boldsymbol{\eta}\right|}  (\left|\boldsymbol{\nu}-\boldsymbol{\eta}\right|!)^{\delta-1}\,
		C_M^3 (3{C_M^2})^{\left|\boldsymbol{\eta}\right|-1} \boldsymbol{b}^{\boldsymbol{\eta}} \left[\tfrac{1}{2} \right]_{\left|\boldsymbol{\eta}\right|}  (\left|\boldsymbol{\eta}\right|!)^{\delta-1}
		\right)\\
		&\quad\quad\quad\leq
		C_M^3\,  {C_M^2} (3{C_M^2})^{\left|\boldsymbol{\nu}\right|-2} \boldsymbol{b}^{\boldsymbol{\nu}} (\left|\boldsymbol{\nu}\right|!)^{\delta-1}
		\sum_{\boldsymbol{0} \leq \boldsymbol{\eta} < \boldsymbol{\nu}}
		\left(\boldsymbol{\nu} \atop \boldsymbol{\eta}\right)
		\left[\tfrac{1}{2} \right]_{\left|\boldsymbol{\nu}-\boldsymbol{\eta}\right|} \,
		\left[\tfrac{1}{2} \right]_{\left|\boldsymbol{\eta}\right|} \\
		&\quad\quad\quad\leq
		C_M^3(3{C_M^2})^{\left|\boldsymbol{\nu}\right|-1}\boldsymbol{b}^{\boldsymbol{\nu}}   \left[\tfrac{1}{2} \right]_{\left|\boldsymbol{\nu}\right|} (\left|\boldsymbol{\nu}\right|!)^{\delta-1}.
	\end{align*}
	In the last step, we use equations \eqref{multiindex-est-1}  and \eqref{multiindex-est-3} to conclude the proof of the second statement.
	
	To prove the third statement, we proceed by induction on the dimension $d$. For the base case $d=1$, we trivially have $\det M=M$, which directly satisfies the desired bound.
	For the inductive step, we assume the statement holds for all matrices of dimension $k$ with $k<d$. We aim to establish it for $d$.
	By applying the Laplace expansion along any column $j$, we obtain
	\begin{align}\label{Laplace iden}
		\det M (\boldsymbol{x})= \sum_{i=1}^{d} (-1)^{i+j} M_{ij}(\boldsymbol{x}) \, \det M_{ij}^{\#}(\boldsymbol{x}),
	\end{align}
	here, the submatrix $M_{ij}^{\#}$ is obtained by deleting the $i$-th row and the $j$-th column of matrix $M$, hence \eqref{M property} implies
\begin{align}\label{matrix M property}
		\left|\partial^{\boldsymbol{\nu}} M_{ij}^{\#}\right|_{F}
		\leq 
		\left|\partial^{\boldsymbol{\nu}} M\right|_{F}
		\leq
		C_M \boldsymbol{b}^{\boldsymbol{\nu}}  \left[\tfrac{1}{2} \right]_{\left|\boldsymbol{\nu}\right|} (\left|\boldsymbol{\nu}\right|!)^{\delta-1}
		\quad\forall \boldsymbol{x}\in D.
	\end{align}
	We apply the inductive hypothesis for the $(d-1)\times(d-1)$-dimensional matrix $M_{ij}^{\#}$ to obtain
	\begin{align}\label{det M hypo}
		\left|\partial^{\boldsymbol{\nu}} \det M_{ij}^{\#}\right|
		\leq
		(d-1)!\,(4C_M)^{d-1}  \boldsymbol{b}^{\boldsymbol{\nu}}  \left[\tfrac{1}{2} \right]_{\left|\boldsymbol{\nu}\right|} (\left|\boldsymbol{\nu}\right|!)^{\delta-1}
		\quad\forall \boldsymbol{x}\in D.
	\end{align}
	Taking the $\boldsymbol{\nu}$-th derivative both sides of \eqref{Laplace iden} and applying Leibniz product rule yields
	\begin{align*}
		{\partial^{\boldsymbol{\nu}}	\det M(\boldsymbol{x})}
		&=
		{ \sum_{i=1}^{d}
			(-1)^{i+j} 
			\sum_{\boldsymbol{0}\leq \boldsymbol{\eta}\leq \boldsymbol{\nu}}
			\left(\boldsymbol{\nu}\atop \boldsymbol{\eta}\right)
			\partial^{\boldsymbol{\nu}-\boldsymbol{\eta}} M_{ij}(\boldsymbol{x}) \,
			\partial^{\boldsymbol{\eta}} \det M_{ij}^{\#}(\boldsymbol{x})}.
\end{align*}
	Taking absolute value both sides and using triangle inequality we have:
	\begin{align*}
		\left|\partial^{\boldsymbol{\nu}}	\det M\right|
&\leq
		\sum_{i=1}^{d}
		\sum_{\boldsymbol{0}\leq \boldsymbol{\eta}\leq \boldsymbol{\nu}}
		\left(\boldsymbol{\nu}\atop \boldsymbol{\eta}\right)
		\left|\partial^{\boldsymbol{\nu}-\boldsymbol{\eta}} M_{ij}\right| \,
		\left|	\partial^{\boldsymbol{\eta}} \det M_{ij}^{\#}\right|.
	\end{align*}
	Taking $L^{\infty}(D)$-norm both sides and noting \eqref{det M hypo} and   \eqref{multiindex-est-1}, we obtain
	\begin{align*}
		&\left\|{\partial^{\boldsymbol{\nu}}	\det M}\right\|_{L^{\infty}(D)}\\
		&\quad\quad\leq
		\sum_{i=1}^{d}
		\sum_{\boldsymbol{0}\leq \boldsymbol{\eta}\leq \boldsymbol{\nu}}
		\left(\boldsymbol{\nu}\atop \boldsymbol{\eta}\right)
		C_M \boldsymbol{b}^{\boldsymbol{\nu}-\boldsymbol{\eta}} \left[\tfrac{1}{2} \right]_{\left|\boldsymbol{\nu}-\boldsymbol{\eta}\right|}  (\left|\boldsymbol{\nu}-\boldsymbol{\eta}\right|!)^{\delta-1}\,
		(d-1)! (4C_M)^{d-1}  \boldsymbol{b}^{\boldsymbol{\eta}} \left[\tfrac{1}{2} \right]_{\left|\boldsymbol{\eta}\right|}  (\left|\boldsymbol{\eta}\right|!)^{\delta-1}\\
		&\quad\quad\leq
		d!C_M^{d} 4^{d-1} \boldsymbol{b}^{\boldsymbol{\nu}} (\left|\boldsymbol{\nu}\right|!)^{\delta-1}
		\sum_{\boldsymbol{0}\leq \boldsymbol{\eta}\leq \boldsymbol{\nu}}
		\left(\boldsymbol{\nu}\atop \boldsymbol{\eta}\right)
		\left[\tfrac{1}{2} \right]_{\left|\boldsymbol{\nu}-\boldsymbol{\eta}\right|} \,
		\left[\tfrac{1}{2} \right]_{\left|\boldsymbol{\eta}\right|}  .
	\end{align*}
	Finally, applying estimates \eqref{multiindex-est-2}, we conclude
	\begin{align*}
		\left\|{\partial^{\boldsymbol{\nu}}	\det M}\right\|_{L^{\infty}(D)}
		&\leq
		d!(4C_M)^{d}  \boldsymbol{b}^{\boldsymbol{\nu}} \left[\tfrac{1}{2} \right]_{\left|\boldsymbol{\nu}\right|} (\left|\boldsymbol{\nu}\right|!)^{\delta-1}.
	\end{align*}
	This completes the proof of the third statement.
\end{proof}
\begin{remark}
Observe that the Gevrey parameter $\delta$ only appears in the exponent of the factorial term in (i)--(iii) in Lemma \ref{matrix dev lem} and has no effect on the other constants.
\end{remark} 
\subsection{Multivariate Fa{\`a} di Bruno formula}\label{sec: faa di bruno}
In the context of domain transform, the given force term in the equation on the reference domain can be expressed as composite functions involving pullback operations. To establish Gevrey class regularity for composite functions, we will apply the multivariate Fa{\`a} di Bruno’s formula, as presented by Constantine and Savits \cite[Theorem 2.1]{CONSTANTINE96}.

For any ${\boldsymbol{\nu}}\in \mathcal F\setminus \left\{\boldsymbol{0}\right\}$, we assume that the vector field $\boldsymbol{g}$ is sufficiently smooth around an open neighborhood of $\boldsymbol{z}$ and that the function $f$ is sufficiently smooth around an open neighborhood $\boldsymbol{g}(\boldsymbol{z})$. Under these conditions, $\partial^{\boldsymbol{\nu}} f(\boldsymbol{g}(\boldsymbol{z}))$  exists and can be explicitly expressed as follows
\begin{align}\label{Faa di Bruno}
	\partial^{\boldsymbol{\nu}} f(\boldsymbol{g}(\boldsymbol{z}))
	=
	\sum_{1\leq\left|\boldsymbol{\lambda}\right|\leq \left|\boldsymbol{\nu}\right|}
	\partial^{\boldsymbol{\lambda}}f
	\sum_{s=1}^n
	\sum_{P_s(\boldsymbol{\nu},\boldsymbol{\lambda})}
	\boldsymbol{\nu}!
	\prod_{j=1}^s
	\frac{(\partial^{\boldsymbol{\ell}_j}\boldsymbol{g}(\boldsymbol{z}))^{\boldsymbol{k}_j}}{\boldsymbol{k}_j! (\boldsymbol{\ell}_j!)^{\left|\boldsymbol{k}_j\right|}},
\end{align}
where 
\begin{align*}
	P_s(\boldsymbol{\nu},\boldsymbol{\lambda})
	&=
	\Biggl\{\left(
		\left\{\boldsymbol{k}_j\right\}_{j=1}^s;
		\left\{\boldsymbol{\ell}_j\right\}_{j=1}^s
	\right)
	: \left|\boldsymbol{k}_j\right|>0,\, \sum_{j=1}^s \boldsymbol{k}_i=\boldsymbol{\lambda}
	\\
	&\quad\quad\quad \boldsymbol{0} \prec \boldsymbol{\ell}_1\prec \dots \prec \boldsymbol{\ell}_s \text{ and } \sum_{j=1}^s \left|\boldsymbol{k}_i\right|\boldsymbol{\ell}_i=\boldsymbol{\nu}
	\Biggl\} .
\end{align*}
To enhance our understanding of the index set $P_s(\boldsymbol{\nu},\boldsymbol{\lambda})$, we present the following lemma. 
\begin{lemma}\label{lem: n k!}
	Let $\boldsymbol{\lambda}$ and $\boldsymbol{\nu}$ be two multi-indices in $\mathcal F\setminus \left\{\boldsymbol{0}\right\}$ such that $\left|\boldsymbol{\lambda}\right|\leq \left|\boldsymbol{\nu}\right|$. Suppose that
	$$\left(
		\left\{\boldsymbol{k}_j\right\}_{j=1}^s;
		\left\{\boldsymbol{\ell}_j\right\}_{j=1}^s
	\right)\in P_s(\boldsymbol{\nu},\boldsymbol{\lambda}),$$
	then for any $\omega\geq 0$, the following inequality holds:
	\begin{align*}
		(\left|\boldsymbol{\lambda}\right|! )^{\omega}
		\prod_{j=1}^s
		(\left|\boldsymbol{\ell}_j\right|!)^{\omega\left|\boldsymbol{k}_j\right|}
		\leq
		(\left|
			\boldsymbol{\nu}\right|!)^{\omega}.
	\end{align*}
\end{lemma}
\begin{proof}
	To proceed, we make use of the following inequalities:
	\begin{align}\label{n! k!}
		n!\,k!\leq (n+k-1)!  \quad\forall n,k\geq 1.
	\end{align}
	Indeed, for every $k,n\geq 1$, we observe that
	\begin{align*}
		n!\,k!
		=
		n!\left(2\cdot 3\dots  k\right)
		\leq
		n! (n+1)\,(n+2)\cdots (n+k-1)
		=
		(n+k-1)!.
	\end{align*}
	Applying \eqref{n! k!} inductively, we obtain the estimate
	\begin{align*}
		\prod_{j=1}^s m_j! \leq \left({\sum_{j=1}^s} m_j - (s-1)\right)!,  \quad\forall \boldsymbol{m}\geq \boldsymbol{1}.
	\end{align*}
	We now apply this bound to show that for any pair of multi-indices $\boldsymbol{\ell}$ and $\boldsymbol{k}$ satisfying $(\boldsymbol{\ell}, \boldsymbol{k}) \in P_s({\boldsymbol{\nu},\boldsymbol{\lambda}})$, the following holds:
	\begin{align*}
		\left|\boldsymbol{\lambda}\right|! \prod_{j=1}^s
		\left|\boldsymbol{\ell}_j\right|!^{\left|\boldsymbol{k}_j\right|}
		&\leq
		\left|\boldsymbol{\lambda}\right|!
		\prod_{j=1}^s
		(\left|\boldsymbol{k}_j\right|\left|\boldsymbol{\ell}_j\right|-(\left|\boldsymbol{k}_j\right|-1))!\\
		&\leq
		\left|\boldsymbol{\lambda}\right|!
		\left(\sum_{j=1}^s \left(\left|\boldsymbol{k}_j\right|\left|\boldsymbol{\ell}_j\right|- (\left|\boldsymbol{k}_j\right|-1)\right) -(s-1)\right)	  !\\
		&=
		\left|\boldsymbol{\lambda}\right|!
		\left(\left|\boldsymbol{\nu}\right|-\left|\boldsymbol{\lambda}\right|+1\right)!
		\leq
		\left|\boldsymbol{\nu}\right|!.
	\end{align*}
	In the final step, we utilize the properties of $P_s({\boldsymbol{\nu},\boldsymbol{\lambda}})$, namely 
	$$\sum_{j=1}^s \left|\boldsymbol{k}_j\right| \left|\boldsymbol{\ell}_j\right|=\left|\boldsymbol{\nu}\right|\quad \mathrm{and}  \quad\sum_{j=1}^s \left|\boldsymbol{k}_j\right| =\left|\boldsymbol{\lambda}\right|.$$ 
	Since the function $x \mapsto x^{\omega}$ is non-decreasing for all $x\geq 1$ and $\omega\geq 0$, it follows that
	$$
	(\left|\boldsymbol{\lambda}\right|! )^{\omega}
	\prod_{j=1}^s
	(\left|\boldsymbol{\ell}_j\right|!)^{\omega\left|\boldsymbol{k}_j\right|}
	=
	\left(\left|\boldsymbol{\lambda}\right|! \prod_{j=1}^s
		(\left|\boldsymbol{\ell}_j\right|!)^{{\left|\boldsymbol{k}_{j}\right|}}\right)^{\omega}
	\leq
	(\left|\boldsymbol{\nu}\right|!)^{\omega}.
	$$
	This completes the proof.
\end{proof}

The following lemma plays a crucial role in the parametric analysis of domain transformation. While this result is not new and was established in {\cite[p. 28-29]{Costabel10}, cf. \cite[Lemma 1.4.1]{Krantz2002}}, it was not explicitly stated as a separate lemma. Therefore, we present it here using the notation as in \eqref{Faa di Bruno} {and give a proof for completeness}.
\begin{lemma}\label{singlevariable FdB lemma}
	For any multi-index $\boldsymbol{\nu}\in \mathcal F\setminus \left\{\boldsymbol{0}\right\}$ and for any integer $d\in\mathbb N$, the following holds
	\begin{align}
		\sum_{1\leq\left|\boldsymbol{\lambda}\right|\leq \left|\boldsymbol{\nu}\right|
			\atop \boldsymbol{\lambda}\in \mathbb N^d}
		\left|\boldsymbol{\lambda}\right|!
		\sum_{s=1}^n
		\sum_{P_s(\boldsymbol{\nu},\boldsymbol{\lambda})}
		\boldsymbol{\nu}!
		\prod_{j=1}^s
		\frac{(\left|\boldsymbol{\ell}_{j}\right|!)^{\left|\boldsymbol{k}_{j}\right|}}{\boldsymbol{k}_{j}! (\boldsymbol{\ell}_{j}!)^{\left|\boldsymbol{k}_j\right|}}
		=
		d\,(d+1)^{\left|\boldsymbol{\nu}\right|-1} \left|\boldsymbol{\nu}\right|!\label{singlevariable FdB form}.
	\end{align}
\end{lemma}
\begin{proof}
	{Following the proof of equation (1.24) in \cite{Costabel10} the statement follows from the Fa{\`a} di Bruno formula applied for specific $f(\boldsymbol{w})$ and $\boldsymbol{g}(\boldsymbol{z})
		=	(g_1(\boldsymbol{z}),\dots,g_d(\boldsymbol{z}))$  defined as follows 
		\begin{align*}
			f(\boldsymbol{w}):=\frac{1}{1-\boldsymbol{w}\cdot \boldsymbol{1}},\qquad 
			g_i(\boldsymbol{z})
			:= f(\boldsymbol{z}) -1 = 
			\frac{1}{1-\boldsymbol{z}\cdot \boldsymbol{1}}-1
		\end{align*}
		and their composition
		\begin{align*}
			h(\boldsymbol{z})&:=f(\boldsymbol{g}(\boldsymbol{z}))
			=
			\frac{1}{1-\boldsymbol{1}\cdot\boldsymbol{g}(\boldsymbol{z})}
			=
			\frac{1}{1-d(\frac{1}{1-\boldsymbol{z}\cdot\boldsymbol{1}}-1)}
			=
			\frac{d(d+1)^{-1}}{1-(d+1)\boldsymbol{z}\cdot\boldsymbol{1}}
			+
			\frac{1}{d+1}.
		\end{align*}	
		The chain rule yields for all $\boldsymbol{\ell},\boldsymbol{\lambda}, \boldsymbol{\nu} \in \mathcal F\setminus \left\{\boldsymbol{0}\right\}$. 
		\[
		\partial^{\boldsymbol{\lambda}} f(\boldsymbol{0})=\left|\boldsymbol{\lambda}\right|!, \qquad \partial^{\boldsymbol{\ell}} g_i(\boldsymbol{0})= \left|\boldsymbol{\ell}\right|!, \qquad 
		\partial^{\boldsymbol{\nu}} h(\boldsymbol{0})=d\,(d+1)^{\left|\boldsymbol{\nu}\right|-1} \left|\boldsymbol{\nu}\right|!.
		\]
		The statement follows now by the Fa{\`a} di Bruno formula and the observation
		\begin{align*}\prod_{i=1}^d(\partial^{\boldsymbol{\ell}}g_i(\boldsymbol{0}))^{k_{i}}=
			\prod_{i=1}^d(\left|\boldsymbol{\ell}\right|!)^{k_{i}}
			=
			(\left|\boldsymbol{\ell}\right|!)^{\left|\boldsymbol{k}\right|}, \qquad
			\forall \boldsymbol{\ell}\in\mathcal F\setminus\left\{\boldsymbol{0}\right\}, ~~\forall \boldsymbol{k}\in \mathbb N^d.
	\end{align*}
}\end{proof}
\section{Well-posedness of the parametric Navier-Stokes equation}\label{sec: well-posedness}

\subsection{Mixed and saddle point problems} \label{sec: NS saddle point}
In this section, we aim to concisely summarize the key properties of the Navier-Stokes equations in the form \eqref{NS gen equation}, defined within a bounded Lipschitz domain $D\in \mathbb R^{d}$ for $d\in \left\{1,2,3\right\}$, as follows
\begin{align*}
		a({\boldsymbol{u}},{\boldsymbol{v}};A(\boldsymbol{y})) &:= \int_{{D}} \text{tr} \big(\nabla {\boldsymbol{u}}^{\top} A(\boldsymbol{y})^{\top} \nabla{\boldsymbol{v}} \big),\\
		b({\boldsymbol{v}},{p};B(\boldsymbol{y})) &:= -\int_{{D}} \text{tr} \big(\nabla{\boldsymbol{v}}  B(\boldsymbol{y}){p} \big),\\
		m({\boldsymbol{u}},{\boldsymbol{v}},{\boldsymbol{w}};M(\boldsymbol{y})) &:= \int_{{D}} (\nabla {\boldsymbol{v}}\, M(\boldsymbol{y}) {\boldsymbol{u}} )^{\top}  {\boldsymbol{w}},\\
		\left\langle {\boldsymbol{f}}(\boldsymbol{y}), \boldsymbol{v} \right\rangle &:= \int_{{D}} {\boldsymbol{f}} \cdot {\boldsymbol{v}} ,\\
		\llangle { {g}(\boldsymbol{y})},{ q}
	\rrangle &:= \int_{{D}} {g} {q} .
\end{align*}
These properties below have been extensively studied in previous works, e.g. \cite{Temam1977, Girault1986, Cohen2018}. In this paper, we provide an independent demonstration {aiming}
to simplify and clarify \ac{the} fundamental properties through an alternative derivation, enhancing the understanding of the Navier-Stokes equations with non-zero source mass term $g$.

 We define the supremum norm as follows:
\begin{align}\label{norm definition}
\left\|{A(\boldsymbol{x},\boldsymbol{y})}\right\|_{\infty}:=
\sup_{\boldsymbol{y}\in U} 
\left\|{\, \left|A(\boldsymbol{x},\boldsymbol{y})\right|_{F} }\right\|_{L^{\infty}(D)}
=
\sup_{\boldsymbol{y}\in U} {\esssup_{\boldsymbol{x}\in D}} \left|A(\boldsymbol{x},\boldsymbol{y})\right|_{F}.
\end{align}
With the definition above, we state an assumption on the matrices $A$, $B$, and $M$ as follows:
\begin{assumption}\label{A,B,M assump}
	For any $(\boldsymbol{x},\boldsymbol{y}) \in D\times U$,  let $A(\boldsymbol{x},\boldsymbol{y}), B(\boldsymbol{x},\boldsymbol{y})$ and $M(\boldsymbol{x},\boldsymbol{y})$ be non-singular matrices and $A(\boldsymbol{x},\boldsymbol{y})$ be symmetric such that $A,B,M\in L^{\infty}(D)^{d\times d}$. Moreover, there exist three positive constants $\overline{a}, \overline{b}$ and $\overline{m}$ such that
\begin{align}\label{NS ab-bounds}
\left\|{A}\right\|_{\infty}\leq \frac{\overline{a}}{2},\quad
\left\|{B}\right\|_{\infty}\leq \frac{\overline{b}}{2},\quad\text{ and }\quad
\left\|{M}\right\|_{\infty}\leq C_4^{-2}\, \frac{\overline{m}}{2},
\end{align}
here, the constant $C_4$  represents the constant associated with the
Sobolev embedding $H^1_0(D)\hookrightarrow L^4(D)$, {i.e. $\|u\|_{L^4(D)} \leq C_4 \|u\|_{H^1_0(D)}$}.  For numerical values of $C_4$, {we} refer, for instance,
to \cite{Mizuguchi2017}.
\end{assumption}
We define the following Banach spaces:
\begin{align*}
	{\mathcal H}:=H^1_0(D)^d,\quad {\mathcal H}^*:=H^{-1}(D)^d,\quad {\mathcal K}:=L^4(D)^d,\quad
	\mathcal L:=L^2({D})\setminus \mathbb R.
\end{align*}
Each space is equipped with the corresponding norm:
\begin{align}
	\left\|{\boldsymbol{u}}\right\|_{{\mathcal H}}
	&:=
	\left(
	\int_{D} \left|\nabla \boldsymbol{u}\right|_F^2\right)^{\frac{1}{2}}
	= \left(
	\sum_{i=1}^d \sum_{j=1}^d \left\|{\frac{\partial u_i}{\partial x_j}}\right\|_{L^2(D)}^2
	\right)^{\frac{1}{2}},\label{H norm def}
	\\
	\left\|{\boldsymbol{f}}\right\|_{{\mathcal H}^*}
	&:=
	\sup_{\boldsymbol{v}\in {\mathcal H}}
	\frac{\left\langle{\boldsymbol{f}},{\boldsymbol{v}}
	\right\rangle}{\left\|{\boldsymbol{v}}\right\|_{{\mathcal H}}},\notag\\
	\left\|{\boldsymbol{u}}\right\|_{{\mathcal K}}&:
	=
	C_4^{-1}\left(\int_D \left|\boldsymbol{u}\right|_F^4 \right)^{\frac{1}{4}}
	=
	C_4^{-1}\left(\int_D \left(\sum_{i=1}^d u_i^2\right)^2 \right)^{\frac{1}{4}},\notag\\
\left\|{p}\right\|_{\mathcal L}
	&:= \min_{c\in \mathbb R} \left\|{p-c}\right\|_{L^2(D)}.\notag
\end{align}
Building upon those definitions, we delve into the key properties of $\mathcal L$, ${\mathcal H}$, and ${\mathcal K}$.
\begin{lemma} \label{lem: L H K prop}
Recall the definition of $\mathcal L$, ${\mathcal H}$, ${\mathcal K}$ and their norms, there holds
\begin{itemize}
\item[(i)] $\mathcal L$ is canonically isomorphic to its dual, i.e. there exists an isomorphism mapping from $\mathcal L$ to its dual. Moreover, for every $q\in \mathcal L$, there holds $\left\|{q}\right\|_{\mathcal L}=\left\|{q}\right\|_{L^2(D)}$.
\item[(ii)] ${\mathcal H} \hookrightarrow{\mathcal K}$, i.e. $\left\|{\boldsymbol{w}}\right\|_{{\mathcal K}}\leq\left\|{\boldsymbol{w}}\right\|_{{\mathcal H}}$ for every $\boldsymbol{w}\in {\mathcal H}$.
\end{itemize} 
\end{lemma}
\begin{proof}
To prove the statement $(i)$, recall the definition of $\left\|{\cdot}\right\|_{\mathcal L}$. For every $q \in \mathcal L$, we have:
\begin{align}\label{L norm prop}
\left\|{q}\right\|_{\mathcal L}^2
=
\min_{c\in \mathbb R} \int_D (q-c)^2
=
\int_D q^2
+
\min_{c\in \mathbb R} 
\left(
\left|D\right| c^2 - 2 \int_D q \, c\right)
=
\int_D q^2
-
\frac{1}{\left|D\right|}\left(\int_D q\right)^2,
\end{align}
where in the last step, we notice that the minimum {is attained} at $c=\left|D\right|^{-1}\int_D q$. Due to $\int_D q =0$ for every $q\in \mathcal L$, equation \eqref{L norm prop} shows $\left\|{q}\right\|_{\mathcal L}=\left\|{q}\right\|_{L^2(D)}$. Based on the definition of the $\mathcal L$-norm, we can define a $\mathcal L$-inner product for all $p,q\in \mathcal L$ as follows
\begin{align}
	\left\langle{p},{q}
	\right\rangle_{\mathcal L}
	:=
	\int_D p\,q - \frac{1}{\left|D\right|} \int_D p\, \int_D q.
\end{align}
Utilizing the Riesz Representation Theorem, we further conclude that the Hilbert space $\mathcal L$ is isometrically isomorphic to its dual, thereby verifying statement $(i)$.

For statement $(ii)$, we now apply H\"older inequality to obtain
\begin{align*}
		\left\|{\boldsymbol{w}}\right\|_{{\mathcal K}}^4
&=
{\frac{1}{C_4^4}}
\int_D \left(\sum_{i=1}^d w_i^2\right)^2
=
{\frac{1}{C_4^4}}
\sum_{i=1}^d \sum_{j=1}^d 
\int_D w_i^2 w_j^2\\
&\leq
{\frac{1}{C_4^4}}
\sum_{i=1}^d \sum_{{j=1}}^d 
\left\|{w_i}\right\|_{L^4(D)}^2 \left\|{w_j}\right\|_{L^4(D)}^2
=
{\frac{1}{C_4^4}}
\left(\sum_{i=1}^d \left\|{ w_i}\right\|_{L^4(D)}^2\right)^{2}.
\end{align*}
By Sobolev embedding $H^1_0(D)^d\hookrightarrow L^4(D)^d$, we have
\begin{align*}
	\left\|{\boldsymbol{w}}\right\|_{{\mathcal K}}
	&\leq
	\frac{1}{C_4}\left(\sum_{i=1}^d \left\|{ w_i}\right\|_{L^4(D)}^2\right)^{\frac{1}{2}}
	\leq
	\frac{1}{C_4}\left(\sum_{i=1}^d C_4^2 \left\|{\nabla w_i}\right\|_{L^2(D)}^2\right)^{\frac{1}{2}}
	\\&=
	\left(\sum_{i=1}^d \sum_{j=1}^d \left\|{\frac{\partial w_i}{\partial x_j}}\right\|_{L^2(D)}^2
	\right)^{\frac{1}{2}}
	=
	\left\|{\boldsymbol{w}}\right\|_{{\mathcal H}}.
\end{align*}
It shows the second statement, and hence completes the proof.
\end{proof}
\begin{remark}
It is important to note that $\mathcal L$ is isometrically isomorphic to its dual. Therefore, without loss of generality, in \eqref{NS gen equation}, we can assume  $g\in\mathcal L$ instead of $g\in\mathcal L^*$. Consequently, the duality inner product $\llangle {\cdot},{\cdot}
	\rrangle$ can be interpreted as the $\mathcal L$-inner product when both arguments belong to  $\mathcal L$. Furthermore, the Cauchy-Schwarz inequality applies to the $\mathcal L$-inner product for all $p,q\in \mathcal L$ satisfying $\int_D p = \int_D q=0$ as shown below:
\begin{align*}
	\llangle {g},{q}
	\rrangle
	&=
	\left\langle{g},{q}
	\right\rangle_{\mathcal L}
	=
	\int_D g\, q - \frac{1}{\left|D\right|} \int_D g \int_D q
	=
		\int_D g\, q
		\leq
		\left\|{g}\right\|_{L^2(D)}
		\left\|{q}\right\|_{L^2(D)}
		=
		\left\|{g}\right\|_{\mathcal L}
		\left\|{q}\right\|_{\mathcal L}.
\end{align*}
\end{remark}
The following assumption ensures that the equation \eqref{NS gen equation} is well-posed.
\begin{assumption}\label{alpha beta assump}
There exist two positive constants $\alpha$ and $\beta$ such that for every $\boldsymbol{y} \in U$ there holds
\begin{align}\label{V-coercive}
\alpha \left\|{\boldsymbol{u}}\right\|_{{\mathcal H}}^2 \leq a(\boldsymbol{u},\boldsymbol{u};A(\boldsymbol{y})),
\quad\forall \boldsymbol{u} \in {\mathcal H},
\end{align}		
\begin{align}\label{inf-sup condition}
	\beta \leq {\inf_{0\neq p\in \mathcal L}}\, \sup_{\boldsymbol{0} \neq \boldsymbol{u}\in {\mathcal H}}\, \frac{b(\boldsymbol{u} ,p;B(\boldsymbol{y}))}{\left\|{\boldsymbol{u}}\right\|_{{\mathcal H}}\left\|{p}\right\|_{\mathcal L}}.
\end{align}
\end{assumption}
According to  \cite{Gallistl2018}, we could numerically estimate the inf-sup constant, it depends on the shape of physical domain $D$. Specially, when $B$ is {the} identity matrix  and $D$ is a square domain, we have $\beta=\sqrt{\frac{1}{2}-\frac{1}{\pi}}$.
Due to Assumption \ref{A,B,M assump}, we define the following operators
\begin{equation}
\label{abm def}
\begin{split}
	\left\langle{{\mathcal A}_{\boldsymbol{y}} \boldsymbol{u}},{\boldsymbol{v}}
	\right\rangle
	&:=
	a(\boldsymbol{u},\boldsymbol{v};A(\boldsymbol{y}))
	=
	\left\langle{{\mathcal A}_{\boldsymbol{y}} \boldsymbol{v}},{\boldsymbol{u}}
	\right\rangle,\\
	\llangle {{\mathcal B}_{\boldsymbol{y}}  \boldsymbol{u}},{p}
	\rrangle
	&:=
	b(\boldsymbol{u},p;B(\boldsymbol{y}))
	=:
	\left\langle{{\mathcal B}_{\boldsymbol{y}}^{\top}p},{  \boldsymbol{u}}
	\right\rangle,
\\
\left\langle{\mathcal M_{\boldsymbol{y}} \boldsymbol{u} \cdot \boldsymbol{v}},{\boldsymbol{w}}
	\right\rangle
	&:=
	m(\boldsymbol{u},\boldsymbol{v},\boldsymbol{w};M(\boldsymbol{y})).
\end{split}
\end{equation}
It also implies that
\begin{align*}
	{\mathcal A}_{\boldsymbol{y}}: {\mathcal H} \mapsto {\mathcal H}^*,\quad 
	{\mathcal B}_{\boldsymbol{y}}: {\mathcal H} \mapsto \mathcal L^*,\quad 
	{\mathcal B}_{\boldsymbol{y}}^{\top}: \mathcal L \mapsto {\mathcal H}^*,\quad 
		\mathcal M_{\boldsymbol{y}}: {\mathcal H} \times {\mathcal H} \mapsto {\mathcal H}^*.
\end{align*}
The H\"older inequality implies that the third equation of \eqref{abm def} is well-defined for $\boldsymbol{u}, \boldsymbol{v}\in L^{4}(D)^d$. By the Sobolev embedding theorem, this is guaranteed for $H^1_0(D)^d$ functions due to $H^1_0(D) \hookrightarrow L^{4}(D)$. For all $ \boldsymbol{v},\boldsymbol{u},\boldsymbol{w}\in {\mathcal H}$ and $p\in \mathcal L$, and noting the properties in \eqref{Frobenius properties}, we have
\begin{align}
a\left(\boldsymbol{w},\boldsymbol{v};A(\boldsymbol{y})\right)
&\leq
\int_D
\left|A(\boldsymbol{y})\right|_F \left|\nabla \boldsymbol{w}\right|_F \left|\nabla \boldsymbol{v}\right|_F
\leq
\frac{\overline{a}}{2}
\left\|{\boldsymbol{w}}\right\|_{{\mathcal H}} \left\|{\boldsymbol{v}}\right\|_{{\mathcal H}},\quad 
 \label{a operator bound}\\
b\left(\boldsymbol{w},p;B(\boldsymbol{y})\right)
&=
b\left(\nabla\boldsymbol{w},B(\boldsymbol{y})p\right)
\leq
\int_D
\left|B(\boldsymbol{y})\right|_F \left|\nabla \boldsymbol{w}\right|_F \left|p\right|
\leq
\frac{\overline{b}}{2}
\left\|{\boldsymbol{w}}\right\|_{{\mathcal H}} \left\|{p}\right\|_{\mathcal L},\quad\label{b operator bound}\\ 
m\left(\boldsymbol{w},\boldsymbol{v},\boldsymbol{u};M(\boldsymbol{y})\right)
&\leq
\int_D \left|M(\boldsymbol{y})\right|_F \left|\nabla \boldsymbol{v}\right|_F \left|\boldsymbol{w}\right|_F \left|\boldsymbol{u}\right|_F
\leq
C_4^2\left\|{M}\right\|_{\infty} \left\|{\boldsymbol{w}}\right\|_{{\mathcal K}}
\left\|{\boldsymbol{v}}\right\|_{{\mathcal H}}
\left\|{\boldsymbol{u}}\right\|_{{\mathcal K}}\notag\\
&\leq
\frac{\overline{m}}{2} \left\|{\boldsymbol{w}}\right\|_{{\mathcal H}}
\left\|{\boldsymbol{v}}\right\|_{{\mathcal H}}
\left\|{\boldsymbol{u}}\right\|_{{\mathcal H}}.
\label{m operator bound}
\end{align}
Drawing inspiration from the arguments expounded in \cite[Chap. I Sec. 2.2]{Girault1986}, we establish essential properties in a general context. Notably, the operators ${\mathcal B}_{\boldsymbol{y}}$ and ${\mathcal B}_{\boldsymbol{y}}^{\top}$ take on roles analogous to $\text{div}$ and $\textbf{grad}$ as presented in the referenced book. Moreover, in this framework, the operator ${\mathcal A}_{\boldsymbol{y}}$ assumes the responsibility of $-\Delta$, the classical Laplace operator associated with Dirichlet's homogeneous problem. For the forthcoming well-posedness analysis, we define the following subspaces of ${\mathcal H}$
\begin{align*}
		\mathcal V:=\left\{\boldsymbol{v}\in {\mathcal H}: \llangle {{\mathcal B}_{\boldsymbol{y}} \boldsymbol{v}},{q}
	\rrangle=0, \forall q\in \mathcal L\right\}\quad \text{and}\quad 
	\mathcal V^{\bot}:=\left\{\boldsymbol{v}\in {\mathcal H}: \left\langle{{\mathcal A}_{\boldsymbol{y}} \boldsymbol{v}},{\boldsymbol{w}}
	\right\rangle=0, \forall \boldsymbol{w}\in \mathcal V\right\}.
\end{align*}
\begin{remark}
In the classical context, as seen in works such as \cite{Girault1986,Temam1977}, $\mathcal V$ serves as the space containing functions with zero divergence, i.e., $\text{div}\, u = 0$. Conversely, $\mathcal V^{\bot}$ assumes the role of the orthogonal space to $\mathcal V$ with respect to the inner product $\left\langle{{\mathcal A}_{\boldsymbol{y}} \,\cdot},{\cdot}
	\right\rangle$.
\end{remark}
Applying the Lax-Milgram Theorem and Riesz Representation Theorem, we establish that ${\mathcal A}_{\boldsymbol{y}}$ is an isomorphism from ${\mathcal H}$ onto ${\mathcal H}^{*}$. Consequently, we obtain the decomposition ${\mathcal H}=\mathcal V \bigoplus \mathcal V^{\bot}$, ensuring that every $\boldsymbol{w}\in {\mathcal H}$ has a unique representation $\boldsymbol{w}=\boldsymbol{w}^* + \boldsymbol{w}^{\bot}$ for some $(\boldsymbol{w}^*,\boldsymbol{w}^{\bot})\in \mathcal V\times \mathcal V^{\bot}$ . Specifically, given any $\boldsymbol{w}$, we can determine $\boldsymbol{w}^*$ and $\boldsymbol{w}^{\bot}$ as follows:
\begin{align*}
\left\langle{{\mathcal A}_{\boldsymbol{y}} \boldsymbol{w}^*},{\boldsymbol{v}}
	\right\rangle= \left\langle{{\mathcal A}_{\boldsymbol{y}} \boldsymbol{w}},{\boldsymbol{v}}
	\right\rangle \quad\forall \boldsymbol{v}\in\mathcal V,
\quad
\text{and} \quad
\boldsymbol{w}^{\bot}:=\boldsymbol{w}-\boldsymbol{w}^*.
\end{align*}
Then, for every $\boldsymbol{v}\in \mathcal V$ we have
\begin{align*}
\left\langle{{\mathcal A}_{\boldsymbol{y}}\boldsymbol{w}^{\bot}},{\boldsymbol{v}}
	\right\rangle
=
\left\langle{{\mathcal A}_{\boldsymbol{y}}(\boldsymbol{w}-\boldsymbol{w}^*)},{\boldsymbol{v}}
	\right\rangle
=
\left\langle{{\mathcal A}_{\boldsymbol{y}}\boldsymbol{w}},{\boldsymbol{v}}
	\right\rangle
-
\left\langle{{\mathcal A}_{\boldsymbol{y}}\boldsymbol{w}^*},{\boldsymbol{v}}
	\right\rangle
=
0.
\end{align*}
It shows $\boldsymbol{w}^{\bot}\in\mathcal V^{\bot}$.

{To further augment the fundamental framework, we will rely on standard results from functional analysis, particularly the Banach Closed Range Theorem. Let us denote the range and the null space of a given linear operator $G$ by $\mathscr{R}(G)$ and $\mathscr{N}(G)$, respectively. The following lemma establishes the properties of the operator ${\mathcal B}_{\boldsymbol{y}}^{\top}$.}

\begin{lemma}\label{lem: BT isomorphism}
Let $\mathcal V^*:=\left\{\boldsymbol{f}\in {\mathcal H}^*: \left\langle{\boldsymbol{f}},{\boldsymbol{v}}
	\right\rangle=0,\forall \boldsymbol{v}\in \mathcal V\right\}$, for every $\boldsymbol{f}\in \mathcal V^*$ there exists a function $p\in \mathcal L$ such that $\boldsymbol{f}={\mathcal B}^{\top} p$.
\end{lemma}
\begin{proof}
{By Lemma \ref{lem: L H K prop}} $\mathcal L$ is canonically isomorphic to its dual. 
 By the definition of ${\mathcal B}_{\boldsymbol{y}}$ and ${\mathcal B}_{\boldsymbol{y}}^{\top}$ in \eqref{abm def}, we have ${\mathcal B}_{\boldsymbol{y}}\in \mathscr{L}({\mathcal H};\mathcal L) $ is the dual operator of ${\mathcal B}_{\boldsymbol{y}}^{\top}\in\mathscr{L}(\mathcal L;{\mathcal H}^*)$. Due to {the equation \eqref{inf-sup condition}}, the range of operator ${\mathcal B}_{\boldsymbol{y}}^{\top}$, denoted by $\mathscr{R}({\mathcal B}_{\boldsymbol{y}}^{\top})$, is a closed subspace of ${\mathcal H}^*$.

{Since $\mathscr{R}({\mathcal B}_{\boldsymbol{y}}^{\top})$ is closed, we can directly apply the Banach Closed Range Theorem \cite[Chap. VII, Sec. 5]{Yosida1995}, which states that the range of ${\mathcal B}_{\boldsymbol{y}}^{\top}$ coincides with the orthogonal complement of the null space of ${\mathcal B}_{\boldsymbol{y}}$. That is, $\mathscr{R}({\mathcal B}_{\boldsymbol{y}}^{\top})=\mathscr{N}({\mathcal B}_{\boldsymbol{y}})^{\bot}$.}

{Observe that the null space of ${\mathcal B}_{\boldsymbol{y}}$ is precisely $\mathcal V$, since $\mathscr{N}({\mathcal B}_{\boldsymbol{y}})=\left\{\boldsymbol{v}\in{\mathcal H}:\llangle {{\mathcal B}_{\boldsymbol{y}}\boldsymbol{v}},{q}
	\rrangle=0,\,\forall q\in\mathcal L\right\}=\mathcal V$. Consequently, we deduce that
	\begin{align}\label{R BT def}
		\mathscr{R}({\mathcal B}_{\boldsymbol{y}}^{\top})
		&= \mathscr{N}({\mathcal B}_{\boldsymbol{y}})^{\bot} = \left\{\boldsymbol{f}\in{\mathcal H}^*:\left\langle{\boldsymbol{f}},{\boldsymbol{v}}
	\right\rangle=0,\,\forall\boldsymbol{v}\in\mathscr{N}({\mathcal B}_{\boldsymbol{y}})\right\}\notag\\
		&= \left\{\boldsymbol{f}\in{\mathcal H}^*:\left\langle{\boldsymbol{f}},{\boldsymbol{v}}
	\right\rangle=0,\,\forall\boldsymbol{v}\in\mathcal V\right\} = \mathcal V^*.
	\end{align}
	This means that $\mathscr{R}({\mathcal B}_{\boldsymbol{y}}^{\top}) = \mathcal V^*$, which completes the proof.}

\end{proof}
To avoid any potential confusion, we emphasize  that $\mathcal V^*$ is not the dual space of $\mathcal V$. Rather, $\mathcal V^*$ serves as the dual space of $\mathcal V^{\bot}$. The following lemma elucidates the properties of $\mathcal V^{\bot}$ and $\mathcal V^*$.
\begin{lemma}\label{lem: VT dual}
	Under Assumption \ref{A,B,M assump} and \ref{alpha beta assump}, recalling the definition of $\mathcal V^*$ in Lemma \ref{lem: BT isomorphism} we have
\begin{itemize}
\item[(i)] 	$	\mathcal V^{\bot}
		=
		\left\{\boldsymbol{w}\in {\mathcal H}: \boldsymbol{w} = {\mathcal A}_{\boldsymbol{y}}^{-1} {\mathcal B}_{\boldsymbol{y}}^{\top} q,\, q\in\mathcal L\right\},
	$
\item[(ii)] $\mathcal V^*$ is the dual of $\mathcal V^{\bot}$ 
.
\end{itemize}
\end{lemma}
\begin{proof}
To prove $(i)$, on the one hand, for every $q\in \mathcal L$ we have to show that $\boldsymbol{w}={\mathcal A}_{\boldsymbol{y}}^{-1} {\mathcal B}_{\boldsymbol{y}}^{\top} q \in \mathcal V^{\bot}$. Indeed, for every $\boldsymbol{v}\in \mathcal V$, there holds
\begin{align*}
\left\langle{{\mathcal A}_{\boldsymbol{y}}\boldsymbol{w}},{\boldsymbol{v}}
	\right\rangle
=
\left\langle{{\mathcal A}_{\boldsymbol{y}}({\mathcal A}_{\boldsymbol{y}}^{-1} {\mathcal B}_{\boldsymbol{y}}^{\top} q)},{\boldsymbol{v}}
	\right\rangle
=
\left\langle{ {\mathcal B}_{\boldsymbol{y}}^{\top} q},{\boldsymbol{v}}
	\right\rangle
=
\llangle {{\mathcal B}_{\boldsymbol{y}} \boldsymbol{v}},{ q}
	\rrangle=0.
\end{align*}
It shows that $	\left\{\boldsymbol{w}\in {\mathcal H}~:~\boldsymbol{w} = {\mathcal A}_{\boldsymbol{y}}^{-1} {\mathcal B}_{\boldsymbol{y}}^{\top} q,\, {q\in\mathcal L}\right\} \subset\mathcal V^{\bot}$. On the other hand, let $\boldsymbol{w}\in \mathcal V^{\bot}$ and consider the linear functional $\boldsymbol{f}\in{\mathcal H}^*$ by $\boldsymbol{f}:={\mathcal A}_{\boldsymbol{y}}\boldsymbol{w}$. By definition of $\mathcal V$ and $\mathcal V^{\bot}$, for every $\boldsymbol{v}\in \mathcal V$ there hold
\begin{align*}
\left\langle{\boldsymbol{f}},{\boldsymbol{v}}
	\right\rangle=\left\langle{{\mathcal A}_{\boldsymbol{y}}\boldsymbol{w}},{\boldsymbol{v}}
	\right\rangle=0
\end{align*}
It implies that $\boldsymbol{f}\in \mathcal V^*$. Therefore, in view of Lemma \ref{lem: BT isomorphism}, there exists a function $q\in\mathcal L$ such that ${\mathcal A}_{\boldsymbol{y}}\boldsymbol{w}={\mathcal B}^{\top}q$. The Lax-Milgram Lemma shows that ${\mathcal A}_{\boldsymbol{y}}$ is isomorphic, then $\boldsymbol{w}$ has indeed the expression $\boldsymbol{w}={\mathcal A}^{-1}{\mathcal B}^{\top}q$, and hence $\mathcal V^{\bot}	\subset \left\{\boldsymbol{w}\in {\mathcal H} ~:~ \boldsymbol{w} = {\mathcal A}_{\boldsymbol{y}}^{-1} {\mathcal B}_{\boldsymbol{y}}^{\top} q,\, q\in\mathcal L\right\} $. This completes the proof for $(i)$.

For statement $(ii)$, notice that ${\mathcal A}_{\boldsymbol{y}}$ is an isomorphism from ${\mathcal H}$ into ${\mathcal H}^*$ and recall definitions of $\mathcal V$, $\mathcal V^{\bot}$ and $\mathcal V^*$, we have
\begin{align*}
\mathcal V^*
&=
\left\{\boldsymbol{f} \in{\mathcal H}^*~:~\boldsymbol{f}={\mathcal A}_{\boldsymbol{y}}\boldsymbol{w},\boldsymbol{w}\in{\mathcal H} \text{ and } \left\langle{{\mathcal A}_{\boldsymbol{y}}\boldsymbol{w}},{\boldsymbol{v}}
	\right\rangle=0,\forall\boldsymbol{v}\in\mathcal V\right\}
\\&=
\left\{\boldsymbol{f} \in{\mathcal H}^*~:~\boldsymbol{f}={\mathcal A}_{\boldsymbol{y}}\boldsymbol{w},\boldsymbol{w}\in\mathcal V^{\bot}\right\}.
\end{align*} 
This shows that $\mathcal V^*$ is dual of $\mathcal V^{\bot}$ and hence completes the proof.
\end{proof}
With the fundamental elements now in place, we can articulate the essential properties encapsulated in the $\inf$-$\sup$ condition specified in the assumption. This condition stands as a pivotal feature, laying the groundwork for the forthcoming analysis of parametric regularity.
\begin{lemma}
	The three following properties are equivalent:
	\begin{enumerate}[(i)]
		\item there exists a constant $\beta>0$ such that
		\begin{align}\label{inf sup const}
			\inf_{0\neq q\in \mathcal L} \,\sup_{\boldsymbol{0} \neq \boldsymbol{v} \in {\mathcal H}} 
			\frac{b( \boldsymbol{v},q;B(\boldsymbol{y}))}{\left\|{q}\right\|_{\mathcal L}\left\|{\boldsymbol{v}}\right\|_{{\mathcal H}}}\geq \beta,
		\end{align}
		\item the operator ${\mathcal B}_{\boldsymbol{y}}^{\top}$ is an isomorphism from $\mathcal L$ onto $\mathcal V^*$ and
		\begin{align}\label{BT isomorphism}
			\left\|{{\mathcal B}_{\boldsymbol{y}}^{\top}q}\right\|_{{\mathcal H}^*}
			\geq
			\beta
			\left\|{q}\right\|_{\mathcal L}\quad \forall q \in \mathcal L,
		\end{align}
		\item the operator ${\mathcal B}_{\boldsymbol{y}}$ is an isomorphism from $\mathcal V^{\bot}$ onto $\mathcal L$ and
		\begin{align}\label{B isomorphism}
			\left\|{{\mathcal B}_{\boldsymbol{y}} \boldsymbol{v}}\right\|_{\mathcal L}
			\geq
			\beta
			\left\|{\boldsymbol{v}}\right\|_{{\mathcal H}}\quad \forall \boldsymbol{v} \in \mathcal V^{\bot}.
		\end{align}
	\end{enumerate}
\end{lemma}
\begin{proof}
Firstly, let us show that the properties $(i)$ and $(ii)$ are equivalent. By definition of the operator ${\mathcal B}_{\boldsymbol{y}}^{\top}$ in \eqref{abm def}, the equation \eqref{inf sup const} implies
\begin{align*}
\left\|{{\mathcal B}_{\boldsymbol{y}}^{\top}q}\right\|_{{\mathcal H}^*}
=
\sup_{\boldsymbol{0} \neq \boldsymbol{v} \in {\mathcal H}}
\frac{\left\langle{{\mathcal B}_{\boldsymbol{y}}^{\top}q},{\boldsymbol{v}}
	\right\rangle}{\left\|{\boldsymbol{v}}\right\|_{{\mathcal H}}}
=
\sup_{\boldsymbol{0} \neq \boldsymbol{v} \in {\mathcal H}}
\frac{b(\boldsymbol{v} ,q;B(\boldsymbol{y}))}{\left\|{\boldsymbol{v}}\right\|_{{\mathcal H}}}
\geq
\beta
\left\|{q}\right\|_{\mathcal L}, \qquad \forall q \in \mathcal L,
\end{align*}
so that \eqref{BT isomorphism} is equivalent to \eqref{inf sup const}. Hence, $(ii)$ implies $(i)$. In order to prove $(i)$ implies $(ii)$, it remains to show that, under the condition \eqref{BT isomorphism}, ${\mathcal B}_{\boldsymbol{y}}^{\top}$ is isomorphism from $\mathcal L$ onto $\mathcal V^*$. Clearly, it follows from \eqref{BT isomorphism} that ${\mathcal B}_{\boldsymbol{y}}$ is a one-to-one mapping from $\mathcal L$ onto its range $\mathscr{R}({\mathcal B}_{\boldsymbol{y}}^{\top})$ with a continuous inverse. Thanks to \eqref{R BT def}, we have $\mathscr{R}({\mathcal B}_{\boldsymbol{y}}^{\top})=\mathcal V^*$, it completes the proof for $(i)\Leftrightarrow(ii)$.

We will now demonstrate that $(ii) \Leftrightarrow (iii)$. Notice Lemma \ref{lem: L H K prop} and Lemma \ref{lem: VT dual}, we have $\mathcal L$ is canonically isomorphic to its dual and $\mathcal V^*$ is dual of $\mathcal V^{\bot}$. Applying Banach closed range theorem \cite[Chap. VII, Sec. 5]{Yosida1995}, we have ${\mathcal B}$ is an isomorphism {from $\mathcal V^{\bot}$ into $\mathcal L$} if and only if ${\mathcal B}^{\top}$ is an isomorphism {from $\mathcal L$ into $\mathcal V^*$}. Consequently, the equations \eqref{BT isomorphism} and \eqref{B isomorphism} are equivalent. We finish the proof.
\end{proof}
 
\subsection{Well-posedness under small given data assumption}\label{sec: NS theory}
We now delve into the well-posedness analysis of the Navier-Stokes equation. For a fixed $\boldsymbol{y} \in U$, the variational formulation of \eqref{NS gen equation} reads
\begin{subequations}
	\label{NS variational form}
	\begin{align}
		a\left(\boldsymbol{u},\boldsymbol{v};A(\boldsymbol{y})\right)
		+
		m\left(\boldsymbol{u},\boldsymbol{u},\boldsymbol{v};M(\boldsymbol{y})\right)
		+
		b\left( \boldsymbol{v} ,p;B(\boldsymbol{y})\right)
		&=
		\left\langle{\boldsymbol{f}(\boldsymbol{y})},{\boldsymbol{v}}
	\right\rangle \label{variational form 1}, \\
		b\left( \boldsymbol{u},q; B(\boldsymbol{y})\right)
		&=
		\llangle {g(\boldsymbol{y})},{q}
	\rrangle ,
		\label{variational form 2}
	\end{align}
\end{subequations}
for every $\boldsymbol{v}\in {\mathcal H}$ and $ q\in \mathcal L$.

	For all $\boldsymbol{u}(\boldsymbol{y})\in {\mathcal H}$, there exists unique $(\boldsymbol{u}^*(\boldsymbol{y}),\boldsymbol{u}^{\bot}(\boldsymbol{y}))\in \mathcal V\times \mathcal V^{\bot}$ such that~{$\boldsymbol{u}(\boldsymbol{y})=\boldsymbol{u}^*(\boldsymbol{y})+\boldsymbol{u}^{\bot}(\boldsymbol{y})$}. We now rewrite the decoupling problem as follows
	\begin{subequations}
		\label{decomposition variational}
			\begin{align}
			a\left( \boldsymbol{u}^*,\boldsymbol{v};A(\boldsymbol{y})\right)
			+
			m\left(\boldsymbol{u},\boldsymbol{u},\boldsymbol{v};M(\boldsymbol{y})\right)
			&=
			\left\langle{\boldsymbol{f}(\boldsymbol{y})},{\boldsymbol{v}}
	\right\rangle, \quad \forall \boldsymbol{v}\in \mathcal V\label{decomposition variational 1}, \\
			b\left(\boldsymbol{u}^{\bot},q;B(\boldsymbol{y})\right)
			&=
			\llangle {g(\boldsymbol{y})},{q}
	\rrangle, \quad \forall q\in \mathcal L.
			\label{decomposition variational 2}  
		\end{align}
	\end{subequations}

	\begin{assumption} \label{small f Assump}
For every $\boldsymbol{y}\in U$, the force term $f$ and source mass term $g$ are uniformly bounded in ${\mathcal H}^*$ and $\mathcal L$, respectively. In other words, for all $\boldsymbol{y}\in U$ there {exist} two constants such that
\begin{align*}
\left\|{f(\boldsymbol{y})}\right\|_{{\mathcal H}^*}\leq \frac{\overline{f}}{2}, \qquad
\left\|{g(\boldsymbol{y})}\right\|_{\mathcal L}\leq \frac{\overline{g}}{2}.
\end{align*} 
Furthermore, the constants $\overline{f},\overline{g}, \overline{m}$ and $\alpha$ must fulfil the following conditions:
		\begin{align*}
		\frac{\overline{f}}{\alpha}
		+
		\frac{\overline{g}}{\beta}
		<
		\frac{\alpha}{\overline{m}}.
		\end{align*} 
	\end{assumption}
\begin{remark}
Assumption \ref{small f Assump} also implies the existence of a positive constant $\gamma\in (0,\frac{\alpha}{\overline{m}})$ such that
		\begin{align}\label{gamma properties 1}
		\frac{\overline{f}}{\alpha}
		+
		\frac{\overline{g}}{\beta}
		=
		\frac{\alpha}{\overline{m}}
		-
		\frac{\alpha}{\overline{m}}\left(1-\frac{\overline{m}}{\alpha}\gamma\right)^2
		\end{align} 
	or
\begin{align}\label{gamma properties 2}
\frac{\overline{f}}{2\alpha}
		+
		\frac{\overline{g}}{2\beta}
+
\frac{\overline{m}}{2\alpha} \gamma^2= \gamma.
\end{align}
\end{remark}
\begin{remark}
	{Assumption \ref{small f Assump} serves as a natural extension of the standard small data conditions typically found in the literature (cf. \cite[Theorem 2.2, p. 287]{Girault1986}, \cite[Theorem 1.3, p. 167]{Temam1977}, and \cite[Eqn. (3.17)]{Cohen2018}) to the broader setting where the source mass term is non-zero, meaning $\left\|{g}\right\|_{\mathcal L}\neq 0$.}
\end{remark}
We now show that Assumption \ref{small f Assump} also guarantees that \eqref{decomposition variational} has a unique solution by means of the Banach fixed-point theorem. Let $\boldsymbol{u}^*_0(\boldsymbol{y})=\boldsymbol{0}$ and $\boldsymbol{u}^{\bot}={\mathcal B}_{\boldsymbol{y}}^{-1}g$, we define {the sequence $\left\{\boldsymbol{u}_n\right\}$ as follows}
\begin{align}\label{un bound}
\boldsymbol{u}_{n}(\boldsymbol{y}):=\boldsymbol{u}^*_{n}+\boldsymbol{u}^{\bot}=\boldsymbol{u}^*_{n}(\boldsymbol{y})+{\mathcal B}_{\boldsymbol{y}}^{-1}g(\boldsymbol{y}),
\quad
\text{and}\quad
\left\|{\boldsymbol{u}_n}\right\|_{{\mathcal H}}
\leq
\left\|{\boldsymbol{u}^*_n}\right\|_{{\mathcal H}}
+
\frac{\overline{g}}{2\beta},
\end{align}
 where $\boldsymbol{u}^*_{n}(\boldsymbol{y})$ is the unique solution of
	\begin{align}\label{contracting mapping}
		a\left( \boldsymbol{u}^*_{n},\boldsymbol{v};A(\boldsymbol{y})\right)
	=
	\left\langle{\boldsymbol{f}(\boldsymbol{y})},{\boldsymbol{v}}
	\right\rangle
	-
	m\left(\boldsymbol{u}_{n-1},\boldsymbol{u}_{n-1},\boldsymbol{v};M(\boldsymbol{y})\right)
	 \quad \forall \boldsymbol{v}\in \mathcal V.
	\end{align}
	The following Lemma shows that the sequence $\{\boldsymbol{u}_n(\boldsymbol{y})\}$ never leaves the closed set 
	$$\boldsymbol{B}(\boldsymbol{0},\gamma):=\left\{\boldsymbol{v}\in {\mathcal H}: \left\|{\boldsymbol{v}}\right\|_{{\mathcal H}}\leq \gamma\right\}$$ 
	and converges {to} a limit in $\boldsymbol{B}(\boldsymbol{0},\gamma)$. Obviously, when the sequence $\left\{\boldsymbol{u}_{n}(\boldsymbol{y})\right\}$ admits a limit point, it would be a solution of \eqref{NS variational form}. Indeed, the following Lemma proves the above statement.
	\begin{lemma}\label{lem: NS u bound}
		For every $\boldsymbol{y}\in U$, the sequence $\left\{\|\boldsymbol{u}_{n}(\boldsymbol{y})\|_{{\mathcal H}}\right\}$ is bounded by $\gamma \in (0, \frac{\alpha}{\overline{m}})$ and $\{\boldsymbol{u}_n(\boldsymbol{y})\}$ converges to a fixed point in $\boldsymbol{B}(\boldsymbol{0},\gamma)$.
	\end{lemma}
	\begin{proof}
		We will prove boundedness of the sequence by induction with respect to $n$. Recall the equation \eqref{un bound} with $\boldsymbol{u}^*_0(\boldsymbol{y})=\boldsymbol{0}$ and \eqref{gamma properties 2}, we obtain
\begin{align}
\left\|{\boldsymbol{u}_0(\boldsymbol{y})}\right\|_{{\mathcal H}}
\leq
\frac{\overline{g}}{2\beta}
\leq
\frac{\overline{f}}{2\alpha}
+
\frac{\overline{g}}{2\beta}
+\frac{\overline{m}}{2\alpha} \gamma^2
= 
\gamma
\end{align}
		and hence ${\boldsymbol{u}_0}(\boldsymbol{y})\in \boldsymbol{B}(\boldsymbol{0},\gamma)$. Assume now that all $\boldsymbol{u}_1, \dots, \boldsymbol{u}_n$ belong to this neighbourhood and prove that the same holds for $\boldsymbol{u}_{n+1}$. For this, we substitute $\boldsymbol{v}={\boldsymbol{u}_{n+1}^*}(\boldsymbol{y})$ into \eqref{contracting mapping} and recall \eqref{NS ab-bounds}, \eqref{m operator bound} and Assumption \ref{small f Assump} to obtain
		\begin{align*}
			\left\|{\boldsymbol{u}^*_{n+1}(\boldsymbol{y})}\right\|_{{\mathcal H}}
			\leq
			\frac{\overline{f}}{2\alpha}
			+
			\frac{\overline{m}}{2\alpha}
			\left\|{u_n(\boldsymbol{y})}\right\|_{{\mathcal H}}^2
			\leq
			\frac{\overline{f}}{2\alpha}
			+
			\frac{\overline{m}}{2\alpha} \gamma^2.
		\end{align*}
		This together with \eqref{un bound} and \eqref{gamma properties 2} show that $\boldsymbol{u}_{n}\in \boldsymbol{B}(\boldsymbol{0},\gamma)$ for all $n$. We now prove that the sequence converges to a limit in ${\mathcal H}$. {We have that $\boldsymbol{u}^*_{n+1}(\boldsymbol{y})$ is the solution of 
		\begin{align*}
			a\left( \boldsymbol{u}^*_{n+1},\boldsymbol{v};A(\boldsymbol{y})\right)
			=
			\left\langle{\boldsymbol{f}(\boldsymbol{y})},{\boldsymbol{v}}
	\right\rangle
			-
			m\left(\boldsymbol{u}_{n},\boldsymbol{u}_{n},\boldsymbol{v};M(\boldsymbol{y})\right)
			\quad \forall \boldsymbol{v}\in \mathcal V.
		\end{align*}
		 Subtract both sides of \eqref{contracting mapping} from the above equation and set $\boldsymbol{v}=\boldsymbol{u}^*_{n+1}(\boldsymbol{y})-u^*_n(\boldsymbol{y})$ to obtain
		 }\begin{align*}
a\left({\boldsymbol{u}^*_{n+1}-\boldsymbol{u}^*_{n}},\boldsymbol{u}^*_{n+1}-\boldsymbol{u}^*_{n};A(\boldsymbol{y})\right)
&=
m\left(\boldsymbol{u}_{n-1},\boldsymbol{u}_{n-1},\boldsymbol{u}^*_{n+1}-\boldsymbol{u}^*_{n};M(\boldsymbol{y})\right)\\
&\quad-
m\left(\boldsymbol{u}_{n},\boldsymbol{u}_{n},\boldsymbol{u}^*_{n+1}-\boldsymbol{u}^*_{n};M(\boldsymbol{y})\right),
		\end{align*} 
		The left-hand side is bounded by $\alpha\|\boldsymbol{u}^*_{n+1} - \boldsymbol{u}^*_n\|_{{\mathcal H}}^2$ from below. To obtain an upper bound for the right-hand side, add and subtract the term $m\left(\boldsymbol{u}_{n-1},\boldsymbol{u}_{n},\boldsymbol{u}^*_{n+1}-\boldsymbol{u}^*_{n};M(\boldsymbol{y})\right)$ and use~\eqref{m operator bound}. This implies
		\begin{equation} \label{b term bound}
			\begin{split}
				\|\boldsymbol{u}^*_{n+1} - \boldsymbol{u}^*_n\|_{{\mathcal H}}^2 &\leq 
				\frac{\overline{m}}{2\alpha} \left\|{\boldsymbol{u}^*_{n+1}-\boldsymbol{u}^*_{n}}\right\|_{{\mathcal H}}
				\left\|{\boldsymbol{u}_{n}-\boldsymbol{u}_{n-1}}\right\|_{{\mathcal H}}
				\left(\left\|{\boldsymbol{u}_n}\right\|_{{\mathcal H}}+\left\|{\boldsymbol{u}_{n-1}}\right\|_{{\mathcal H}}\right)
				\\
				&\leq
				\frac{\overline{m}\,\gamma}{\alpha}
\left\|{\boldsymbol{u}^*_{n+1}-\boldsymbol{u}^*_{n}}\right\|_{{\mathcal H}}
				\left\|{\boldsymbol{u}_{n}-\boldsymbol{u}_{n-1}}\right\|_{{\mathcal H}},
			\end{split}
		\end{equation}
		where $\left\{\boldsymbol{u}_{n}\right\}\subset \boldsymbol{B}(\boldsymbol{0},\gamma)$ has been used in the last step. Due to $\boldsymbol{u}_{n+1}=\boldsymbol{u}^*_{n+1}+\boldsymbol{u}^{\bot}$ for all $n\in \mathbb N$, we have $\boldsymbol{u}_{n+1} - \boldsymbol{u}_n=\boldsymbol{u}^*_{n+1} - \boldsymbol{u}^*_n$. Now we cancel  $ \left\|{\boldsymbol{u}_{n+1}-\boldsymbol{u}_{n}}\right\|_{{\mathcal H}}$ from both sides imply the contraction property and obtain
		\begin{align*}
			\left\|{\boldsymbol{u}_{n+1}-\boldsymbol{u}_{n}}\right\|_{{\mathcal H}}
\leq
			\frac{\overline{m}\,\gamma}{\alpha}
			\left\|{\boldsymbol{u}_{n}-\boldsymbol{u}_{n-1}}\right\|_{{\mathcal H}}.
		\end{align*}
		Since $\gamma\in (0,\frac{\alpha}{\overline{m}})$ and $\frac{\overline{m}\,\gamma}{\alpha}<1$, the sequence $\left\{\boldsymbol{u}_{n}(\boldsymbol{y})\right\}$ converges to a fixed point in $ \boldsymbol{B}(\boldsymbol{0},\gamma)$ by the Banach fixed point theorem.
	\end{proof}	 We can now establish the upper bounds for the norms of the velocity field and the pressure, as stated in the following lemma.
\begin{lemma}\label{u p bound}
For every $\boldsymbol{y}\in U$ and $(\boldsymbol{u}(\boldsymbol{y}),p(\boldsymbol{y}))$ is solution of \eqref{NS variational form}, we have
\begin{align}
\left\|{\boldsymbol{u}(\boldsymbol{y})}\right\|_{{\mathcal H}}
&\leq
\gamma=:\overline{u},\label{def u bar} 
\\
\left\|{p(\boldsymbol{y})}\right\|_{\mathcal L}
&\leq
\frac{1}{\beta}
\left(\frac{\overline{f}}{2}
+
\frac{\overline{a}\gamma}{2}
+
\frac{\overline{m}\gamma^2}{2}
\right)
=: \overline{p}.\label{def p bar} 
\end{align}
\end{lemma}
\begin{proof}
From Lemma \ref{lem: NS u bound}, we have that $\boldsymbol{u} \in \boldsymbol{B}(\boldsymbol{0},\gamma)$, it directly shows \eqref{def u bar}. In order to show \eqref{def p bar}, notice \eqref{variational form 1} we obtain
\begin{align*}
		b\left( \boldsymbol{v} ,p(\boldsymbol{y});B(\boldsymbol{y})\right)
		=
		\left\langle{\boldsymbol{f}(\boldsymbol{y})},{\boldsymbol{v}}
	\right\rangle
		-
a\left( \boldsymbol{u}(\boldsymbol{y}),\boldsymbol{v};A(\boldsymbol{y})\right)
		-
		m\left(\boldsymbol{u}(\boldsymbol{y}),\boldsymbol{u}(\boldsymbol{y}),\boldsymbol{v};M(\boldsymbol{y})\right)
		.
\end{align*}
This together with \eqref{B isomorphism} and Assumption \ref{A,B,M assump} implies \eqref{def p bar}. We finish the proof.
\end{proof}

	For each $\boldsymbol{y}\in U$, we denote by $t_{\boldsymbol{y}}(\boldsymbol{u}, \boldsymbol{w},\boldsymbol{v})$ the linearization of \eqref{variational form 1} mapping from $ \mathcal V \times \mathcal V\times \mathcal V \rightarrow \mathbb R$ as
	\begin{align} \label{Atilde-def}
		t_{\boldsymbol{y}}(\boldsymbol{u},\boldsymbol{w},\boldsymbol{v})
		= a( \boldsymbol{w}, \boldsymbol{v};A(\boldsymbol{y})) 
		+
		m(\boldsymbol{u},\boldsymbol{w},\boldsymbol{v};M(\boldsymbol{y}))+m(\boldsymbol{w},\boldsymbol{u},\boldsymbol{v};M(\boldsymbol{y})).
	\end{align}
The following lemma demonstrates the coercivity of $t_{\boldsymbol{y}}$, which is essential for the regularity proof in Section \ref{sec: NS main result}.
	\begin{lemma}\label{lem: NS coercive-type}
		Under Assumption \ref{A,B,M assump}, Assumption \ref{alpha beta assump} and Assumption \ref{small f Assump}, the operator $t(\boldsymbol{y})$ is uniformly coercive in $\boldsymbol{y}$, i.e.
		\begin{align}\label{coercive lower bound}
			t_{\boldsymbol{y}}(\boldsymbol{u}, \boldsymbol{v},\boldsymbol{v} )\geq
			\widetilde{\alpha} \left\|{\boldsymbol{v}}\right\|_{{\mathcal H}}^2,
			\quad\forall \boldsymbol{v} \in {\mathcal H} \text{ and } \forall \boldsymbol{u} \in \boldsymbol{B}(0,\gamma),
		\end{align}
		where $\widetilde{\alpha}:= \alpha-\overline{m} \gamma>0$.
	\end{lemma}
	\begin{proof}
	Employing \eqref{V-coercive}, \eqref{m operator bound} and notice that $\boldsymbol{u}\in\boldsymbol{B}(0,\gamma)$  we have
		\begin{align*}
			t_{\boldsymbol{y}}(\boldsymbol{u}, \boldsymbol{v},\boldsymbol{v} )\geq
			\left(\alpha-
				\overline{m}\left\|{\boldsymbol{u}}\right\|_{{\mathcal H}}\right)\left\|{\boldsymbol{v}}\right\|_{{\mathcal H}}^2
\geq
\left(\alpha-
				\overline{m}\gamma\right)\left\|{\boldsymbol{v}}\right\|_{{\mathcal H}}^2.
		\end{align*}
		This together with $\gamma\in(0,\frac{\alpha}{\overline{m}})$ shows that $\widetilde{\alpha}:= \alpha-\overline{m} \gamma>0$ and finishes the proof.
	\end{proof}

\section{Parametric regularity for Navier-Stokes equation}\label{sec: NS main result}
We now make an assumption on the coefficients, which, in particular, ensure that the solution of the Navier-Stokes problem \eqref{NS gen equation} is Gevrey-class regular.
	\begin{assumption} \label{NS Assumption}
		For all fixed values $\boldsymbol{y}\in U {\subset} \mathbb R^s$ with {$s<\infty$, the} matrices $A(\boldsymbol{y})$, $B(\boldsymbol{y})$ and $M(\boldsymbol{y})$ are of Gevrey class $G^\delta(U,L^{\infty}(D)^{d\times d})$, $\boldsymbol{f}(\boldsymbol{y})$ is of Gevrey class $G^\delta(U,{\mathcal H}^*)$ and $g(\boldsymbol{y})$ is of Gevrey class $G^\delta(U,\mathcal L)$, i.e. for all $\boldsymbol{\nu}\in \mathbb N^s$ there exist $\boldsymbol{R}$ independent of $s$ such that
\begin{align*}
\max\left\{
\frac{\left\|{\partial^{\boldsymbol{\nu}} A(\boldsymbol{y})}\right\|_{\infty}}{\overline{a}},
\frac{\left\|{\partial^{\boldsymbol{\nu}} B(\boldsymbol{y})}\right\|_{\infty}}{\overline{b}},
\frac{\left\|{\partial^{\boldsymbol{\nu}} M(\boldsymbol{y})}\right\|_{\infty}}{C_4^{-2}\overline{m}}
\right\}
\leq
\frac{1}{2}
\frac{(\left|\boldsymbol{\nu}\right|!)^\delta }{(2\boldsymbol{R})^{\boldsymbol{\nu}}}
\end{align*}
and
		\begin{align*}
			\left\|{\partial^{\boldsymbol{\nu}} \boldsymbol{f}(\boldsymbol{y})}\right\|_{{\mathcal H}^*}
			\leq
			\frac{\overline{f}}{2}
			\frac{(\left|\boldsymbol{\nu}\right|!)^\delta }{(2\boldsymbol{R})^{\boldsymbol{\nu}}} , \qquad
			\left\|{\partial^{\boldsymbol{\nu}} g(\boldsymbol{y})}\right\|_{\mathcal L}
			\leq
			\frac{\overline{g}}{2}
			\frac{(\left|\boldsymbol{\nu}\right|!)^\delta }{(2\boldsymbol{R})^{\boldsymbol{\nu}}}  .
		\end{align*} 
	\end{assumption}
	Notice that for $\boldsymbol{\nu} = \boldsymbol{0}$, Assumption \ref{NS Assumption} agrees with the upper bounds in \eqref{NS ab-bounds}. Notice also that the components of $\boldsymbol{R}$ are readily scaled by the factor of $2$. This leads to no loss of generality, but helps to shorten the forthcoming expressions. 
	For example, in view of \eqref{multiindex-est-2}, Assumption \ref{NS Assumption} immediately implies
\begin{equation}\label{abf assumption}
\begin{split}
&\max\left\{
\frac{\left\|{\partial^{\boldsymbol{\nu}} A(\boldsymbol{y})}\right\|_{\infty}}{\overline{a}},
\frac{\left\|{\partial^{\boldsymbol{\nu}} B(\boldsymbol{y})}\right\|_{\infty}}{\overline{b}},
\frac{\left\|{\partial^{\boldsymbol{\nu}} M(\boldsymbol{y})}\right\|_{\infty}}{C_4^{-2}\overline{m}}
\right\}
\leq
\frac{
 \left[\tfrac{1}{2} \right]_{\left|\boldsymbol{\nu}\right|} }{\boldsymbol{R}^{\boldsymbol{\nu}}}(\left|\boldsymbol{\nu}\right|!)^{\delta-1} ,\\
&
			\left\|{\partial^{\boldsymbol{\nu}} \boldsymbol{f}(\boldsymbol{y})}\right\|_{{\mathcal H}^*}
			\leq
			\frac{{\overline{f}}
 \left[\tfrac{1}{2} \right]_{\left|\boldsymbol{\nu}\right|} }{\boldsymbol{R}^{\boldsymbol{\nu}}}(\left|\boldsymbol{\nu}\right|!)^{\delta-1}  , \qquad
			\left\|{\partial^{\boldsymbol{\nu}} g(\boldsymbol{y})}\right\|_{\mathcal L}
			\leq
			\frac{{\overline{g}}
 \left[\tfrac{1}{2} \right]_{\left|\boldsymbol{\nu}\right|} }{\boldsymbol{R}^{\boldsymbol{\nu}}}(\left|\boldsymbol{\nu}\right|!)^{\delta-1} .
\end{split}
\end{equation}
The subsequent theorem constitutes a crucial element for establishing the main result in this paper.
\begin{theorem}\label{gevrey regularity for NS}
	Let the matrices $A,B,M$, the force term $\boldsymbol{f}$ and the source mass term $g$ of \eqref{NS variational form} satisfy  Assumption \ref{small f Assump} and  Assumption \ref{NS Assumption} for some $\delta \geq 1$. Additionally, assume that both Assumption \ref{A,B,M assump} and Assumption \ref{alpha beta assump} also  hold. The following estimates hold for all $\boldsymbol{\nu} \in \mathcal F\setminus\left\{\boldsymbol{0}\right\}$ and $\boldsymbol{y} \in U$
	\begin{align}\label{u H norm}
		\left\|{\partial^{\boldsymbol{\nu}} \boldsymbol{u}( \boldsymbol{y})}\right\|_{{\mathcal H}}
		\leq
		\frac{C_u \rho^{\left|\boldsymbol{\nu}\right|-1} \left[\tfrac{1}{2} \right]_{\left|\boldsymbol{\nu}\right|} }{\boldsymbol{R}^{{\boldsymbol{\nu}}}} (\left|\boldsymbol{\nu}\right|!)^{\delta-1}
	\end{align}
	and
	\begin{align}\label{p L norm}
		\left\|{\partial^{\boldsymbol{\nu}} p( \boldsymbol{y})}\right\|_{\mathcal L}
		\leq
		\frac{C_p \rho^{\left|\boldsymbol{\nu}\right|-1} \left[\tfrac{1}{2} \right]_{\left|\boldsymbol{\nu}\right|} }{\boldsymbol{R}^{{\boldsymbol{\nu}}}} (\left|\boldsymbol{\nu}\right|!)^{\delta-1}.
	\end{align}
	The constants in the above bounds are explicitly determined as $C_u:=\gamma \sigma_u$,  $C_p:=\overline{p} \sigma_p$, $\rho:=\max\left\{1,\rho_u,\rho_p\right\}$ and
{\small	
\begin{align} \label{NV sigma-rho-def}
\sigma_u
&:=
\frac{1}{\gamma\widetilde{\alpha}}\left(\overline{f}
+
{\overline{b}\overline{p}}
+
\gamma
\left(
\overline{m}
\frac{(\overline{g}+\overline{b}\,\gamma)}{\beta}
+
\overline{a}
+
\overline{m}\,\gamma
\right)\right)
+
\frac{(\overline{g}+\overline{b}\,\gamma)}{\gamma\beta}\notag,\\
\sigma_p
&:=\frac{1}{\beta \overline{p}}
\left(\overline{f}
+
\overline{b}\overline{p}
+
\overline{a}\gamma
\left(1+\frac{\sigma_u}{2}\right)
+
\overline{m}\gamma^2(1+\sigma_u)\right),\\
\rho_u
&:=
\frac{\left(
1
+
\frac{\gamma\,
		\overline{m}}{\widetilde \alpha}
\right)\left(\overline{g}+\overline{b} \,\gamma+2C_u\overline{b}\right)}
{\beta C_u}
+
\frac{
\overline{f}+\overline{b}\,\overline{p}+
\gamma
\left(
\overline{a}+\overline{m}\,\gamma
\right)
+
2C_p\overline{b}
}{\widetilde \alpha C_u}
+
\frac{
			2\overline{a}
			+
			4\overline{m}\,\gamma 
			+
			7\overline{m}C_u}{\widetilde \alpha}
\notag,\\
\rho_p
&:=
\frac{\left(
\overline{f}+\overline{b}\,\overline{p}+
\gamma
\left(
\overline{a}+\overline{m}\,\gamma
\right)
+
2C_u\overline{a}
+
{2C_p \overline{b}}
+
4\overline{m}\,\gamma C_u
+
7\overline{m}C_u^2
\right)}{\beta C_p}
+
\frac{C_u\rho_u}{\beta C_p}
\left(
\frac{\overline{a}}{2}
+
\gamma\,
\overline{m}\right)\notag.
\end{align}}
\end{theorem}
{The proof strategy begins by isolating the highest-order derivative using the recurrence relation \eqref{var form vnu derivative}. We then project the highest derivative of the velocity field $\partial^{\boldsymbol{\nu}+\boldsymbol{e}} \boldsymbol{u}$ onto the spaces $\mathcal V$ and $\mathcal V^{\top}$, denoted by $\boldsymbol{w}^*$ and $\boldsymbol{w}^{\bot}$, respectively. In Lemma \ref{lem:wT-bnd}, we employ the isomorphism property of the operators $B_{\boldsymbol{y}}$ to bound $\boldsymbol{w}^{\bot}$ in terms of lower-order terms. Next, in Lemma \ref{lem:w-bnd}, the coercivity of the bilinear forms $t_{\boldsymbol{y}}(\boldsymbol{u},\cdot,\cdot)$  allows us to bound $\boldsymbol{w}^*$	using the bound for  $\boldsymbol{w}^{\bot}$ and additional lower-order contributions. Finally, in Lemma \ref{lem:p-bnd}, the isomorphism property of the operators $B_{\boldsymbol{y}}^{\top}$	is used to bound the pressure derivative $\partial^{\boldsymbol{\nu}+\boldsymbol{e}} p$ in terms of the velocity derivative bounds and the remaining lower-order pressure derivatives.}
\begin{lemma}\label{lem:wT-bnd}
	For sufficiently regular solutions of \eqref{decomposition variational}, let $\partial^{\boldsymbol{\nu}+\boldsymbol{e}} \boldsymbol{u}=\boldsymbol{w}^*+\boldsymbol{w}^{\bot}$ with $\boldsymbol{w}^*\in \mathcal V$ and $\boldsymbol{w}^{\bot}\in \mathcal V^{\bot}$, where $\boldsymbol{e}$ is a unit multi-index in $\mathcal F$, i.e. $\left|\boldsymbol{e}\right|=1$. There holds
	\begin{align}\label{wT norm upper}
			\beta\left\|{  \boldsymbol{w}^{\bot}}\right\|_{{\mathcal H}}
&\leq
\mathfrak{F}_1(\boldsymbol{\nu})
+
\mathfrak{B}_1(\boldsymbol{\nu})
	\end{align}
	where
	\begin{align*}
		\mathfrak{F}_1(\boldsymbol{\nu})
&:=
\left\|{\partial^{\boldsymbol{\nu}+\boldsymbol{e}}g}\right\|_{\mathcal L}
+
\overline{u}
			\left\|{\partial^{\boldsymbol{\nu}+\boldsymbol{e}} B}\right\|_{\infty},\\
		\mathfrak{B}_1(\boldsymbol{\nu})
		&:=
			\sum_{\boldsymbol{0}<\boldsymbol{\eta} \leq \boldsymbol{\nu}}
			\left(\boldsymbol{\nu} \atop \boldsymbol{\eta}\right)
\left(\left\|{\partial^{\boldsymbol{\nu}+\boldsymbol{e}-\boldsymbol{\eta}} B}\right\|_{\infty}
			\left\|{	\partial^{\boldsymbol{\eta}} \boldsymbol{u} }\right\|_{{\mathcal H}}
+
			\left\|{\partial^{\boldsymbol{\eta}} B}\right\|_{\infty}
			\left\|{	 \partial^{\boldsymbol{\nu}+\boldsymbol{e}-\boldsymbol{\eta}} \boldsymbol{u} }\right\|_{{\mathcal H}}\right).
	\end{align*}
\end{lemma}
\begin{proof}
	We recall variational formulation \eqref{NS variational form}. for the bound of $\boldsymbol{w}^{\bot}$, taking the $\boldsymbol{e}$-th derivative both sides of \eqref{decomposition variational 2} to obtain
	\begin{align*}
		b\left(
		\partial^{\boldsymbol{e}} \boldsymbol{u} ,q;B\right)
		+
		b\left(
			 \boldsymbol{u} ,q;\partial^{\boldsymbol{e}}B\right)
			=
			\llangle {\partial^{\boldsymbol{e}} g},{q}
	\rrangle, \quad\forall q\in \mathcal L.
	\end{align*}
	Since $b\left(
		 \boldsymbol{v},q;B\right)=0$ for all $\boldsymbol{v}\in \mathcal V$, we have that $b\left(
		\partial^{\boldsymbol{\nu}+\boldsymbol{e}} \boldsymbol{u},q;B\right)=b\left(
		 \boldsymbol{w}^{\bot},q;B\right)$. Taking the $\boldsymbol{\nu}$-derivative and collecting the highest order $b\left(
		 \boldsymbol{w},q;B\right)$ on the left-hand side to derive
\begin{align*}
			b\left(  \boldsymbol{w}^{\bot},q;B\right)
	=
		\llangle[\Big] {\widehat{f}_1(\boldsymbol{\nu})},{q}
	\rrangle[\Big]
	-
	\llangle[\Big] {\widehat{b}_1(\boldsymbol{\nu})},{q}
	\rrangle[\Big]	 \quad\forall q\in \mathcal L
\end{align*}
where
\begin{align*}
	\llangle[\Big] {\widehat{f}_1(\boldsymbol{\nu})},{q}
	\rrangle[\Big]
&:=
\llangle[\big] {\partial^{\boldsymbol{\nu}+\boldsymbol{e}} g},{q}
	\rrangle[\big]
-
b\left( \boldsymbol{u},q;\partial^{\boldsymbol{\nu}+\boldsymbol{e}}B\right),
\\
	\llangle[\Big] {\widehat{b}_1(\boldsymbol{\nu})},{q}
	\rrangle[\Big]
	&:=
	\sum_{\boldsymbol{0}<\boldsymbol{\eta} \leq \boldsymbol{\nu}}
	\left(\boldsymbol{\nu} \atop \boldsymbol{\eta}\right)
\left(
	b\left(
		\partial^{\boldsymbol{\nu}+\boldsymbol{e}-\boldsymbol{\eta}}  \boldsymbol{u} ,q;\partial^{\boldsymbol{\eta}} B\right)
+
	b\left(
		\partial^{\boldsymbol{\eta}}  \boldsymbol{u},q; \partial^{\boldsymbol{\nu}+\boldsymbol{e}-\boldsymbol{\eta}} B\right)\right).
\end{align*}
Noting \eqref{b operator bound} , we estimate $\left\|{\widehat{f}_1(\boldsymbol{\nu})}\right\|_{\mathcal L}$ and $\left\|{\widehat{b}_1(\boldsymbol{\nu})}\right\|_{\mathcal L}$ as follows
\begin{align*}
\left\|{\widehat{f}_1(\boldsymbol{\nu})}\right\|_{\mathcal L}
	&\leq
\left\|{\partial^{\boldsymbol{\nu}+\boldsymbol{e}} g}\right\|_{\mathcal L}
+
\overline{u}
\left\|{\partial^{\boldsymbol{\nu}+\boldsymbol{e}} B}\right\|_{\infty}
=
\mathfrak{F}_1(\boldsymbol{\nu}),
\\
	\left\|{\widehat{b}_1(\boldsymbol{\nu})}\right\|_{\mathcal L}
	&\leq
	\sum_{\boldsymbol{0}<\boldsymbol{\eta} \leq \boldsymbol{\nu}}
	\left(\boldsymbol{\nu} \atop \boldsymbol{\eta}\right)
\left(
	\left\|{\partial^{\boldsymbol{\nu}+\boldsymbol{e}-\boldsymbol{\eta}} B}\right\|_{\infty}
	\left\|{	\partial^{\boldsymbol{\eta}} \boldsymbol{u} }\right\|_{{\mathcal H}}
+
	\left\|{\partial^{\boldsymbol{\eta}} B}\right\|_{\infty}
	\left\|{	 \partial^{\boldsymbol{\nu}+\boldsymbol{e}-\boldsymbol{\eta}} \boldsymbol{u} }\right\|_{{\mathcal H}}
\right)
	=
\mathfrak{B}_1(\boldsymbol{\nu}).
\end{align*}
	This together with inf-sup condition \eqref{B isomorphism} shows \eqref{wT norm upper}.
\end{proof}
\begin{lemma}\label{lem:w-bnd}
	For sufficiently regular solutions of \eqref{decomposition variational}, and recalling $\partial^{\boldsymbol{\nu}+\boldsymbol{e}} \boldsymbol{u}=\boldsymbol{w}^*+\boldsymbol{w}^{\bot}$ with $\boldsymbol{w}^*\in \mathcal V$, $\boldsymbol{w}^{\bot}\in \mathcal V^{\bot}$, and $\boldsymbol{e}\in \mathcal F$ a unit multi-index, the following holds:
	\begin{align}\label{w* norm upper}
		\widetilde{\alpha}\left\|{  \boldsymbol{w}^*}\right\|_{{\mathcal H}}
		\leq
		\mathfrak{F}_2(\boldsymbol{\nu}) +		
\mathfrak{A}(\boldsymbol{\nu})
		+
		\mathfrak{B}_2(\boldsymbol{\nu})
		+
		\mathfrak{M}(\boldsymbol{\nu})
+
\mathfrak{N}(\boldsymbol{\nu})
+
\overline{u}\,
		\overline{m}
	\left\|{\boldsymbol{w}^{\bot}}\right\|_{{\mathcal H}}
	\end{align}
	where{\small
	\begin{align*}
		\mathfrak{F}_2(\boldsymbol{\nu})
&:=
\left\|{\partial^{\boldsymbol{\nu}+\boldsymbol{e}} \boldsymbol{f}}\right\|_{{\mathcal H}^*}
+
\overline{p}
			\left\|{\partial^{\boldsymbol{\nu}+\boldsymbol{e}}B}\right\|_{\infty}
		+
\overline{u}
\left\|{\partial^{\boldsymbol{\nu}+\boldsymbol{e}}A }\right\|_{\infty}
+		
{C_4^2}\overline{u}^2
\left\|{\partial^{\boldsymbol{\nu}+\boldsymbol{e}}M}\right\|_{\infty},
\\
		\mathfrak{A}(\boldsymbol{\nu})
		&:=
			\sum_{\boldsymbol{0}< \boldsymbol{\eta}\leq \boldsymbol{\nu}}
			\left(\boldsymbol{\nu}\atop \boldsymbol{\eta}\right)
\left(
			\left\|{\partial^{\boldsymbol{\nu}+\boldsymbol{e}-\boldsymbol{\eta}}A }\right\|_{\infty}
			\left\|{\partial^{\boldsymbol{\eta}}  \boldsymbol{u}}\right\|_{{\mathcal H}}
+
			\left\|{\partial^{\boldsymbol{\eta}}A }\right\|_{\infty}
			\left\|{\partial^{\boldsymbol{\nu}+\boldsymbol{e}-\boldsymbol{\eta}} \boldsymbol{u}}\right\|_{{\mathcal H}}\right),\\
		\mathfrak{B}_2(\boldsymbol{\nu})
		&:=
			\sum_{\boldsymbol{0}< \boldsymbol{\eta}\leq \boldsymbol{\nu}}
			\left(\boldsymbol{\nu}\atop \boldsymbol{\eta}\right)
\left(
			\left\|{\partial^{\boldsymbol{\nu}+\boldsymbol{e}-\boldsymbol{\eta}}B}\right\|_{\infty}
			\left\|{\partial^{\boldsymbol{\eta}}p }\right\|_{\mathcal L}
+
			\left\|{\partial^{\boldsymbol{\eta}}B}\right\|_{\infty}
			\left\|{ \partial^{\boldsymbol{\nu}+\boldsymbol{e}-\boldsymbol{\eta}}p }\right\|_{\mathcal L}\right),\\
		\mathfrak{M}(\boldsymbol{\nu})
		&:=
		2\,{C_4^2}\sum_{\boldsymbol{0}< \boldsymbol{\eta}\leq \boldsymbol{\nu}}
		\left(\boldsymbol{\nu}\atop \boldsymbol{\eta}\right)
\overline{u}
		\left\|{\partial^{\boldsymbol{\nu}+\boldsymbol{e}-\boldsymbol{\eta}}M}\right\|_{\infty}
\left\|{\partial^{\boldsymbol{\eta}}\boldsymbol{u}}\right\|_{{\mathcal H}}\\
&\qquad\qquad\qquad
+		2\sum_{\boldsymbol{0}< \boldsymbol{\eta}\leq \boldsymbol{\nu}}
\left(\boldsymbol{\nu}\atop \boldsymbol{\eta}\right)
		\left\|{\partial^{\boldsymbol{\nu}+\boldsymbol{e}-\boldsymbol{\eta}}\boldsymbol{u}}\right\|_{{\mathcal H}}
\left(
		\frac{\overline{m}}{2}
\left\|{\partial^{\boldsymbol{\eta}}\boldsymbol{u}}\right\|_{{\mathcal H}}
+
{C_4^2}\, \overline{u}
		\left\|{\partial^{\boldsymbol{\eta}} M}\right\|_{\infty}\right),\\
		\mathfrak{N}(\boldsymbol{\nu})
		&:=
		2\, {C_4^2}\sum_{\boldsymbol{0}< \boldsymbol{\eta}\leq \boldsymbol{\nu}}
		\sum_{\boldsymbol{0}< \boldsymbol{\ell} < \boldsymbol{\eta}}
		\left(\boldsymbol{\nu}\atop \boldsymbol{\eta}\right)
		\left(\boldsymbol{\eta}\atop \boldsymbol{\ell}\right)
		\left\|{\partial^{\boldsymbol{\ell}}\boldsymbol{u}}\right\|_{{\mathcal H}}
		\left\|{\partial^{\boldsymbol{\eta}-\boldsymbol{\ell}} M}\right\|_{\infty}
		\left\|{\partial^{\boldsymbol{\nu}+\boldsymbol{e}-\boldsymbol{\eta}}\boldsymbol{u}}\right\|_{{\mathcal H}}
\\&\qquad\qquad\qquad+
{C_4^2}\sum_{\boldsymbol{0}< \boldsymbol{\eta}\leq \boldsymbol{\nu}}
\sum_{\boldsymbol{0}< \boldsymbol{\ell} < \boldsymbol{\eta}}
\left(\boldsymbol{\nu}\atop \boldsymbol{\eta}\right)
\left(\boldsymbol{\eta}\atop \boldsymbol{\ell}\right)
\left\|{\partial^{\boldsymbol{\ell}}\boldsymbol{u}}\right\|_{{\mathcal H}}
		\left\|{\partial^{\boldsymbol{\nu}+\boldsymbol{e}-\boldsymbol{\eta}}M}\right\|_{\infty}
		\left\|{\partial^{\boldsymbol{\eta}-\boldsymbol{\ell}}\boldsymbol{u}}\right\|_{{\mathcal H}}
.
	\end{align*}}
\end{lemma}
\begin{proof}
For bound of $\boldsymbol{w}^{*}$, taking the $\boldsymbol{e}$-th derivative both sides of \eqref{variational form 1} and collect the terms with $\partial^{\boldsymbol{e}} \boldsymbol{u}$ on the left-hand side to obtain 
		\begin{align}\label{ve derivative}
		&a(\partial^{\boldsymbol{e}} \boldsymbol{u}, \boldsymbol{v};A) 
		+
		m(\boldsymbol{u},\partial^{\boldsymbol{e}}\boldsymbol{u},\boldsymbol{v};M)
		+m(\partial^{\boldsymbol{e}}\boldsymbol{u},\boldsymbol{u} ,\boldsymbol{v};M)
-
b( \boldsymbol{v},\partial^{\boldsymbol{e}}p ;B)\notag\\
		&\qquad=
		\left\langle{\partial^{\boldsymbol{e}}\boldsymbol{f}},{\boldsymbol{v}}
	\right\rangle
+
b( \boldsymbol{v} , p;\partial^{\boldsymbol{e}} B)		
-
		a(\boldsymbol{u},\boldsymbol{v};\partial^{\boldsymbol{e}}A)
		-
		m(\boldsymbol{u},\boldsymbol{u},\boldsymbol{v};\partial^{\boldsymbol{e}}M),
\qquad 
\forall\boldsymbol{v}\in{\mathcal H} 
.
	\end{align}
We now take the $\boldsymbol{\nu}$-derivative in both sides and keep the highest order of $\boldsymbol{u}$ and $p$ on the left-hand side to obtain
\begin{align}\label{var form vnu derivative}
	&a( \partial^{\boldsymbol{\nu}+\boldsymbol{e}} \boldsymbol{u}, \boldsymbol{v};A) 
		+
		m(\boldsymbol{u},\partial^{\boldsymbol{\nu}+\boldsymbol{e}}\boldsymbol{u},\boldsymbol{v};M)
		+m(\partial^{\boldsymbol{\nu}+\boldsymbol{e}}\boldsymbol{u},\boldsymbol{u} ,\boldsymbol{v};M)
-
b( \boldsymbol{v},\partial^{\boldsymbol{\nu}+\boldsymbol{e}} p;B)\notag\\
&\quad=
\left\langle{\widehat{f}_2(\boldsymbol{\nu})},{\boldsymbol{v}}
	\right\rangle
+
\left\langle{\widehat{b}_2(\boldsymbol{\nu})},{\boldsymbol{v}}
	\right\rangle
-
\left\langle{\widehat{a}_2(\boldsymbol{\nu})},{\boldsymbol{v}}
	\right\rangle
-
\left\langle{\widehat{m}_2(\boldsymbol{\nu})},{\boldsymbol{v}}
	\right\rangle,
\qquad\qquad\qquad \forall\boldsymbol{v}\in{\mathcal H},
\end{align}
where
\begin{align*}
\left\langle{\widehat{f}_2(\boldsymbol{\nu})},{\boldsymbol{v}}
	\right\rangle
&:=
\left\langle{\partial^{\boldsymbol{\nu}+\boldsymbol{e}}\boldsymbol{f}},{\boldsymbol{v}}
	\right\rangle
-
a(\boldsymbol{u} , \boldsymbol{v};\partial^{\boldsymbol{\nu}+\boldsymbol{e}}A )
+
b( \boldsymbol{v} ,p;\partial^{\boldsymbol{\nu}+\boldsymbol{e}}B )
-
m(\boldsymbol{u},\boldsymbol{u} ,\boldsymbol{v};\partial^{\boldsymbol{\nu}+\boldsymbol{e}}M),
\\
\left\langle{\widehat{a}_2(\boldsymbol{\nu})},{\boldsymbol{v}}
	\right\rangle
&:=
\sum_{\boldsymbol{0}< \boldsymbol{\eta}\leq \boldsymbol{\nu}}
	\left(\boldsymbol{\nu}\atop \boldsymbol{\eta}\right)
\left(a( \partial^{\boldsymbol{\eta}}  \boldsymbol{u} , \boldsymbol{v};\partial^{\boldsymbol{\nu}+\boldsymbol{e}-\boldsymbol{\eta}}A  )
+
a( \partial^{\boldsymbol{\nu}+\boldsymbol{e}-\boldsymbol{\eta}}  \boldsymbol{u} , \boldsymbol{v} ;\partial^{\boldsymbol{\eta}}A )\right),
\\
\left\langle{\widehat{b}_2(\boldsymbol{\nu})},{\boldsymbol{v}}
	\right\rangle
&:=
\sum_{\boldsymbol{0}< \boldsymbol{\eta}\leq \boldsymbol{\nu}}
	\left(\boldsymbol{\nu}\atop \boldsymbol{\eta}\right)
\left(b( \boldsymbol{v} ,\partial^{\boldsymbol{\eta}}p;\partial^{\boldsymbol{\nu}+\boldsymbol{e}-\boldsymbol{\eta}}B )
+
b( \boldsymbol{v} , \partial^{\boldsymbol{\nu}+\boldsymbol{e}-\boldsymbol{\eta}}p;\partial^{\boldsymbol{\eta}} B )\right),\\
\left\langle{\widehat{m}_2(\boldsymbol{\nu})},{\boldsymbol{v}}
	\right\rangle
&:=
\sum_{\boldsymbol{0}< \boldsymbol{\eta}\leq \boldsymbol{\nu}}
\sum_{\boldsymbol{0}\leq \boldsymbol{\ell}\leq \boldsymbol{\eta}}
	\left(\boldsymbol{\nu}\atop \boldsymbol{\eta}\right)
\left(\boldsymbol{\eta}\atop \boldsymbol{\ell}\right)
(
m(\partial^{\boldsymbol{\ell}}\boldsymbol{u}, \partial^{\boldsymbol{\nu}+\boldsymbol{e}-\boldsymbol{\eta}}\boldsymbol{u} ,\boldsymbol{v};\partial^{\boldsymbol{\eta}-\boldsymbol{\ell}}M)
\notag\\
&\qquad +
	m(\partial^{\boldsymbol{\nu}+\boldsymbol{e}-\boldsymbol{\eta}}\boldsymbol{u},\partial^{\boldsymbol{\ell}}\boldsymbol{u} ,\boldsymbol{v};\partial^{\boldsymbol{\eta}-\boldsymbol{\ell}}M)
+
	m(\partial^{\boldsymbol{\eta}-\boldsymbol{\ell}}\boldsymbol{u},\partial^{\boldsymbol{\ell}}\boldsymbol{u} ,\boldsymbol{v};\partial^{\boldsymbol{\nu}+\boldsymbol{e}-\boldsymbol{\eta}}M)
).
\end{align*}
Notice \eqref{a operator bound},\eqref{b operator bound} and \eqref{m operator bound}, we estimate the ${\mathcal H}^*$-norm of $\widehat{f}_2(\boldsymbol{\nu}),\widehat{a}_2(\boldsymbol{\nu}),\widehat{b}_2(\boldsymbol{\nu})$ and $\widehat{m}_2(\boldsymbol{\nu})$ as following
\begin{align}
\left\|{\widehat{f}_2(\boldsymbol{\nu})}\right\|_{{\mathcal H}^*}
&\leq
\left\|{\partial^{\boldsymbol{\nu}+\boldsymbol{e}}\boldsymbol{f}}\right\|_{{\mathcal H}^*}
+
\overline{u}
\left\|{\partial^{\boldsymbol{\nu}+\boldsymbol{e}}A }\right\|_{\infty}
+
\left\|{\partial^{\boldsymbol{\nu}+\boldsymbol{e}}B}\right\|_{\infty}
\left\|{p}\right\|_{\mathcal L}
+
{C_4^2} \overline{u}^2
\left\|{\partial^{\boldsymbol{\nu}+\boldsymbol{e}}M}\right\|_{\infty}
=
\mathfrak{F}_2(\boldsymbol{\nu}) \label{f2 bnd},
\\
\left\|{\widehat{a}_2(\boldsymbol{\nu})}\right\|_{{\mathcal H}^*}
&\leq
\sum_{\boldsymbol{0}< \boldsymbol{\eta}\leq \boldsymbol{\nu}}
\left(\boldsymbol{\nu}\atop \boldsymbol{\eta}\right)
\left(
\left\|{\partial^{\boldsymbol{\nu}+\boldsymbol{e}-\boldsymbol{\eta}}A }\right\|_{\infty}
\left\|{\partial^{\boldsymbol{\eta}}  \boldsymbol{u}}\right\|_{{\mathcal H}}
+
\left\|{\partial^{\boldsymbol{\eta}}A }\right\|_{\infty}
\left\|{\partial^{\boldsymbol{\nu}+\boldsymbol{e}-\boldsymbol{\eta}} \boldsymbol{u}}\right\|_{{\mathcal H}}\right)
=
\mathfrak{A}(\boldsymbol{\nu}),\label{a2 bnd}\\
\left\|{\widehat{b}_2(\boldsymbol{\nu})}\right\|_{{\mathcal H}^*}
&\leq
\sum_{\boldsymbol{0}< \boldsymbol{\eta}\leq \boldsymbol{\nu}}
\left(\boldsymbol{\nu}\atop \boldsymbol{\eta}\right)
\left(
\left\|{\partial^{\boldsymbol{\nu}+\boldsymbol{e}-\boldsymbol{\eta}}B}\right\|_{\infty}
\left\|{\partial^{\boldsymbol{\eta}}p }\right\|_{\mathcal L}
+
\left\|{\partial^{\boldsymbol{\eta}}B}\right\|_{\infty}
\left\|{ \partial^{\boldsymbol{\nu}+\boldsymbol{e}-\boldsymbol{\eta}}p }\right\|_{\mathcal L}\right)
=
\mathfrak{B}_2(\boldsymbol{\nu})\label{b2 bnd},\\
\left\|{\widehat{m}_2(\boldsymbol{\nu})}\right\|_{{\mathcal H}^*}
&\leq
\sum_{\boldsymbol{0}< \boldsymbol{\eta}\leq \boldsymbol{\nu}}
\sum_{\boldsymbol{0}\leq \boldsymbol{\ell}\leq \boldsymbol{\eta}}
\left(\boldsymbol{\nu}\atop \boldsymbol{\eta}\right)
\left(\boldsymbol{\eta}\atop \boldsymbol{\ell}\right)
\Big(
2\,{C_4^2}\left\|{\partial^{\boldsymbol{\eta}-\boldsymbol{\ell}} M}\right\|_{\infty}
\left\|{\partial^{\boldsymbol{\nu}+\boldsymbol{e}-\boldsymbol{\eta}}\boldsymbol{u}}\right\|_{{\mathcal H}}
\left\|{\partial^{\boldsymbol{\ell}}\boldsymbol{u}}\right\|_{{\mathcal H}}\notag
\\&\qquad\qquad
+
{C_4^2}\left\|{\partial^{\boldsymbol{\nu}+\boldsymbol{e}-\boldsymbol{\eta}}M}\right\|_{\infty}
\left\|{\partial^{\boldsymbol{\eta}-\boldsymbol{\ell}}\boldsymbol{u}}\right\|_{{\mathcal H}}
\left\|{\partial^{\boldsymbol{\ell}}\boldsymbol{u}}\right\|_{{\mathcal H}}
\Big)
=
		\mathfrak{M}(\boldsymbol{\nu})
+
\mathfrak{N}(\boldsymbol{\nu})\label{m2 bnd}.
\end{align}
In the last equation, we split the second sum into two parts: $\mathfrak{M}(\boldsymbol{\nu})$ collects the terms corresponding to the multi-indices $\boldsymbol{\ell} = \boldsymbol{0}$ and $\boldsymbol{\ell} = \boldsymbol{\eta}$, while $\mathfrak{N}(\boldsymbol{\nu})$ contains the remaining terms, where $\boldsymbol{0} < \boldsymbol{\ell}<\boldsymbol{\eta}$.

Due to $\boldsymbol{w}^*\in \mathcal V$, we have $a( \partial^{\boldsymbol{\nu}+\boldsymbol{e}} \boldsymbol{u}, \boldsymbol{w}^*;A) = a( \boldsymbol{w}^* , \boldsymbol{w}^*;A) $ and $b( \boldsymbol{w}^* ,\partial^{\boldsymbol{\nu}+\boldsymbol{e}} p;B)=0$. We substitute $\boldsymbol{v}=\boldsymbol{w}^*$ into \eqref{var form vnu derivative} and notice the definition of $t_{\boldsymbol{y}}$ in \eqref{Atilde-def} to derive
\begin{align*}
t_{\boldsymbol{y}}(\boldsymbol{u},\boldsymbol{w}^*,\boldsymbol{w}^*)
&=
\left\langle{\widehat{f}_2(\boldsymbol{\nu})},{\boldsymbol{w}^*}
	\right\rangle
+
\left\langle{\widehat{b}_2(\boldsymbol{\nu})},{\boldsymbol{w}^*}
	\right\rangle
-
\left\langle{\widehat{a}_2(\boldsymbol{\nu})},{\boldsymbol{w}^*}
	\right\rangle
-
\left\langle{\widehat{m}_2(\boldsymbol{\nu})},{\boldsymbol{w}^*}
	\right\rangle\\
&\qquad-
		m(\boldsymbol{u}, \boldsymbol{w}^{\bot},\boldsymbol{w}^*;M)
-
		m(\boldsymbol{w}^{\bot}, \boldsymbol{u},\boldsymbol{w}^*;M)\\
&\leq
\left(
\mathfrak{F}_2(\boldsymbol{\nu})
+
\mathfrak{A}(\boldsymbol{\nu})
+
\mathfrak{B}_2(\boldsymbol{\nu})
+
		\mathfrak{M}(\boldsymbol{\nu})
+
\mathfrak{N}(\boldsymbol{\nu})
+
\overline{m}\left\|{\boldsymbol{u}}\right\|_{{\mathcal H}}\left\|{\boldsymbol{w}^{\bot}}\right\|_{{\mathcal H}}
\right)\left\|{\boldsymbol{w}^*}\right\|_{{\mathcal H}}.
\end{align*}
Notice the coercivity of $t_{\boldsymbol{y}}$ in \eqref{coercive lower bound}, we cancel $\left\|{\boldsymbol{w}^*}\right\|_{{\mathcal H}}$ in both sides to obtain \eqref{w* norm upper} and finish the proof.
\end{proof} \begin{lemma}\label{lem:p-bnd}
	For sufficiently regular solutions of \eqref{decomposition variational}, and recalling the definitions of $\mathfrak{F}_2(\boldsymbol{\nu})$, ${\mathfrak{A}}(\boldsymbol{\nu})$, $\mathfrak{B}_2(\boldsymbol{\nu})$, $\mathfrak{M}(\boldsymbol{\nu})$ and $\mathfrak{N}(\boldsymbol{\nu})$ from Lemma \ref{lem:w-bnd}. The following holds:
	\begin{align}
\label{p norm upper}
			\beta\left\|{\partial^{\boldsymbol{\nu}+\boldsymbol{e}} p}\right\|_{\mathcal L}
\leq
\mathfrak{F}_2(\boldsymbol{\nu})
+
\mathfrak{A}(\boldsymbol{\nu})
+
\mathfrak{B}_2(\boldsymbol{\nu})
		+
		\mathfrak{M}(\boldsymbol{\nu})
+
\mathfrak{N}(\boldsymbol{\nu})
+
\left(
\frac{\overline{a}}{2}
+
\overline{u}\,
\overline{m}\right)
\left\|{\partial^{\boldsymbol{\nu}+\boldsymbol{e}}\boldsymbol{u}}\right\|_{{\mathcal H}}.
	\end{align}
\end{lemma}
\begin{proof}
Taking the $\boldsymbol{\nu}$-derivative both side of \eqref{ve derivative} and collecting the term $b\left( \boldsymbol{v} ,\partial^{\boldsymbol{\nu}+\boldsymbol{e}} p;B\right)$ on the left-hand side, we obtain
\begin{align*}
	b\left( \boldsymbol{v},\partial^{\boldsymbol{\nu}+\boldsymbol{e}} p;B\right)
=	
	-
	\left\langle{\widehat{f}_3(\boldsymbol{\nu})},{\boldsymbol{v}}
	\right\rangle	
	+
	\left\langle{\widehat{a}_3(\boldsymbol{\nu})},{\boldsymbol{v}}
	\right\rangle
		-
	\left\langle{\widehat{b}_3(\boldsymbol{\nu})},{\boldsymbol{v}}
	\right\rangle
	+
	\left\langle{\widehat{m}_3(\boldsymbol{\nu})},{\boldsymbol{v}}
	\right\rangle
\end{align*}
where
\begin{align*}
	\left\langle{\widehat{f}_3(\boldsymbol{\nu})},{\boldsymbol{v}}
	\right\rangle
&:=
	\left\langle{\partial^{\boldsymbol{\nu}+\boldsymbol{e}}\boldsymbol{f}},{\boldsymbol{v}}
	\right\rangle
-
	a(  \boldsymbol{u} , \boldsymbol{v};\partial^{\boldsymbol{\nu}+\boldsymbol{e}}A )
+
	b( \boldsymbol{v} ,p ;\partial^{\boldsymbol{\nu}+\boldsymbol{e}}B)
-
m(\boldsymbol{u},\boldsymbol{u} ,\boldsymbol{v};\partial^{\boldsymbol{\nu}+\boldsymbol{e}}M)
\\
	&\quad-
a(  \partial^{\boldsymbol{\nu}+\boldsymbol{e}} \boldsymbol{u} , \boldsymbol{v} ;A)
-
m(\boldsymbol{u}, \partial^{\boldsymbol{\nu}+\boldsymbol{e}} \boldsymbol{u} ,\boldsymbol{v};M)
-
m(\partial^{\boldsymbol{\nu}+\boldsymbol{e}} \boldsymbol{u},\boldsymbol{u} ,\boldsymbol{v};M)
\\
&\,=
\left\langle{\widehat{f}_2(\boldsymbol{\nu})},{\boldsymbol{v}}
	\right\rangle
-
a(\partial^{\boldsymbol{\nu}+\boldsymbol{e}} \boldsymbol{u} , \boldsymbol{v} ;A)
-
m(\boldsymbol{u}, \partial^{\boldsymbol{\nu}+\boldsymbol{e}} \boldsymbol{u} ,\boldsymbol{v};M)
-
m(\partial^{\boldsymbol{\nu}+\boldsymbol{e}}\boldsymbol{u},\boldsymbol{u} ,\boldsymbol{v};M),
\\
	\left\langle{\widehat{a}_3(\boldsymbol{\nu})},{\boldsymbol{v}}
	\right\rangle
	&:=
	\sum_{\boldsymbol{0}< \boldsymbol{\eta}\leq \boldsymbol{\nu}}
	\left(\boldsymbol{\nu}\atop \boldsymbol{\eta}\right)
	\left(
a(  \partial^{\boldsymbol{\eta}}  \boldsymbol{u} , \boldsymbol{v};\partial^{\boldsymbol{\nu}+\boldsymbol{e}-\boldsymbol{\eta}}A )
	+
	a( \partial^{\boldsymbol{\nu}+\boldsymbol{e}-\boldsymbol{\eta}}  \boldsymbol{u} , \boldsymbol{v} ;\partial^{\boldsymbol{\eta}}A )\right)
=
\left\langle{\widehat{a}_2(\boldsymbol{\nu})},{\boldsymbol{v}}
	\right\rangle,
	\\
	\left\langle{\widehat{b}_3(\boldsymbol{\nu})},{\boldsymbol{v}}
	\right\rangle
	&:=
	\sum_{\boldsymbol{0}<  \boldsymbol{\eta}\leq \boldsymbol{\nu}}
	\left(\boldsymbol{\nu}\atop \boldsymbol{\eta}\right)
\left(	b( \boldsymbol{v} ,\partial^{\boldsymbol{\eta}}p;\partial^{\boldsymbol{\nu}+\boldsymbol{e}-\boldsymbol{\eta}}B )
	+
	b( \boldsymbol{v}, \partial^{\boldsymbol{\nu}+\boldsymbol{e}-\boldsymbol{\eta}}p; \partial^{\boldsymbol{\eta}} B )\right)
=
\left\langle{\widehat{b}_2(\boldsymbol{\nu})},{\boldsymbol{v}}
	\right\rangle,\\
	\left\langle{\widehat{m}_3(\boldsymbol{\nu})},{\boldsymbol{v}}
	\right\rangle
	&:=
\sum_{\boldsymbol{0}< \boldsymbol{\eta}\leq \boldsymbol{\nu}}
\sum_{\boldsymbol{0}\leq \boldsymbol{\ell}\leq \boldsymbol{\eta}}
	\left(\boldsymbol{\nu}\atop \boldsymbol{\eta}\right)
\left(\boldsymbol{\eta}\atop \boldsymbol{\ell}\right)
(
m(\partial^{\boldsymbol{\ell}}\boldsymbol{u}, \partial^{\boldsymbol{\nu}+\boldsymbol{e}-\boldsymbol{\eta}}\boldsymbol{u} ,\boldsymbol{v};\partial^{\boldsymbol{\eta}-\boldsymbol{\ell}}M)
\notag\\
& +
	m(\partial^{\boldsymbol{\nu}+\boldsymbol{e}-\boldsymbol{\eta}}\boldsymbol{u},\partial^{\boldsymbol{\ell}}\boldsymbol{u} ,\boldsymbol{v};\partial^{\boldsymbol{\eta}-\boldsymbol{\ell}}M)
+
	m(\partial^{\boldsymbol{\eta}-\boldsymbol{\ell}}\boldsymbol{u},\partial^{\boldsymbol{\ell}}\boldsymbol{u} ,\boldsymbol{v};\partial^{\boldsymbol{\nu}+\boldsymbol{e}-\boldsymbol{\eta}}M)
)
=
\left\langle{\widehat{m}_2(\boldsymbol{\nu})},{\boldsymbol{v}}
	\right\rangle.
\end{align*}
Observe that $\widehat{a}_3(\boldsymbol{\nu})=\widehat{a}_2(\boldsymbol{\nu}),\widehat{b}_3(\boldsymbol{\nu})=\widehat{b}_2(\boldsymbol{\nu})$ and $\widehat{m}_3(\boldsymbol{\nu})=\widehat{m}_2(\boldsymbol{\nu})$, we now only need to estimate ${\mathcal H}^*$-norm of $\widehat{f}_3(\boldsymbol{\nu})$. Notice \eqref{f2 bnd},\eqref{a2 bnd},\eqref{b2 bnd} and \eqref{m2 bnd}, we get
\begin{align*}
\left\|{\widehat{a}_3(\boldsymbol{\nu})}\right\|_{{\mathcal H}^*}\leq
\mathfrak{A}(\boldsymbol{\nu}),\quad
\left\|{\widehat{b}_3(\boldsymbol{\nu})}\right\|_{{\mathcal H}^*}\leq
\mathfrak{B}_2(\boldsymbol{\nu}),\quad
\left\|{\widehat{m}_3(\boldsymbol{\nu})}\right\|_{{\mathcal H}^*}\leq
\mathfrak{M}(\boldsymbol{\nu})
+
\mathfrak{N}(\boldsymbol{\nu})
\end{align*}
and
\begin{align*}
\left\|{\widehat{f}_3(\boldsymbol{\nu})}\right\|_{{\mathcal H}^*}
&\leq
\left\|{\widehat{f}_2(\boldsymbol{\nu})}\right\|_{{\mathcal H}^*}
+
\frac{\overline{a}}{2}
\left\|{ \partial^{\boldsymbol{\nu}+\boldsymbol{e}} \boldsymbol{u}}\right\|_{{\mathcal H}}
+
\overline{m}
\left\|{\partial^{\boldsymbol{\nu}+\boldsymbol{e}}\boldsymbol{u}}\right\|_{{\mathcal H}}
\left\|{\boldsymbol{u}}\right\|_{{\mathcal H}}\\
&=
\mathfrak{F}_2(\boldsymbol{\nu})
+
\left(
\frac{\overline{a}}{2}
+
\overline{u}\,
\overline{m}\right)
\left\|{\partial^{\boldsymbol{\nu}+\boldsymbol{e}}\boldsymbol{u}}\right\|_{{\mathcal H}}.
\end{align*}
This together with \eqref{BT isomorphism} shows \eqref{p norm upper}. We finish the proof.
\end{proof}
 Before proving Theorem \ref{gevrey regularity for NS}, we need to bound all the derivative terms appearing in the lemmas above. The following lemma provides explicit bounds for these terms.

\begin{lemma}\label{lem:F A B M bound}
	For sufficiently regular solutions of \eqref{decomposition variational}, and for any multi-index $\boldsymbol{\nu} \in \mathcal F\setminus\left\{\boldsymbol{0}\right\}$ and $\delta \geq 1$, we assume the existence of constants $C_u>0$ and $\rho\geq 1$ such that for all $\boldsymbol{\eta} \leq \boldsymbol{\nu}$  and $\left|\boldsymbol{\eta}\right|\geq 1$ 
\begin{align}\label{deri u assum}
		\left\|{\partial^{\boldsymbol{\eta}} \boldsymbol{u}( \boldsymbol{y})}\right\|_{{\mathcal H}}
		&\leq
		\frac{C_u \rho^{\left|\boldsymbol{\eta}\right|-1} \left[\tfrac{1}{2} \right]_{\left|{\boldsymbol{\eta}}\right|} }{\boldsymbol{R}^{{\boldsymbol{\eta}}}} (\left|\boldsymbol{\eta}\right|!)^{\delta-1},
\end{align}
and
\begin{align}\label{deri p assum}
		\left\|{\partial^{\boldsymbol{\eta}} p( \boldsymbol{y})}\right\|_{{\mathcal L}}
		&\leq
		\frac{C_p\, \rho^{\left|\boldsymbol{\eta}\right|-1} \left[\tfrac{1}{2} \right]_{\left|{\boldsymbol{\eta}}\right|} }{\boldsymbol{R}^{{\boldsymbol{\eta}}}} (\left|\boldsymbol{\eta}\right|!)^{\delta-1}.
\end{align}
We make use of the definitions of $\mathfrak{F}_1(\boldsymbol{\nu})$ and $\mathfrak{B}_1(\boldsymbol{\nu})$ given in Lemma \ref{lem:wT-bnd}, along with those of $\mathfrak{F}_2(\boldsymbol{\nu})$, ${\mathfrak{A}}(\boldsymbol{\nu})$, $\mathfrak{B}_2(\boldsymbol{\nu})$, $\mathfrak{M}(\boldsymbol{\nu})$ and $\mathfrak{N}(\boldsymbol{\nu})$ from Lemma \ref{lem:w-bnd} . Then, the following holds
\begin{align}\label{F_1 bound}
		\mathfrak{F}_1(\boldsymbol{\nu})
&\leq
\left(\overline{g}+ \overline{b}\, \overline{u}\right) \frac{ \left[\tfrac{1}{2} \right]_{\left|\boldsymbol{\nu}+\boldsymbol{e}\right|} }{\boldsymbol{R}^{\boldsymbol{\nu}+\boldsymbol{e}}(\left|\boldsymbol{\nu}+\boldsymbol{e}\right|!)^{1-\delta}},\\
	\mathfrak{F}_2(\boldsymbol{\nu})
&\leq
\left(\overline{f}+\overline{b}\,\overline{p}+
\overline{u}
\left(
\overline{a}+\overline{m}\,\overline{u}
\right)\right)
\frac{ \left[\tfrac{1}{2} \right]_{\left|\boldsymbol{\nu}+\boldsymbol{e}\right|} }{\boldsymbol{R}^{\boldsymbol{\nu}+\boldsymbol{e}}(\left|\boldsymbol{\nu}+\boldsymbol{e}\right|!)^{1-\delta}}
\label{F_2 bound},\\	
\mathfrak{A}(\boldsymbol{\nu})
		&\leq
		2C_u \overline{a}
\frac{\rho^{\left|\boldsymbol{\nu}\right|-1} \left[\tfrac{1}{2} \right]_{\left|\boldsymbol{\nu}+\boldsymbol{e}\right|} }{\boldsymbol{R}^{\boldsymbol{\nu}+\boldsymbol{e}}(\left|\boldsymbol{\nu}+\boldsymbol{e}\right|!)^{1-\delta}}\label{A bound},\\
		\mathfrak{B}_1(\boldsymbol{\nu})
		&\leq
2C_u \overline{b}
\frac{\rho^{\left|\boldsymbol{\nu}\right|-1} \left[\tfrac{1}{2} \right]_{\left|\boldsymbol{\nu}+\boldsymbol{e}\right|} }{\boldsymbol{R}^{\boldsymbol{\nu}+\boldsymbol{e}}(\left|\boldsymbol{\nu}+\boldsymbol{e}\right|!)^{1-\delta}}\label{B_1 bound},
\\
		\mathfrak{B}_2(\boldsymbol{\nu})
		&\leq
		2C_p \overline{b}
\frac{\rho^{\left|\boldsymbol{\nu}\right|-1} \left[\tfrac{1}{2} \right]_{\left|\boldsymbol{\nu}+\boldsymbol{e}\right|} }{\boldsymbol{R}^{\boldsymbol{\nu}+\boldsymbol{e}}(\left|\boldsymbol{\nu}+\boldsymbol{e}\right|!)^{1-\delta}}\label{B_2 bound},\\
		\mathfrak{M}(\boldsymbol{\nu})
		&\leq
\left(4\overline{m}\,\overline{u}C_u +\overline{m}C_u^2\right)
\frac{\rho^{\left|\boldsymbol{\nu}\right|-1} \left[\tfrac{1}{2} \right]_{\left|\boldsymbol{\nu}+\boldsymbol{e}\right|} }{\boldsymbol{R}^{\boldsymbol{\nu}+\boldsymbol{e}}(\left|\boldsymbol{\nu}+\boldsymbol{e}\right|!)^{1-\delta}}\label{M bound},\\
		\mathfrak{N}(\boldsymbol{\nu})
		&\leq
			6\overline{m}C_u^2
\frac{\rho^{\left|\boldsymbol{\nu}\right|-1} \left[\tfrac{1}{2} \right]_{\left|\boldsymbol{\nu}+\boldsymbol{e}\right|} }{\boldsymbol{R}^{\boldsymbol{\nu}+\boldsymbol{e}}(\left|\boldsymbol{\nu}+\boldsymbol{e}\right|!)^{1-\delta}}\label{N bound}.
\end{align}
\end{lemma}
\begin{proof}
To establish the bounds above, we make use of the definitions of $\mathfrak{F}_1(\boldsymbol{\nu})$ and $\mathfrak{B}_1(\boldsymbol{\nu})$ from Lemma \ref{lem:wT-bnd}, as well as $\mathfrak{F}_2(\boldsymbol{\nu})$, $\mathfrak{A}(\boldsymbol{\nu})$, $\mathfrak{B}_2(\boldsymbol{\nu})$, $\mathfrak{M}(\boldsymbol{\nu})$, and $\mathfrak{N}(\boldsymbol{\nu})$ from Lemma \ref{lem:w-bnd}, together with the derivative bounds for $A$, $B$, $M$, $\boldsymbol{f}$, and $g$ given in \eqref{abf assumption}. We obtain the bound for $\mathfrak{F}_1(\boldsymbol{\nu})$ as follows
\begin{align*}
\mathfrak{F}_1(\boldsymbol{\nu})
&\leq
\frac{\overline{g}\left[\tfrac{1}{2} \right]_{\left|\boldsymbol{\nu}+\boldsymbol{e}\right|} (\left|\boldsymbol{\nu}+\boldsymbol{e}\right|!)^{\delta-1}}{\boldsymbol{R}^{\boldsymbol{\nu}+\boldsymbol{e}}}
+
\frac{\overline{b}\, \overline{u} \left[\tfrac{1}{2} \right]_{\left|\boldsymbol{\nu}+\boldsymbol{e}\right|} (\left|\boldsymbol{\nu}+\boldsymbol{e}\right|!)^{\delta-1}}{\boldsymbol{R}^{\boldsymbol{\nu}+\boldsymbol{e}}}\\
&\leq
\left(\overline{g}+ \overline{b}\, \overline{u}\right) \frac{ \left[\tfrac{1}{2} \right]_{\left|\boldsymbol{\nu}+\boldsymbol{e}\right|} (\left|\boldsymbol{\nu}+\boldsymbol{e}\right|!)^{\delta-1}}{\boldsymbol{R}^{\boldsymbol{\nu}+\boldsymbol{e}}}.
\end{align*} 
It shows \eqref{F_1 bound}. 
For the bounds of $\mathfrak{F}_2(\boldsymbol{\nu})$, we have
\begin{align*}
\mathfrak{F}_2(\boldsymbol{\nu})
&\leq
\frac{\overline{f}\left[\tfrac{1}{2} \right]_{\left|\boldsymbol{\nu}+\boldsymbol{e}\right|} (\left|\boldsymbol{\nu}+\boldsymbol{e}\right|!)^{\delta-1}}{\boldsymbol{R}^{\boldsymbol{\nu}+\boldsymbol{e}}}
+
\overline{p}\,
\frac{\overline{b} \left[\tfrac{1}{2} \right]_{\left|\boldsymbol{\nu}+\boldsymbol{e}\right|} (\left|\boldsymbol{\nu}+\boldsymbol{e}\right|!)^{\delta-1}}{\boldsymbol{R}^{\boldsymbol{\nu}+\boldsymbol{e}}}\\
&\quad+
\overline{u}\,
\frac{\overline{a} \left[\tfrac{1}{2} \right]_{\left|\boldsymbol{\nu}+\boldsymbol{e}\right|} (\left|\boldsymbol{\nu}+\boldsymbol{e}\right|!)^{\delta-1}}{\boldsymbol{R}^{\boldsymbol{\nu}+\boldsymbol{e}}}
+
\overline{u}^2
\frac{\overline{m} \left[\tfrac{1}{2} \right]_{\left|\boldsymbol{\nu}+\boldsymbol{e}\right|} (\left|\boldsymbol{\nu}+\boldsymbol{e}\right|!)^{\delta-1}}{\boldsymbol{R}^{\boldsymbol{\nu}+\boldsymbol{e}}}\\
&\leq
\left(\overline{f}+\overline{b}\,\overline{p}+
\overline{u}
\left(
\overline{a}+\overline{m}\,\overline{u}
\right)\right)
\frac{ \left[\tfrac{1}{2} \right]_{\left|\boldsymbol{\nu}+\boldsymbol{e}\right|} }{\boldsymbol{R}^{\boldsymbol{\nu}+\boldsymbol{e}}(\left|\boldsymbol{\nu}+\boldsymbol{e}\right|!)^{1-\delta}}.
\end{align*} 
This is precisely \eqref{F_2 bound}. For the bound of $\mathfrak{A}(\boldsymbol{\nu})$, notice  $(\left|\boldsymbol{\nu}+\boldsymbol{e}-\boldsymbol{\eta}\right|!\left|\boldsymbol{\eta}\right|!)^{\delta-1}\leq (\left|\boldsymbol{\nu}+\boldsymbol{e}\right|!)^{\delta-1}$ and \eqref{multiindex-est-6}, we obtain
\begin{align*}
\mathfrak{A}(\boldsymbol{\nu})
&\leq
(\left|\boldsymbol{\nu}+\boldsymbol{e}\right|!)^{\delta-1}
\sum_{0<\boldsymbol{\eta}\leq \boldsymbol{\nu}}
\left(\boldsymbol{\nu}\atop \boldsymbol{\eta}\right)
\frac{\overline{a}\left[\tfrac{1}{2} \right]_{\left|\boldsymbol{\nu}+\boldsymbol{e}-\boldsymbol{\eta}\right|} }{\boldsymbol{R}^{\boldsymbol{\nu}+\boldsymbol{e}-\boldsymbol{\eta}}}
\,
\frac{C_u \rho^{\left|\boldsymbol{\eta}\right|-1}\left[\tfrac{1}{2} \right]_{\left|\boldsymbol{\eta}\right|} }{\boldsymbol{R}^{\boldsymbol{\eta}}}\\
&\qquad+
(\left|\boldsymbol{\nu}+\boldsymbol{e}\right|!)^{\delta-1}
\sum_{0<\boldsymbol{\eta}\leq \boldsymbol{\nu}}
\left(\boldsymbol{\nu}\atop \boldsymbol{\eta}\right)
\frac{\overline{a}\left[\tfrac{1}{2} \right]_{\left|\boldsymbol{\eta}\right|} }{\boldsymbol{R}^{\boldsymbol{\eta}}}
\,
\frac{C_u \rho^{\left|\boldsymbol{\nu}+\boldsymbol{e}-\boldsymbol{\eta}\right|-1}\left[\tfrac{1}{2} \right]_{\left|\boldsymbol{\nu}+\boldsymbol{e}-\boldsymbol{\eta}\right|} }{\boldsymbol{R}^{\boldsymbol{\nu}+\boldsymbol{e}-\boldsymbol{\eta}}}\\
&\leq
\frac{2C_u \overline{a} \rho^{\left|\boldsymbol{\nu}\right|-1}}{\boldsymbol{R}^{\boldsymbol{\nu}+\boldsymbol{e}}}
(\left|\boldsymbol{\nu}+\boldsymbol{e}\right|!)^{\delta-1}
\sum_{0<\boldsymbol{\eta}\leq \boldsymbol{\nu}}
\left(\boldsymbol{\nu}\atop \boldsymbol{\eta}\right)
\left[\tfrac{1}{2} \right]_{\left|\boldsymbol{\nu}+\boldsymbol{e}-\boldsymbol{\eta}\right|} 
\left[\tfrac{1}{2} \right]_{\left|\boldsymbol{\eta}\right|} 
\leq
2C_u \overline{a}
\frac{\rho^{\left|\boldsymbol{\nu}\right|-1} \left[\tfrac{1}{2} \right]_{\left|\boldsymbol{\nu}+\boldsymbol{e}\right|} }{\boldsymbol{R}^{\boldsymbol{\nu}+\boldsymbol{e}}(\left|\boldsymbol{\nu}+\boldsymbol{e}\right|!)^{1-\delta}}.
\end{align*}
It shows \eqref{A bound}. To estimate $\mathfrak{B}_1(\boldsymbol{\nu})$ and $\mathfrak{B}_2(\boldsymbol{\nu})$, we also employ \eqref{multiindex-est-6} to obtain
\begin{align*}
\mathfrak{B}_1(\boldsymbol{\nu})
&\leq
(\left|\boldsymbol{\nu}+\boldsymbol{e}\right|!)^{\delta-1}
\sum_{0<\boldsymbol{\eta}\leq \boldsymbol{\nu}}
\left(\boldsymbol{\nu}\atop \boldsymbol{\eta}\right)
\frac{\overline{b}\left[\tfrac{1}{2} \right]_{\left|\boldsymbol{\nu}+\boldsymbol{e}-\boldsymbol{\eta}\right|} }{\boldsymbol{R}^{\boldsymbol{\nu}+\boldsymbol{e}-\boldsymbol{\eta}}}
\,
\frac{C_u \rho^{\left|\boldsymbol{\eta}\right|-1}\left[\tfrac{1}{2} \right]_{\left|\boldsymbol{\eta}\right|} }{\boldsymbol{R}^{\boldsymbol{\eta}}}\\
&\qquad+
(\left|\boldsymbol{\nu}+\boldsymbol{e}\right|!)^{\delta-1}
\sum_{0<\boldsymbol{\eta}\leq \boldsymbol{\nu}}
\left(\boldsymbol{\nu}\atop \boldsymbol{\eta}\right)
\frac{\overline{b}\left[\tfrac{1}{2} \right]_{\left|\boldsymbol{\eta}\right|} }{\boldsymbol{R}^{\boldsymbol{\eta}}}
\,
\frac{C_u \rho^{\left|\boldsymbol{\nu}+\boldsymbol{e}-\boldsymbol{\eta}\right|-1}\left[\tfrac{1}{2} \right]_{\left|\boldsymbol{\nu}+\boldsymbol{e}-\boldsymbol{\eta}\right|} }{\boldsymbol{R}^{\boldsymbol{\nu}+\boldsymbol{e}-\boldsymbol{\eta}}}\\
&=
\frac{2C_u \overline{b} \rho^{\left|\boldsymbol{\nu}\right|-1}}{\boldsymbol{R}^{\boldsymbol{\nu}+\boldsymbol{e}}}
(\left|\boldsymbol{\nu}+\boldsymbol{e}\right|!)^{\delta-1}
\sum_{0<\boldsymbol{\eta}\leq \boldsymbol{\nu}}
\left(\boldsymbol{\nu}\atop \boldsymbol{\eta}\right)
\left[\tfrac{1}{2} \right]_{\left|\boldsymbol{\nu}+\boldsymbol{e}-\boldsymbol{\eta}\right|} 
\left[\tfrac{1}{2} \right]_{\left|\boldsymbol{\eta}\right|} 
\leq
2C_u \overline{b}
\frac{\rho^{\left|\boldsymbol{\nu}\right|-1} \left[\tfrac{1}{2} \right]_{\left|\boldsymbol{\nu}+\boldsymbol{e}\right|} }{\boldsymbol{R}^{\boldsymbol{\nu}+\boldsymbol{e}}(\left|\boldsymbol{\nu}+\boldsymbol{e}\right|!)^{1-\delta}}.
\end{align*}
and 
{\small
\begin{align*}
\mathfrak{B}_2(\boldsymbol{\nu})
&\leq
(\left|\boldsymbol{\nu}+\boldsymbol{e}\right|!)^{\delta-1}
\sum_{0<\boldsymbol{\eta}\leq \boldsymbol{\nu}}
\left(\boldsymbol{\nu}\atop \boldsymbol{\eta}\right)
\frac{\overline{b}\left[\tfrac{1}{2} \right]_{\left|\boldsymbol{\nu}+\boldsymbol{e}-\boldsymbol{\eta}\right|} }{\boldsymbol{R}^{\boldsymbol{\nu}+\boldsymbol{e}-\boldsymbol{\eta}}}
\,
\frac{C_p \rho^{\left|\boldsymbol{\eta}\right|-1}\left[\tfrac{1}{2} \right]_{\left|\boldsymbol{\eta}\right|} }{\boldsymbol{R}^{\boldsymbol{\eta}}}\\
&\qquad+
(\left|\boldsymbol{\nu}+\boldsymbol{e}\right|!)^{\delta-1}
\sum_{0<\boldsymbol{\eta}\leq \boldsymbol{\nu}}
\left(\boldsymbol{\nu}\atop \boldsymbol{\eta}\right)
\frac{\overline{b}\left[\tfrac{1}{2} \right]_{\left|\boldsymbol{\eta}\right|} }{\boldsymbol{R}^{\boldsymbol{\eta}}}
\,
\frac{C_p \rho^{\left|\boldsymbol{\nu}+\boldsymbol{e}-\boldsymbol{\eta}\right|-1}\left[\tfrac{1}{2} \right]_{\left|\boldsymbol{\nu}+\boldsymbol{e}-\boldsymbol{\eta}\right|} }{\boldsymbol{R}^{\boldsymbol{\nu}+\boldsymbol{e}-\boldsymbol{\eta}}}\\
&=
\frac{2C_p \overline{b} \rho^{\left|\boldsymbol{\nu}\right|-1}}{\boldsymbol{R}^{\boldsymbol{\nu}+\boldsymbol{e}}}
(\left|\boldsymbol{\nu}+\boldsymbol{e}\right|!)^{\delta-1}
\sum_{0<\boldsymbol{\eta}\leq \boldsymbol{\nu}}
\left(\boldsymbol{\nu}\atop \boldsymbol{\eta}\right)
\left[\tfrac{1}{2} \right]_{\left|\boldsymbol{\nu}+\boldsymbol{e}-\boldsymbol{\eta}\right|} 
\left[\tfrac{1}{2} \right]_{\left|\boldsymbol{\eta}\right|} 
\leq
2C_p \overline{b}
\frac{\rho^{\left|\boldsymbol{\nu}\right|-1} \left[\tfrac{1}{2} \right]_{\left|\boldsymbol{\nu}+\boldsymbol{e}\right|} }{\boldsymbol{R}^{\boldsymbol{\nu}+\boldsymbol{e}}(\left|\boldsymbol{\nu}+\boldsymbol{e}\right|!)^{1-\delta}}.
\end{align*}
}It shows \eqref{B_1 bound} and \eqref{B_2 bound}. Similarly, notice \eqref{multiindex-est-6}, we derive the bound for $\mathfrak{M}(\boldsymbol{\nu})$  as follows
\begin{align*}
\mathfrak{M}(\boldsymbol{\nu})
&\leq
	2(\left|\boldsymbol{\nu}+\boldsymbol{e}\right|!)^{\delta-1}\overline{u}\sum_{\boldsymbol{0}< \boldsymbol{\eta}\leq \boldsymbol{\nu}}
		\left(\boldsymbol{\nu}\atop \boldsymbol{\eta}\right)
\frac{\overline{m}\left[\tfrac{1}{2} \right]_{\left|\boldsymbol{\nu}+\boldsymbol{e}-\boldsymbol{\eta}\right|} }{\boldsymbol{R}^{\boldsymbol{\nu}+\boldsymbol{e}-\boldsymbol{\eta}}}
\,
\frac{C_u \rho^{\left|\boldsymbol{\eta}\right|-1}\left[\tfrac{1}{2} \right]_{\left|\boldsymbol{\eta}\right|} }{\boldsymbol{R}^{\boldsymbol{\eta}}}\\
&\qquad+\overline{m}(\left|\boldsymbol{\nu}+\boldsymbol{e}\right|!)^{\delta-1}
\sum_{\boldsymbol{0}< \boldsymbol{\eta}\leq \boldsymbol{\nu}}
		\left(\boldsymbol{\nu}\atop \boldsymbol{\eta}\right)
\frac{C_u \rho^{\left|\boldsymbol{\nu}-\boldsymbol{\eta}\right|}\left[\tfrac{1}{2} \right]_{\left|\boldsymbol{\nu}+\boldsymbol{e}-\boldsymbol{\eta}\right|} }{\boldsymbol{R}^{\boldsymbol{\nu}+\boldsymbol{e}-\boldsymbol{\eta}}}
\frac{C_u \rho^{\left|\boldsymbol{\eta}\right|-1}\left[\tfrac{1}{2} \right]_{\left|\boldsymbol{\eta}\right|} }{\boldsymbol{R}^{\boldsymbol{\eta}}}\\
&\qquad+
2\overline{u}(\left|\boldsymbol{\nu}+\boldsymbol{e}\right|!)^{\delta-1}
\sum_{\boldsymbol{0}< \boldsymbol{\eta}\leq \boldsymbol{\nu}}
		\left(\boldsymbol{\nu}\atop \boldsymbol{\eta}\right)
\frac{\overline{m}\left[\tfrac{1}{2} \right]_{\left|\boldsymbol{\eta}\right|} }{\boldsymbol{R}^{\boldsymbol{\eta}}}
\frac{C_u \rho^{\left|\boldsymbol{\nu}-\boldsymbol{\eta}\right|}\left[\tfrac{1}{2} \right]_{\left|\boldsymbol{\nu}+\boldsymbol{e}-\boldsymbol{\eta}\right|} }{\boldsymbol{R}^{\boldsymbol{\nu}+\boldsymbol{e}-\boldsymbol{\eta}}}
\\
&\leq
\left(4\overline{m}\,\overline{u}C_u +\overline{m}C_u^2\right)
\frac{\rho^{\left|\boldsymbol{\nu}\right|-1}}{\boldsymbol{R}^{\boldsymbol{\nu}+\boldsymbol{e}}(\left|\boldsymbol{\nu}+\boldsymbol{e}\right|!)^{1-\delta}}
\sum_{\boldsymbol{0}< \boldsymbol{\eta}\leq \boldsymbol{\nu}}
		\left(\boldsymbol{\nu}\atop \boldsymbol{\eta}\right)
\left[\tfrac{1}{2} \right]_{\left|\boldsymbol{\eta}\right|} 
\left[\tfrac{1}{2} \right]_{\left|\boldsymbol{\nu}+\boldsymbol{e}-\boldsymbol{\eta}\right|} \\
&\leq
\left(4\overline{m}\,\overline{u}C_u +\overline{m}C_u^2\right)
\frac{\rho^{\left|\boldsymbol{\nu}\right|-1} \left[\tfrac{1}{2} \right]_{\left|\boldsymbol{\nu}+\boldsymbol{e}\right|} }{\boldsymbol{R}^{\boldsymbol{\nu}+\boldsymbol{e}}(\left|\boldsymbol{\nu}+\boldsymbol{e}\right|!)^{1-\delta}}
.
\end{align*}
For bound of $\mathfrak{N}(\boldsymbol{\nu})$, we notice \eqref{multiindex-est-8} to obtain
\begin{align*}
\mathfrak{N}(\boldsymbol{\nu})
&\leq
(\left|\boldsymbol{\nu}+\boldsymbol{e}\right|!)^{\delta-1}
\Bigg(
\sum_{\boldsymbol{0}< \boldsymbol{\eta}\leq \boldsymbol{\nu}}
\sum_{\boldsymbol{0}< \boldsymbol{\ell} < \boldsymbol{\eta}}
\left(\boldsymbol{\nu}\atop \boldsymbol{\eta}\right)
\left(\boldsymbol{\eta}\atop \boldsymbol{\ell}\right)
\frac{\overline{m}\left[\tfrac{1}{2} \right]_{\left|\boldsymbol{\nu}+\boldsymbol{e}-\boldsymbol{\eta}\right|} }{\boldsymbol{R}^{\boldsymbol{\nu}+\boldsymbol{e}-\boldsymbol{\eta}}}
\frac{C_u \left[\tfrac{1}{2} \right]_{\left|\boldsymbol{\ell}\right|} }{\rho^{1-\left|\boldsymbol{\ell}\right|}\boldsymbol{R}^{\boldsymbol{\ell}}}
\frac{C_u \left[\tfrac{1}{2} \right]_{\left|\boldsymbol{\eta}-\boldsymbol{\ell}\right|} }{\rho^{\left|\boldsymbol{\ell}\right|-\left|\boldsymbol{\eta}\right|{+1}}\boldsymbol{R}^{\boldsymbol{\eta}-\boldsymbol{\ell}}}
\\
&\qquad+
2\sum_{\boldsymbol{0}< \boldsymbol{\eta}\leq \boldsymbol{\nu}}
\sum_{\boldsymbol{0}< \boldsymbol{\ell} < \boldsymbol{\eta}}
\left(\boldsymbol{\nu}\atop \boldsymbol{\eta}\right)
\left(\boldsymbol{\eta}\atop \boldsymbol{\ell}\right)
\frac{\overline{m}\left[\tfrac{1}{2} \right]_{\left|\boldsymbol{\eta}-\boldsymbol{\ell}\right|} }{\boldsymbol{R}^{\boldsymbol{\eta}-\boldsymbol{\ell}}}
\frac{C_u \left[\tfrac{1}{2} \right]_{\left|\boldsymbol{\ell}\right|} }{\rho^{1-\left|\boldsymbol{\ell}\right|}\boldsymbol{R}^{\boldsymbol{\ell}}}
\frac{C_u \left[\tfrac{1}{2} \right]_{\left|\boldsymbol{\nu}+\boldsymbol{e}-\boldsymbol{\eta}\right|} }{\rho^{\left|\boldsymbol{\eta}\right|-\left|\boldsymbol{\nu}\right|}\boldsymbol{R}^{\boldsymbol{\nu}+\boldsymbol{e}-\boldsymbol{\eta}}}
\Bigg)\\
&\leq
3\overline{m}C_u^2(\left|\boldsymbol{\nu}+\boldsymbol{e}\right|!)^{\delta-1}\frac{\rho^{\left|\boldsymbol{\nu}\right|-1}}{\boldsymbol{R}^{\boldsymbol{\nu}+\boldsymbol{e}}}
\sum_{\boldsymbol{0}< \boldsymbol{\eta}\leq \boldsymbol{\nu}}
	\sum_{\boldsymbol{0}< \boldsymbol{\ell} < \boldsymbol{\eta}}
		\left(\boldsymbol{\nu}\atop \boldsymbol{\eta}\right)
		\left(\boldsymbol{\eta}\atop \boldsymbol{\ell}\right)
\left[\tfrac{1}{2} \right]_{\left|\boldsymbol{\nu}+\boldsymbol{e}-\boldsymbol{\eta}\right|} 
\left[\tfrac{1}{2} \right]_{\left|\boldsymbol{\ell}\right|} 
\left[\tfrac{1}{2} \right]_{\left|\boldsymbol{\eta}-\boldsymbol{\ell}\right|} 
\\
&
\leq 
6\overline{m}C_u^2\frac{\rho^{\left|\boldsymbol{\nu}\right|-1} \left[\tfrac{1}{2} \right]_{\left|\boldsymbol{\nu}+\boldsymbol{e}\right|} }{\boldsymbol{R}^{\boldsymbol{\nu}+\boldsymbol{e}}(\left|\boldsymbol{\nu}+\boldsymbol{e}\right|!)^{1-\delta}},
\end{align*}
where in the last step, the double sum is bounded by $2\left[\tfrac{1}{2} \right]_{\left|\boldsymbol{\nu}+\boldsymbol{e}\right|} $, and hence completes the proof.
\end{proof} 
\begin{proof}[Proof of Theorem \ref{gevrey regularity for NS}]
We argue by induction with respect to the order of the derivative $\boldsymbol{\nu}$. For the first-order derivatives we use decomposition $\partial^{\boldsymbol{e}} \boldsymbol{u}=\boldsymbol{w}_{\boldsymbol{e}}^*+\boldsymbol{w}_{\boldsymbol{e}}^{\bot}$, where $\boldsymbol{w}_{\boldsymbol{e}}^*\in\mathcal V$ and $\boldsymbol{w}_{\boldsymbol{e}}^{\bot}\in\mathcal V^{\bot}$. For the bound of $\boldsymbol{w}_{\boldsymbol{e}}^{\bot}$, we use \eqref{wT norm upper} with $\boldsymbol{\nu} = \boldsymbol{0}$  and \eqref{F_1 bound} to obtain
	\begin{align*}
		\left\|{\boldsymbol{w}_{\boldsymbol{e}}^{\bot}}\right\|_{{\mathcal H}}
\leq
\frac{\mathfrak{F}_1(\boldsymbol{0})}{\beta}
=
\frac{1}{\beta}
\left(\left\|{\partial^{\boldsymbol{e}} g}\right\|_{\mathcal L}+
\left\|{\partial^{\boldsymbol{e}} B}\right\|_{\infty}\left\|{\boldsymbol{u}}\right\|_{{\mathcal H}}\right)
\leq
\frac{(\overline{g}+\overline{b}\,\overline{u})}{\beta}
\frac{\left[\tfrac{1}{2} \right]_{1} }{\boldsymbol{R}^{\boldsymbol{e}}}
.
	\end{align*} 
For the bound of $\boldsymbol{w}_{\boldsymbol{e}}^{*}$, we \eqref{wT norm upper}, \eqref{F_2 bound} and notice bound of $\left\|{\boldsymbol{w}_{\boldsymbol{e}}^{\bot}}\right\|_{{\mathcal H}}$ to obtain
\begin{align*}
	\left\|{\boldsymbol{w}_{\boldsymbol{e}}^{*}}\right\|_{{\mathcal H}}
&\leq
\frac{1}{\widetilde{\alpha}}
\left(
\mathfrak{F}_2(\boldsymbol{0})
+
\overline{u}\,
		\overline{m}
	\left\|{\boldsymbol{w}^{\bot}}\right\|_{{\mathcal H}}\right)
\leq
\frac{1}{\widetilde{\alpha}}
\left(
\overline{f}+\overline{b}\,\overline{p}+
\overline{u}
\left(
\overline{a}+\overline{m}\,\overline{u}
\right)
+
\frac{\overline{u}\,\overline{m}(\overline{g}+\overline{b}\,\overline{u})}{\beta}
\right)\frac{\left[\tfrac{1}{2} \right]_{1} }{\boldsymbol{R}^{\boldsymbol{e}}}.
\end{align*}
Notice $\overline{u}=\gamma$, it implies that
\begin{align*}
	\left\|{\partial^{\boldsymbol{e}}\boldsymbol{u}}\right\|_{{\mathcal H}}
\leq
	\left\|{\boldsymbol{w}_{\boldsymbol{e}}^*}\right\|_{{\mathcal H}}
+
	\left\|{\boldsymbol{w}_{\boldsymbol{e}}^{\bot}}\right\|_{{\mathcal H}}
=
\frac{\overline{u}\, \sigma_u \left[\tfrac{1}{2} \right]_{1} }{\boldsymbol{R}^{\boldsymbol{e}}}
=
\frac{C_u  \left[\tfrac{1}{2} \right]_{1} }{\boldsymbol{R}^{\boldsymbol{e}}},
\end{align*}
where we notice the definition of $\sigma_u$ in \eqref{NV sigma-rho-def}. 
For the bound of $\partial^{\boldsymbol{e}} p$, we apply \eqref{p norm upper} with $\boldsymbol{\nu} = \boldsymbol{0}$ and notice bound of $\left\|{\partial^{\boldsymbol{e}}\boldsymbol{u}}\right\|_{{\mathcal H}}$ above to obtain
\begin{align*}
\left\|{\partial^{\boldsymbol{e}} p}\right\|_{\mathcal L}
&\leq
\frac{1}{\beta}
\left(
\mathfrak{F}_2(\boldsymbol{0})
+
\left(\frac{\overline{a}}{2}
+
\overline{u}\,\overline{m}
\right)
\left\|{\partial^{\boldsymbol{e}}\boldsymbol{u}}\right\|_{{\mathcal H}}\right)
\\&\leq
\frac{1}{\beta}
\left(\overline{f}+\overline{b}\,\overline{p}+
\overline{u}
\left(
\overline{a}+\overline{m}\,\overline{u}
\right)
+
\left(\frac{\overline{a}}{2}
+
\overline{u}\,\overline{m}
\right)\overline{u}\sigma_u\right)
\frac{\left[\tfrac{1}{2} \right]_{1} }{\boldsymbol{R}^{\boldsymbol{e}}}
\\
&
=
\frac{\overline{p}\sigma_p \left[\tfrac{1}{2} \right]_{1} }{\boldsymbol{R}^{\boldsymbol{e}}}
=
\frac{C_p \left[\tfrac{1}{2} \right]_{1} }{\boldsymbol{R}^{\boldsymbol{e}}}.
\end{align*}
In the last step, we have used the definition \eqref{NV sigma-rho-def}. 
Thus, we have satisfied the base of induction for the constants $C_u, C_p, \sigma_u$ and $\sigma_p$ defined in \eqref{NV sigma-rho-def}. Now, assuming the validity of \eqref{u H norm} and \eqref{p L norm} for the $\boldsymbol{\nu}$-th derivative, our objective is to establish the same bounds for the $(\boldsymbol{\nu} + \boldsymbol{e})$-th order derivative, where $\boldsymbol{e}$ is a unit multi-index. To bound $\left\|{\partial^{\boldsymbol{\nu}+\boldsymbol{e}}\boldsymbol{u}}\right\|_{{\mathcal H}}$, we denote by $\boldsymbol{w}^*$ and $\boldsymbol{w}^{\bot}$ the projections of $\partial^{\boldsymbol{\nu}+\boldsymbol{e}}\boldsymbol{u}$ into $\mathcal V$ and $\mathcal V^{\bot}$, respectively.

In the upcoming proof, it is crucial to acknowledge that $\rho \geq 1$. Consequently, we will multiply by a sufficient power of $\rho$ on the right-hand side of the bounds to achieve the desired result. Notice \eqref{wT norm upper} and Lemma \ref{lem:F A B M bound}, we get
\begin{align}\label{wT factorial bound}
\left\|{\boldsymbol{w}^{\bot}}\right\|_{{\mathcal H}}
\leq
\frac{1}{\beta}
\left(\mathfrak{F}_1(\boldsymbol{\nu})
		+
		\mathfrak{B}_1(\boldsymbol{\nu})\right)
\leq
\frac{1}{\beta}\left(\overline{g}+\overline{b} \,\overline{u}+2C_u\overline{b}\right)
\frac{\rho^{\left|\boldsymbol{\nu}\right|-1} \left[\tfrac{1}{2} \right]_{\left|\boldsymbol{\nu}+\boldsymbol{e}\right|} }{\boldsymbol{R}^{\boldsymbol{\nu}+\boldsymbol{e}}(\left|\boldsymbol{\nu}+\boldsymbol{e}\right|!)^{1-\delta}}.
\end{align}
By orthogonality of $\mathcal V$ and $\mathcal V^{\bot}$, it implies $\partial^{\boldsymbol{\nu}+\boldsymbol{e}}\boldsymbol{u}=\boldsymbol{w}^*+\boldsymbol{w}^{\bot}$. Notice \eqref{w* norm upper}, \eqref{wT factorial bound}, \eqref{equ: F_2+A+B+M+N bound} and \eqref{NV sigma-rho-def}, we derive
	\begin{align*}
\left\|{\partial^{\boldsymbol{\nu}+\boldsymbol{e}}\boldsymbol{u}}\right\|_{{\mathcal H}}
&\leq
\left\|{\boldsymbol{w}^{\bot}}\right\|_{{\mathcal H}}
+
\left\|{\boldsymbol{w}^{*}}\right\|_{{\mathcal H}}\notag\\
&\leq
\left(
1
+
\frac{\overline{u}\,
		\overline{m}}{\widetilde \alpha}
\right)\left\|{\boldsymbol{w}^{\bot}}\right\|_{{\mathcal H}}
+
\frac{\mathfrak{F}_2(\boldsymbol{\nu}) +		
\mathfrak{A}(\boldsymbol{\nu})
		+
		\mathfrak{B}_2(\boldsymbol{\nu})
		+
		\mathfrak{M}(\boldsymbol{\nu})
+
\mathfrak{N}(\boldsymbol{\nu})}{\widetilde \alpha}
\end{align*}
On the other hand, notice equations \eqref{F_2 bound},\eqref{A bound}, \eqref{B_2 bound}, \eqref{M bound} and \eqref{N bound}, we obtain
\begin{align}\label{equ: F_2+A+B+M+N bound}
	\mathfrak{F}_2(\boldsymbol{\nu}) +		
	&\mathfrak{A}(\boldsymbol{\nu})
	+
	\mathfrak{B}_2(\boldsymbol{\nu})
	+
	\mathfrak{M}(\boldsymbol{\nu})
	+
	\mathfrak{N}(\boldsymbol{\nu})\\
	&\leq
	\left(
		\overline{f}+\overline{b}\,\overline{p}+
		\overline{u}
		\left(
			\overline{a}+\overline{m}\,\overline{u}
		\right)
		2C_u\overline{a}
		+
		2C_p\overline{b}
		+
		4\overline{m}\,\overline{u}C_u
		+
		7\overline{m}C_u^2
	\right)
	\frac{\rho^{\left|\boldsymbol{\nu}\right|-1} \left[\tfrac{1}{2} \right]_{\left|\boldsymbol{\nu}+\boldsymbol{e}\right|} }{\boldsymbol{R}^{\boldsymbol{\nu}+\boldsymbol{e}}(\left|\boldsymbol{\nu}+\boldsymbol{e}\right|!)^{1-\delta}}\notag.
\end{align}
Combining this with the definitions of the constants $C_u,\rho_u$ and $ \rho$ in \eqref{NV sigma-rho-def}, we arrive
\begin{align}
	\left\|{\partial^{\boldsymbol{\nu}+\boldsymbol{e}}\boldsymbol{u}}\right\|_{{\mathcal H}}
	\leq
	C_u \rho_u
	\frac{\rho^{\left|\boldsymbol{\nu}\right|-1} \left[\tfrac{1}{2} \right]_{\left|\boldsymbol{\nu}+\boldsymbol{e}\right|} }{\boldsymbol{R}^{\boldsymbol{\nu}+\boldsymbol{e}}(\left|\boldsymbol{\nu}+\boldsymbol{e}\right|!)^{1-\delta}}.
	\leq
	\frac{C_u \rho^{\left|\boldsymbol{\nu}\right|} \left[\tfrac{1}{2} \right]_{\left|\boldsymbol{\nu}+\boldsymbol{e}\right|} }{\boldsymbol{R}^{\boldsymbol{\nu}+\boldsymbol{e}}(\left|\boldsymbol{\nu}+\boldsymbol{e}\right|!)^{1-\delta}}
	 \label{deri u norm step}
\end{align}
This establishes the desired bound for $\left\|{\partial^{\boldsymbol{\nu}+\boldsymbol{e}}\boldsymbol{u}}\right\|_{{\mathcal H}}$. To derive a similar estimate for $\left\|{\partial^{\boldsymbol{\nu}+\boldsymbol{e}} p}\right\|_{\mathcal L}$, we refer to equations \eqref{p norm upper}, \eqref{equ: F_2+A+B+M+N bound}, and \eqref{deri u norm step}, which yield
\begin{align*}
\left\|{\partial^{\boldsymbol{\nu}+\boldsymbol{e}} p}\right\|_{\mathcal L}
&\leq
\frac{1}{\beta}
\left(
\mathfrak{F}_2(\boldsymbol{\nu}) +		
\mathfrak{A}(\boldsymbol{\nu})
		+
		\mathfrak{B}_2(\boldsymbol{\nu})
		+
		\mathfrak{M}(\boldsymbol{\nu})
+
\mathfrak{N}(\boldsymbol{\nu})\right)
+
\frac{1}{\beta}
\left(
\frac{\overline{a}}{2}
+
\overline{u}\,
\overline{m}\right)
\left\|{\partial^{\boldsymbol{\nu}+\boldsymbol{e}}\boldsymbol{u}}\right\|_{{\mathcal H}}\\
&\leq
C_p \rho_p
\frac{\rho^{\left|\boldsymbol{\nu}\right|-1} \left[\tfrac{1}{2} \right]_{\left|\boldsymbol{\nu}+\boldsymbol{e}\right|} }{\boldsymbol{R}^{\boldsymbol{\nu}+\boldsymbol{e}}(\left|\boldsymbol{\nu}+\boldsymbol{e}\right|!)^{1-\delta}}
\leq
\frac{C_p \rho^{\left|\boldsymbol{\nu}\right|} \left[\tfrac{1}{2} \right]_{\left|\boldsymbol{\nu}+\boldsymbol{e}\right|} }{\boldsymbol{R}^{\boldsymbol{\nu}+\boldsymbol{e}}(\left|\boldsymbol{\nu}+\boldsymbol{e}\right|!)^{1-\delta}}.
\end{align*}
In the last step, we use the definition of $C_p,\rho_p$ and $\rho$ in \eqref{NV sigma-rho-def}. We complete the inductive argument and thereby the proof of the theorem.
\end{proof}
We now possess all the necessary components for the main result of this paper, which is formulated as follows.
\begin{theorem}\label{gevrey regularity for NS original}
		Let the matrices $A,B,M$, the force term $\boldsymbol{f}$ and the source mass term $g$ of \eqref{NS variational form} satisfy  Assumption \ref{small f Assump} and  Assumption \ref{NS Assumption} for some $\delta \geq 1$. Additionally, assume that both Assumption \ref{A,B,M assump} and Assumption \ref{alpha beta assump} also  hold. Then, the solution $(\boldsymbol{u},p)$ to \eqref{NS variational form} is of class Gevrey-$\delta$. More precisely, there {exist} two positive constants $\widetilde  C_u$, $\widetilde  C_p$ and a sequence $\widetilde{\boldsymbol{R}}$ with positive components $\widetilde R_i>0$ such that for all $\boldsymbol{y}\in U$ and $\boldsymbol{\nu}\in \mathcal F$ the derivative of the solution $(\boldsymbol{u},p)$ of \eqref{NS variational form} is bounded by
	\begin{align*}\left\|{\partial^{\boldsymbol{\nu}} \boldsymbol{u}( \boldsymbol{y})}\right\|_{{\mathcal H}}
		\leq
		\frac{\widetilde  C_u (\left|\boldsymbol{\nu}\right|!)^{\delta}}{{\widetilde{\boldsymbol{R}}}^{{\boldsymbol{\nu}}}} 
	\end{align*}
	and
	\begin{align*}\left\|{\partial^{\boldsymbol{\nu}}\, p( \boldsymbol{y})}\right\|_{\mathcal L}
		\leq
		\frac{\widetilde C_p  (\left|\boldsymbol{\nu}\right|!)^{\delta}}{{\widetilde{\boldsymbol{R}}}^{{\boldsymbol{\nu}}}} .
	\end{align*}
\end{theorem}
\begin{proof}
Recall that $\left[\tfrac{1}{2} \right]_{\left|\boldsymbol{\nu}\right|}  \leq |\boldsymbol{\nu}|!$. Then the statement follows directly from the Lemma \ref{u p bound} and equations \eqref{u H norm},\eqref{p L norm} by choosing the constants  $\widetilde C_u=\max\left\{\overline{u},C_u/\rho \right\}$, ${\widetilde C_p}=\max\left\{\overline{p},C_p/\rho \right\}$ and $\widetilde R_i = R_i/ \rho$ with the scaling factors $\rho \geq 1$ from Theorem~\ref{gevrey regularity for NS}.
\end{proof}

\section{Parametric regularity for the plain pullback of the given data} \label{sec: Domain transform}
In this section, we analyze the parametric regularity of the plain pullback transformation building on Example \ref{Ex: dom trans}. {In particular, Theorems \ref{thm:reg_A_B_M} and \ref{thm:reg_f_g} below verify Assumption \ref{NS Assumption} from the previous section.} We assume that the reference domain $\widehat{D}\in \mathbb R^d$, where $d\in \left\{2,3\right\}$, is non-empty, bounded, and has a Lipschitz regular boundary.  Throughout this section, we denote by $\boldsymbol{T}_{\boldsymbol{y}}: \mathrm{cl}\left(\widehat{D}\right) \rightarrow \mathbb R^d$ a bijective bi-Lipschitz map between these Lipschitz domains, i.e. $\boldsymbol{T}_{\boldsymbol{y}}\in W^{1,\infty}\left(\widehat{D},D_{\boldsymbol{y}}\right)$ and $\boldsymbol{T}_{\boldsymbol{y}}^{-1}\in W^{1,\infty}\left(D_{\boldsymbol{y}},\widehat{D}\right)$ for all $\boldsymbol{y}\in U$. Moreover, the random perturbation field $\boldsymbol{T}_{\boldsymbol{y}}$ is parameterized by a countable sequence $\boldsymbol{y}=(y_j)_{{j}\geq 1}\in U$ of independent and identically distributed (i.i.d.) uniform random variables $y_j$. \\
	We define the family of admissible domains $\left\{D_{\boldsymbol{y}}\right\}_{\boldsymbol{y}\in U}$, where each domain is given by $D_{\boldsymbol{y}}=\boldsymbol{T}_{\boldsymbol{y}}\left(\widehat{D}\right)$. The \emph{hold-all domain} is then defined as
\begin{align*}
	\boldsymbol{\rm{D}} =\bigcup_{{\boldsymbol{y}}\in U} D_{\boldsymbol{y}}
	\end{align*}
We show that the condition \eqref{abf assumption} is preserved under domain transformation. This requires establishing the parametric regularity of the matrices $A,B$ and $M$ as well as the input data $\boldsymbol{f}$ and $g$, all of which are affected by the domain transformation. For brevity, in this section, we use the notation $\left\|{\cdot}\right\|_{\infty}$ as in \eqref{norm definition} to denote $\left\|{\cdot}\right\|_{L^{\infty}(\widehat{D})}$. The key modeling assumptions regarding the random perturbation field are outlined below. 

\begin{assumption}\label{domain transform assump}
	We make the following assumptions regarding the vector field $\boldsymbol{T}_{\boldsymbol{y}}: \widehat{D}\rightarrow D_{\boldsymbol{y}}$ for all $\boldsymbol{y}\in U$
	
	\begin{enumerate}[(i)]
		\item \textbf{Boundedness and Invertibility:} the vector fields $\boldsymbol{T}_{\boldsymbol{y}}$ is a continuously differentiable vector field with a bounded inverse. In particular, let  $d\boldsymbol{T}_{\boldsymbol{y}}$ denotes the differential of $\boldsymbol{T}_{\boldsymbol{y}}$, there exists a constant $C_{\boldsymbol{T}}\geq 1$ such that:
			\begin{align}\label{T assumption}
			\max\left\{ \left\|{\boldsymbol{T}_{\boldsymbol{y}}}\right\|_{\infty},
				\left\|{\boldsymbol{T}_{\boldsymbol{y}}^{-1}}\right\|_{\infty},
				\left\|{d\boldsymbol{T}_{\boldsymbol{y}}}\right\|_{\infty},
				\left\|{d\boldsymbol{T}_{\boldsymbol{y}}^{-1}}\right\|_{\infty},
			\right\}\leq C_{\boldsymbol{T}}.
		\end{align}
		\item \textbf{Gevrey-$\delta$ Regularity:} The domain transform $\boldsymbol{T}_{\boldsymbol{y}}$ belongs to the Gevrey class $G^\delta(U, W^{1,\infty}(\widehat{D}))$ with $\delta\geq 1$. This means that for all multi-indices $\boldsymbol{\nu}$, there exists a sequence of positive numbers $\boldsymbol{R}$ such that:
		\begin{align}\label{dT falling factorial assum}
			\max\left\{	\left\|{\partial^{\boldsymbol{\nu}} \boldsymbol{T}_{\boldsymbol{y}}}\right\|_{\infty},\left\|{\partial^{\boldsymbol{\nu}}  d\boldsymbol{T}_{\boldsymbol{y}}}\right\|_{\infty}\right\}
			\leq
			\frac{C_{\boldsymbol{T}}\left[\tfrac{1}{2} \right]_{\left|\boldsymbol{\nu}\right|} }{\boldsymbol{R}^{\boldsymbol{\nu}}}  (\left|\boldsymbol{\nu}\right|!)^{\delta-1}.
		\end{align}
	\end{enumerate}
\end{assumption}
Furthermore, additional assumptions on the force $\boldsymbol{f}$ and source mass terms $g$ are required to ensure the parametric regularity of the domain transformation.
\begin{assumption}\label{RHS assumption}
Let the force term $\boldsymbol{f}: \boldsymbol{\rm{D}} \to \mathbb{R}^d$ and the source mass term $g: \boldsymbol{\rm{D}} \to \mathbb{R}$ be of class Gevrey $G^{\delta}(\boldsymbol{\rm{D}},L^{\infty}(\boldsymbol{\rm{D}}))$. Assume there exist constants $C_{\boldsymbol{f}}, C_g > 0$ and a sequence of positive real numbers $\boldsymbol{\tau} = (\tau_j)_{j=1}^d$ with $\boldsymbol{\tau} \cdot \boldsymbol{1} \geq 1$, such that for all multi-indices $\boldsymbol{\lambda} \in \mathbb{N}_0^d$, the following estimate holds: 
	\begin{align*}
		\max\left\{
			\frac{\left\|{\partial_{\boldsymbol{x}}^{\boldsymbol{\lambda}} \boldsymbol{f}}\right\|_{L^{\infty}(\boldsymbol{\rm{D}})}}{C_{\boldsymbol{f}}},
			\frac{\left\|{\partial_{\boldsymbol{x}}^{\boldsymbol{\lambda}} g}\right\|_{L^{\infty}(\boldsymbol{\rm{D}})}}{C_{g}}
			\right\}
		\leq
		\boldsymbol{\tau}^{\boldsymbol{\lambda}}\,(\left|\boldsymbol{\lambda}\right|!)^{\delta},
	\end{align*}
\end{assumption}
As a direct consequence of the above assumption, the same Gevrey-$\delta$ regularity holds pointwise in $U$, i.e. for every $\boldsymbol{y} \in U$, we have
	\begin{align*}
		\max\left\{
			\frac{\left\|{\partial_{\boldsymbol{x}}^{\boldsymbol{\lambda}} \boldsymbol{f}}\right\|_{L^{\infty}(D_{\boldsymbol{y}})}}{C_{\boldsymbol{f}}},
			\frac{\left\|{\partial_{\boldsymbol{x}}^{\boldsymbol{\lambda}} g}\right\|_{L^{\infty}(D_{\boldsymbol{y}})}}{C_{g}}
		\right\}
	\leq
	\boldsymbol{\tau}^{\boldsymbol{\lambda}}\,(\left|\boldsymbol{\lambda}\right|!)^{\delta}.
\end{align*}
The following theorem provides parametric regularity for the transformed matrices $A,B$ and $M$.
\begin{theorem} \label{thm:reg_A_B_M}
	Under Assumption \ref{domain transform assump}, the matrix $$A(\boldsymbol{y})=d\boldsymbol{T}^{-1}d\boldsymbol{T}^{-\top}\det (d\boldsymbol{T})
	\quad \text{ and }\quad
	B(\boldsymbol{y})=M(\boldsymbol{y})={d\boldsymbol{T}^{-1}}\det (d\boldsymbol{T})$$
	are of Gevrey class $G^{\delta}(U,L^{\infty}(\widehat{D}))$ with $\delta\geq1$, i.e.  for all $\boldsymbol{y}\in U$ and $\boldsymbol{\nu}\in\mathcal F$ the following holds
	\begin{align*}
		\max\left\{\frac{		\left\|{\partial^{\boldsymbol{\nu}}A}\right\|_{\infty}}{\overline{a}},
		\frac{\left\|{\partial^{\boldsymbol{\nu}}B}\right\|_{\infty}}{\overline{b}},
		\frac{\left\|{\partial^{\boldsymbol{\nu}}M}\right\|_{\infty}}{{C_4^{-2}}\overline{m}}
		\right\}
		\leq
		\frac{\rho_{\mathrm{inv}}^{\left|\boldsymbol{\nu}\right|}\left[\tfrac{1}{2} \right]_{\left|\boldsymbol{\nu}\right|} }{\boldsymbol{R}^{\boldsymbol{\nu}}}  (\left|\boldsymbol{\nu}\right|!)^{\delta-1},
	\end{align*}
	where 
	\begin{align*}
		\overline{a}&={4^{d+2}}\,d! {C_{\boldsymbol{T}}^{d+6}},\quad
	\overline{b}={4^{d+1}} \,d! {C_{\boldsymbol{T}}^{d+3}},\\
	{\overline{m}}&={4^{d+1}} C_4^{2} \,d! {C_{\boldsymbol{T}}^{d+3}},\quad
	\rho_{\mathrm{inv}}={3{C_{\boldsymbol{T}}^2}}.
	\end{align*}
\end{theorem}
\begin{proof}
	Since $B=M={d\boldsymbol{T}^{-1}}\det (d\boldsymbol{T})$ and {$A=B\, d\boldsymbol{T}^{-\top}\,$}, we first establish the parametric regularity of $d\boldsymbol{T}^{-1}$ and $J=\det (d\boldsymbol{T})$.
	
	Noting \eqref{T assumption} and \eqref{dT falling factorial assum} and applying the inverse matrix estimate and determinant estimate from Lemma \ref{matrix dev lem}, we obtain
		\begin{align*}
		\left\|{\partial^{\boldsymbol{\nu}}d\boldsymbol{T}^{-1}}\right\|_{\infty}
		\leq 
		\frac{C_{\text{inv}} \rho_{\text{inv}}^{\left|\boldsymbol{\nu}\right|-1} \left[\tfrac{1}{2} \right]_{\left|\boldsymbol{\nu}\right|} }{\boldsymbol{R}^{\boldsymbol{\nu}}}  (\left|\boldsymbol{\nu}\right|!)^{\delta-1},
	\end{align*}
	and
	\begin{align}\label{derivative J bound}
		\left\|{\partial^{\boldsymbol{\nu}}J}\right\|_{\infty}
		\leq 
		\frac{C_J \left[\tfrac{1}{2} \right]_{\left|\boldsymbol{\nu}\right|} }{\boldsymbol{R}^{\boldsymbol{\nu}}}  (\left|\boldsymbol{\nu}\right|!)^{\delta-1},
	\end{align}
	where 
	$$
	C_{\text{inv}}=C_{\boldsymbol{T}}^3,
	\quad \rho_{\text{inv}}=3{C_{\boldsymbol{T}}^2}, \quad 
	C_J={(4 C_{\boldsymbol{T}})^d}\, d!.$$
Noting that $\left\|{\partial^{\boldsymbol{\nu}}d\boldsymbol{T}^{-1}}\right\|_{\infty}=\left\|{\partial^{\boldsymbol{\nu}}d\boldsymbol{T}^{-\top}}\right\|_{\infty}$
for all $\boldsymbol{\nu}\in \mathcal F$, which holds trivially due to Frobenius norm. Applying the product rule from Lemma \ref{matrix dev lem} to  $d\boldsymbol{T}^{-1}d\boldsymbol{T}^{-\top}$, we get
			\begin{align*}
		\left\|{\partial^{\boldsymbol{\nu}}(d\boldsymbol{T}^{-1}d\boldsymbol{T}^{-\top})}\right\|_{\infty}
		\leq 
		\frac{4C_{\text{inv}}^2 \rho_{\text{inv}}^{{\left|\boldsymbol{\nu}\right|}} \left[\tfrac{1}{2} \right]_{\left|\boldsymbol{\nu}\right|} }{\boldsymbol{R}^{\boldsymbol{\nu}}}  (\left|\boldsymbol{\nu}\right|!)^{\delta-1},
	\end{align*}
Now, we establish the Gevrey regularity of $B$. Applying the product rule again, we obtain
		\begin{align*}
		\left\|{\partial^{\boldsymbol{\nu}}B}\right\|_{\infty}
		\leq 
		\frac{	\overline{b}\,  \rho_{\text{inv}}^{{\left|\boldsymbol{\nu}\right|}} \left[\tfrac{1}{2} \right]_{\left|\boldsymbol{\nu}\right|} }{\boldsymbol{R}^{\boldsymbol{\nu}}} (\left|\boldsymbol{\nu}\right|!)^{\delta-1}
{\quad \text{ and } \quad
		\left\|{\partial^{\boldsymbol{\nu}}M}\right\|_{\infty}
		\leq 
		\frac{C_4^{-2}	\overline{m}\,  \rho_{\text{inv}}^{{\left|\boldsymbol{\nu}\right|}} \left[\tfrac{1}{2} \right]_{\left|\boldsymbol{\nu}\right|} }{\boldsymbol{R}^{\boldsymbol{\nu}}} (\left|\boldsymbol{\nu}\right|!)^{\delta-1},}
	\end{align*}
	where $\overline{b}=4C_{\text{inv}} C_J$ {and $\overline{m}=4C_4^2 C_{\text{inv}} C_J$}. Similarly, applying the product rule for {$A=B\, d\boldsymbol{T}^{-\top}$}, we derive
		\begin{align*}
		\left\|{\partial^{\boldsymbol{\nu}}A}\right\|_{\infty}
		\leq 
		\frac{	\overline{a}  \rho_{\text{inv}}^{{\left|\boldsymbol{\nu}\right|}} \left[\tfrac{1}{2} \right]_{\left|\boldsymbol{\nu}\right|} }{\boldsymbol{R}^{\boldsymbol{\nu}}} (\left|\boldsymbol{\nu}\right|!)^{\delta-1},
	\end{align*}
	where $\overline{a}=16C_{\text{inv}}^2 C_J$. This completes the proof.
\end{proof}

We now establish the parametric regularity of the plain pullback force terms $\widehat{\boldsymbol{f}}$ and $\widehat{g}$ as follows 
\begin{lemma}
	Under Assumption \ref{domain transform assump} and Assumption \ref{RHS assumption},	for all $\boldsymbol{\nu}\neq \boldsymbol{0}$ and $\boldsymbol{y}\in U$. Let 	$\widehat{\boldsymbol{f}}(\widehat{\boldsymbol{x}})=\boldsymbol{f}(\boldsymbol{x})=\boldsymbol{f}(\boldsymbol{T}_{\boldsymbol{y}}(\widehat{\boldsymbol{x}}))$
	and $\widehat{g}(\widehat{\boldsymbol{x}})=g(\boldsymbol{x})=g(\boldsymbol{T}_{\boldsymbol{y}}(\widehat{\boldsymbol{x}}))$, the following bounds hold:
	\begin{align}
		\left\|{\partial^{\boldsymbol{\nu}}\widehat{\boldsymbol{f}} }\right\|_{\infty}
		&\leq
		\frac{C_{\widehat{\boldsymbol{f}}}}{2}
		\frac{\rho_{\mathrm{RHS}}^{\left|\boldsymbol{\nu}\right|}}{(2\boldsymbol{R})^{\boldsymbol{\nu}}}
		(\left|\boldsymbol{\nu}\right|!)^{\delta},\label{f hat bound}\\
		\left\|{\partial^{\boldsymbol{\nu}} 
		 \widehat{g}}\right\|_{\infty}
		&\leq
		\frac{C_{\widehat{g}}}{2}
		\frac{\rho_\mathrm{RHS}^{\left|\boldsymbol{\nu}\right|}}{(2\boldsymbol{R})^{\boldsymbol{\nu}}}
		(\left|\boldsymbol{\nu}\right|!)^{\delta},\label{g hat bound}
	\end{align}
	where 
$$
	C_{\widehat{\boldsymbol{f}}}=\frac{C_{\boldsymbol{f}}d\sqrt{d}}{d+1},
	\quad
		C_{\widehat{g}}=\frac{ C_{g}d\sqrt{d}}{d+1},\quad
	\rho_\mathrm{RHS}={ {C_{\boldsymbol{T}}(d+1)\boldsymbol{\tau}\cdot\boldsymbol{1}}}.$$
\end{lemma}
\begin{proof}
	Using definition of $L^{\infty}(\widehat{D})$-norm and applying triangle inequality, we obtain
	$$\left\|{\partial^{\boldsymbol{\nu}} \boldsymbol{f}(\boldsymbol{T}_{\boldsymbol{y}}(\widehat{\boldsymbol{x}}))}\right\|_{\infty}
	=
	\left\|{\left(\sum_{i=1}^d
	(	\partial^{\boldsymbol{\nu}} f_i(\boldsymbol{T}_{\boldsymbol{y}}(\widehat{\boldsymbol{x}})))^2\right)^{\frac{1}{2}}
	}\right\|_{\infty}
	\leq
	\sqrt{d}
	\max_{1\leq i \leq d}\left\|{\partial^{\boldsymbol{\nu}} f_i(\boldsymbol{T}_{\boldsymbol{y}}(\widehat{\boldsymbol{x}}))}\right\|_{\infty}.$$ 
	For any $i$ such that $1\leq i \leq d$, we apply Fa{\`a} di Bruno formula to obtain
	\begin{align*}
		\partial^{\boldsymbol{\nu}} f_i(\boldsymbol{T}_{\boldsymbol{y}}(\widehat{\boldsymbol{x}}))
		=
		\sum_{1\leq\left|\boldsymbol{\lambda}\right|\leq \left|\boldsymbol{\nu}\right|}
		\partial^{\boldsymbol{\lambda}}_{\boldsymbol{x}}f_i(\boldsymbol{T}_{\boldsymbol{y}}(\widehat{\boldsymbol{x}}))
		\sum_{s=1}^n
		\sum_{P_s(\boldsymbol{\nu},\boldsymbol{\lambda})}
		\boldsymbol{\nu}!
		\prod_{j=1}^s
		\frac{(\partial^{\boldsymbol{\ell}_j}\boldsymbol{T}_{\boldsymbol{y}}(\widehat{\boldsymbol{x}}))^{\boldsymbol{k}_j}}{\boldsymbol{k}_j! (\boldsymbol{\ell}_j!)^{\left|\boldsymbol{k}_j\right|}}.
	\end{align*}
	Noting that $\boldsymbol{T}_{\boldsymbol{y}}(\widehat{\boldsymbol{x}})=(T_{\boldsymbol{y},m}(\widehat{\boldsymbol{x}}))_{m=1}^d$, for all $m=1,\dots,d$ and multi-index $\boldsymbol{\ell}_j$, we have
	$$
	\left\|{\partial^{\boldsymbol{\ell}_j}T_{\boldsymbol{y},m}(\widehat{\boldsymbol{x}})}\right\|_{\infty}
	\leq
	\left\|{\partial^{\boldsymbol{\ell}_j}\boldsymbol{T}_{\boldsymbol{y}}(\widehat{\boldsymbol{x}})}\right\|_{\infty}
	\leq
	\frac{C_{\boldsymbol{T}}}{2} \frac{(\left|\boldsymbol{\ell}_j\right|!)^{\delta}}{(2\boldsymbol{R})^{\boldsymbol{\ell}_j}}.
	$$
	{Taking the $L^{\infty}(\widehat{D})$-norm} both sides and applying the triangle inequality, we derive
	\begin{align*}
		\left\|{\partial^{\boldsymbol{\nu}} f_i(\boldsymbol{T}_{\boldsymbol{y}}(\widehat{\boldsymbol{x}}))}\right\|_{\infty}
		&\leq
		\sum_{1\leq\left|\boldsymbol{\lambda}\right|\leq \left|\boldsymbol{\nu}\right|}
		\left\|{\partial^{\boldsymbol{\lambda}}_{\boldsymbol{x}}f_i(\boldsymbol{x})}\right\|_{L^{\infty}(D_{\boldsymbol{y}})}\\
		&\qquad\qquad\qquad\quad \times
		\sum_{s=1}^n
		\sum_{P_s(\boldsymbol{\nu},\boldsymbol{\lambda})}
		\boldsymbol{\nu}!
		\prod_{j=1}^s
		\frac{{\displaystyle \prod_{m=1}^d}\left\|{\partial^{\boldsymbol{\ell}_j}T_{\boldsymbol{y},m}(\widehat{\boldsymbol{x}})}\right\|_{\infty}^{k_{j,m}}}{\boldsymbol{k}_j! (\boldsymbol{\ell}_j!)^{\left|\boldsymbol{k}_j\right|}}\\
		&\leq
		\sum_{1\leq\left|\boldsymbol{\lambda}\right|\leq \left|\boldsymbol{\nu}\right|}
		C_{\boldsymbol{f}}\,\boldsymbol{\tau}^{\boldsymbol{\lambda}}\,(\left|\boldsymbol{\lambda}\right|!)^{\delta}
		\sum_{s=1}^n
		\sum_{P_s(\boldsymbol{\nu},\boldsymbol{\lambda})}
		\boldsymbol{\nu}!
		\prod_{j=1}^s
		\frac{{\displaystyle \prod_{m=1}^d}\left(	\frac{C_{\boldsymbol{T}}}{2} \frac{(\left|\boldsymbol{\ell}_j\right|!)^{\delta}}{(2\boldsymbol{R})^{\boldsymbol{\ell}_j}}\right)^{k_{j,m}}}{\boldsymbol{k}_j! (\boldsymbol{\ell}_j!)^{\left|\boldsymbol{k}_j\right|}}.
	\end{align*}
Factorizing constant and applying Lemma \ref{lem: n k!} with $\omega=\delta-1{\geq 0}$, we derive
		\begin{align*}
		\left\|{\partial^{\boldsymbol{\nu}} f_i(\boldsymbol{T}_{\boldsymbol{y}}(\widehat{\boldsymbol{x}}))}\right\|_{\infty}
\leq
				\frac{C_{\boldsymbol{f}} \left(C_{\boldsymbol{T}}\,\boldsymbol{\tau}\cdot\boldsymbol{1}\right)^{\left|\boldsymbol{\nu}\right|}(\left|\boldsymbol{\nu}\right|!)^{\delta-1}	}{2\,(2\boldsymbol{R})^{\boldsymbol{\nu}}}	
				\sum_{1\leq\left|\boldsymbol{\lambda}\right|\leq \left|\boldsymbol{\nu}\right|\atop \boldsymbol{\lambda} \in \mathbb N^d}
				\left|\boldsymbol{\lambda}\right|!
				\sum_{s=1}^n
				\sum_{P_s(\boldsymbol{\nu},\boldsymbol{\lambda})}
				\boldsymbol{\nu}!
				\prod_{j=1}^s
				\frac{(\left|\boldsymbol{\ell}_j\right|!)^{\left|\boldsymbol{k}_j\right|}}{\boldsymbol{k}_j! (\boldsymbol{\ell}_j!)^{\left|\boldsymbol{k}_j\right|}}.
	\end{align*}
	Employing Lemma \ref{singlevariable FdB lemma}, we obtain
	\begin{align*}
		\left\|{\partial^{\boldsymbol{\nu}} \boldsymbol{f}(\boldsymbol{T}_{\boldsymbol{y}}(\widehat{\boldsymbol{x}}))}\right\|_{\infty}
		&\leq
		\sqrt{d}
		\max_{1\leq i \leq d}\left\|{\partial^{\boldsymbol{\nu}} f_i(\boldsymbol{T}_{\boldsymbol{y}}(\widehat{\boldsymbol{x}}))}\right\|_{\infty}
		\leq
		\sqrt{d}\left(\frac{C_{\boldsymbol{f}}d}{2(d+1)}
		\left(\frac{ {C_{\boldsymbol{T}}(d+1)\boldsymbol{\tau}\cdot\boldsymbol{1}}}{2\boldsymbol{R}}\right)^{\boldsymbol{\nu}}	(\left|\boldsymbol{\nu}\right|!)^{\delta}	\right).
	\end{align*}
	Noting that $C_{\widehat{\boldsymbol{f}}}=\frac{ C_{\boldsymbol{f}}d\sqrt{d}}{{2}(d+1)}$ and $\rho_\text{RHS}={ {C_{\boldsymbol{T}}(d+1)\,\boldsymbol{\tau}\cdot\boldsymbol{1}}}$, this completes the proof for the force term $\boldsymbol{f}$. Similarly, we can derive the same result for the source mass term with a different constant $C_{\widehat{g}}=\frac{C_{g}d\sqrt{d}}{{2}(d+1)}$ .
\end{proof}
We now can derive parametric regularity for the transformed function $\boldsymbol{f}, g$. It is as follows
\begin{theorem} \label{thm:reg_f_g}
	Under Assumption \ref{domain transform assump} and Assumption \ref{RHS assumption}. Let 	$\widetilde{\boldsymbol{f}}(\boldsymbol{y})=\widehat{\boldsymbol{f}}\, \det (d \boldsymbol{T}_{\boldsymbol{y}})$
	and $\widetilde{g}(\boldsymbol{y})=\widehat{g} \, \det (d \boldsymbol{T}_{\boldsymbol{y}})$, the following bound holds
	\begin{align*}
		\max\left\{\frac{		\left\|{\partial^{\boldsymbol{\nu}}\widetilde{\boldsymbol{f}}}\right\|_{\infty}}{\overline{f}},
			\frac{\left\|{\partial^{\boldsymbol{\nu}}\widetilde{g}}\right\|_{\infty}}{\overline{g}}
		\right\}
		\leq
		\frac{\rho_\text{RHS}^{\left|\boldsymbol{\nu}\right|}\left[\tfrac{1}{2} \right]_{\left|\boldsymbol{\nu}\right|} }{\boldsymbol{R}^{\boldsymbol{\nu}}}  (\left|\boldsymbol{\nu}\right|!)^{\delta-1},
	\end{align*}
	where 
	$$
	 \overline{f}={4^{d+1}}\, C_{\boldsymbol{f}} C_{\boldsymbol{T}}^d \frac{d\sqrt{d}}{d+1}\, d!,\,
\quad
\overline{g}={4^{d+1}}\, C_{g} C_{\boldsymbol{T}}^d \frac{d\sqrt{d}}{d+1}\, d!.
	$$
\end{theorem}
\begin{proof}
	Noting \eqref{f hat bound}, \eqref{g hat bound} and \eqref{ff-estimates}, we obtain the following bounds for all $\boldsymbol{\nu}\neq \boldsymbol{0}$
	\begin{align*}
		\left\|{\partial^{\boldsymbol{\nu}}\widehat{\boldsymbol{f}} }\right\|_{\infty}
&\leq
C_{\widehat{\boldsymbol{f}}}\,
\frac{  \rho_{\text{RHS}}^{\left|\boldsymbol{\nu}\right|}\left[\tfrac{1}{2} \right]_{\left|\boldsymbol{\nu}\right|} }{\boldsymbol{R}^{\boldsymbol{\nu}}}
(\left|\boldsymbol{\nu}\right|!)^{\delta-1},\\
\left\|{\partial^{\boldsymbol{\nu}} 
	\widehat{g}}\right\|_{\infty}
&\leq
C_{\widehat{g}}\,
\frac{  \rho_{\text{RHS}}^{\left|\boldsymbol{\nu}\right|}\left[\tfrac{1}{2} \right]_{\left|\boldsymbol{\nu}\right|} }{\boldsymbol{R}^{\boldsymbol{\nu}}}
(\left|\boldsymbol{\nu}\right|!)^{\delta-1},
	\end{align*}
	 By applying the product rule from Lemma \ref{matrix dev lem} and using the bound in \eqref{derivative J bound}, while noting that $\rho_{\text{RHS}}\geq 1$, we complete the proof.
\end{proof}

\section{Applications and numerical experiments}\label{sec: NS num exp}
In this section, we conduct a series of numerical experiments to approximate the expected values of solutions to the Navier-Stokes equations under domain uncertainty. To achieve this, we employ Gauss-Legendre quadrature and QMC methods. The uncertain domains are parameterized as
$$
D_{\boldsymbol{y}}=\left\{\boldsymbol{x}\in \mathbb R^2: 0\leq x_1 \leq 1 \, \text{ and }\, 0\leq x_2 \leq   T(x_1,\boldsymbol{y})\right\},
$$
here, the perturbation function $T(x_1,\boldsymbol{y})$ would be changed in following subsections. The reference domain is fixed as the unit square, $\widehat{D}= [0,1]^2$. The source mass term $g\equiv 0$ and the force term is $\boldsymbol{f}=-\Delta \boldsymbol{w} + \nabla q$, where
\begin{align*}
	\boldsymbol{w}(\boldsymbol{x})
	=
	\left(\sin^2(\pi x_1)\sin(2 \pi x_2)
		\atop
		-\sin(2\pi x_1)\sin^2(\pi x_2)
	\right)
	\quad \text{ and } \quad
	q(\boldsymbol{x})
	=
	\sin(2\pi x_1)\sin(2\pi x_2).
\end{align*}
The matrices $A,B$ and $M$ are defined as in Example \ref{Ex: dom trans} where the domain transform $\boldsymbol{T}_{y}: \widehat{D}\mapsto D_{\boldsymbol{y}}$ and its differential are given as following expressions
\begin{align*}
\boldsymbol{T}_{\boldsymbol{y}}(\widehat{x}_1,\widehat{x}_2)
=
 \begin{pmatrix}
	\widehat{x}_1 \\
\,	\widehat{x}_2T(\widehat{x}_1,\boldsymbol{y})\,
\end{pmatrix}
\quad
\text{ and }\quad
d\boldsymbol{T}_{\boldsymbol{y}}(\widehat{x}_1,\widehat{x}_2)
=
\begin{pmatrix}
	1 & 0 \\
	\widehat{x}_2\widehat{\partial}_1 T(\widehat{x}_1,\boldsymbol{y})\,
	&
	T(\widehat{x}_1,\boldsymbol{y})
\end{pmatrix}
,
\end{align*}
where $\widehat{\partial}_1$ denotes the partial derivative with respect to the variable $\widehat{x}_1$.

\subsection{Gauss-Legendre quadrature} \label{sec: Gauss-Legendre Quadrature}
 In this section, we compare the convergence behaviour for the analytic ($\delta = 1$) and Gevrey non-analytic ($\delta > 1$) cases, which validate the theoretical predictions
of Theorems \ref{gevrey regularity for NS original}. 
 The domain perturbation function $T(x_1,y)$ is one of the following functions
\begin{align}
T^{(1)}(x_1,y)&=1+0.15\cos(5 \pi \, x_1\, y ) ,
\label{def: phi 1}\\
T^{(2)}(x_1,y)&=1+0.35\exp\left(-\,\frac{1}{\sin^2(5 \pi \, x_1\, y )}\right).
\label{def: phi 2}
\end{align}
Here, $y$ is a scalar real random variable uniformly distributed in $[-1,1]$. It is worth noting that $T^{(1)}$ is analytic for all $y$ in $[-1,1]$. Our current focus is the integral of the squared $\mathcal L$-norm of the pressure
\begin{equation}\label{NS Elambda}
I({p}) =  \int_{-1}^1 \left\|{p}\right\|_{\mathcal L}^2 \, dy.
\end{equation}
Since the pressure $p$ is not available in closed form, the integral \eqref{NS Elambda} can be approximated by numerical quadrature, e.g. the Gauss-Legendre quadrature which will be efficient for the case of a single real-valued parameter. Let $\{\xi_i,w_i\}_{i=1}^n$ be the nodes and weights of the $n$-point Gauss-Legendre quadrature
\begin{equation}
Q_n[{p}] =  \sum_{i=1}^n w_i \left\|{p(\xi_i)}\right\|_{\mathcal L}^2.
\end{equation}
We are interested in the behaviour of the quadrature error
\begin{equation}
\varepsilon_n = \frac{\left|I(p) - Q_n[p]\right|}{\left|I(p)\right| }
\end{equation}
with increasing $n$.  It is known that the convergence of $\varepsilon_n$ is strongly related to the regularity of $p$ with respect to $y$. In particular, \cite[Theorem 5.2]{ChernovSchwab2012} implies that
\begin{equation}
\varepsilon_n \leq C \exp( - r n^{1/\delta}).
\end{equation}
with positive constants $C$ and $r$ independent of $n$.
\begin{itemize}
\item[(1)]
In the case of the analytic perturbation $T^{(1)}$ as in \eqref{def: phi 1}, Theorem \ref{gevrey regularity for NS original} implies that $\left\|{p}\right\|_{\mathcal L}$ is analytic in $y$, i.e. $\delta = 1$, and therefore we expect
\begin{equation}\label{eps1}
\varepsilon_n^{(1)} \leq C \exp (-rn).
\end{equation}
\item[(2)] The perturbation $T^{(2)}$ is not analytic at $y=0$, but is Gevrey-$\delta$ uniformly for all $y \in [-1,1]$ with $\delta \geq \frac{3}{2}$, see \cite{ChenRodino1996} and \cite[Section 6]{ChernovSchwab2012}. Theorem \ref{gevrey regularity for NS original} implies that $\left\|{p}\right\|_{\mathcal L}$ is Gevrey-$\delta$ with the same $\delta = \frac{3}{2}$ and hence we expect
\begin{equation}\label{eps2}
\varepsilon_n^{(2)} \leq C \exp (-rn^{2/3}).
\end{equation}
\end{itemize}
\begin{figure}[h]
\centerline{
\begin{tikzpicture}
\pgfdeclareimage[width=\textwidth]{navier_stoke_NS_GL}{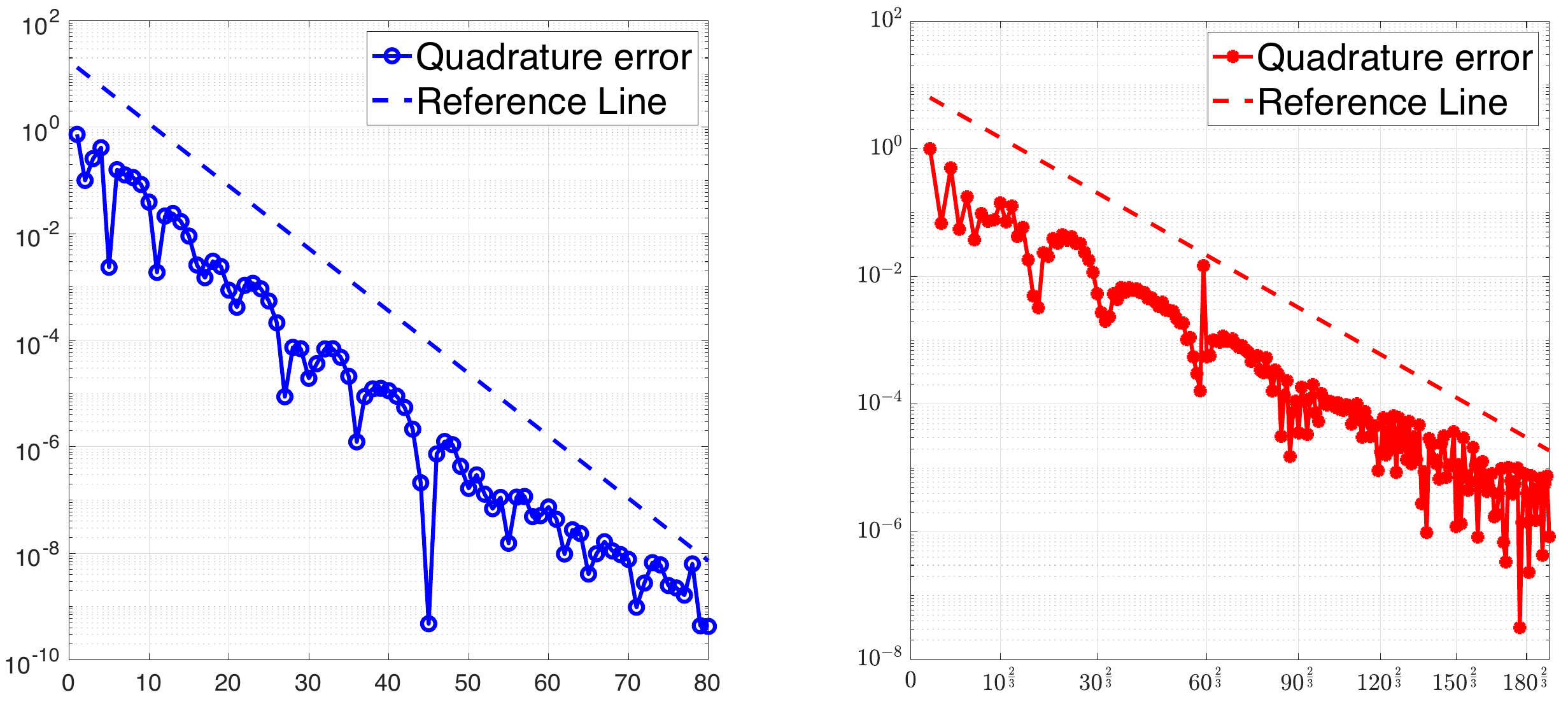};
\node (navier_stoke_NS_GL)
at (0,0) {\pgfuseimage{navier_stoke_NS_GL}};
\end{tikzpicture}
}
	\caption{The absolute error of Gauss-Legendre quadrature: Navier-Stokes equation}{Quadrature error $\varepsilon^{(1)}_n$ (left) with respect to the number $n$ of quadrature points and Quadrature error $\varepsilon^{(2)}_n$ (right) with respect to $m = n^{2/3}$.}\label{NS fig:GL}
\end{figure}

\begin{figure}[h]
	\centerline{
		\begin{tikzpicture}
\node at (0,0) {\includegraphics[width=0.5\textwidth]{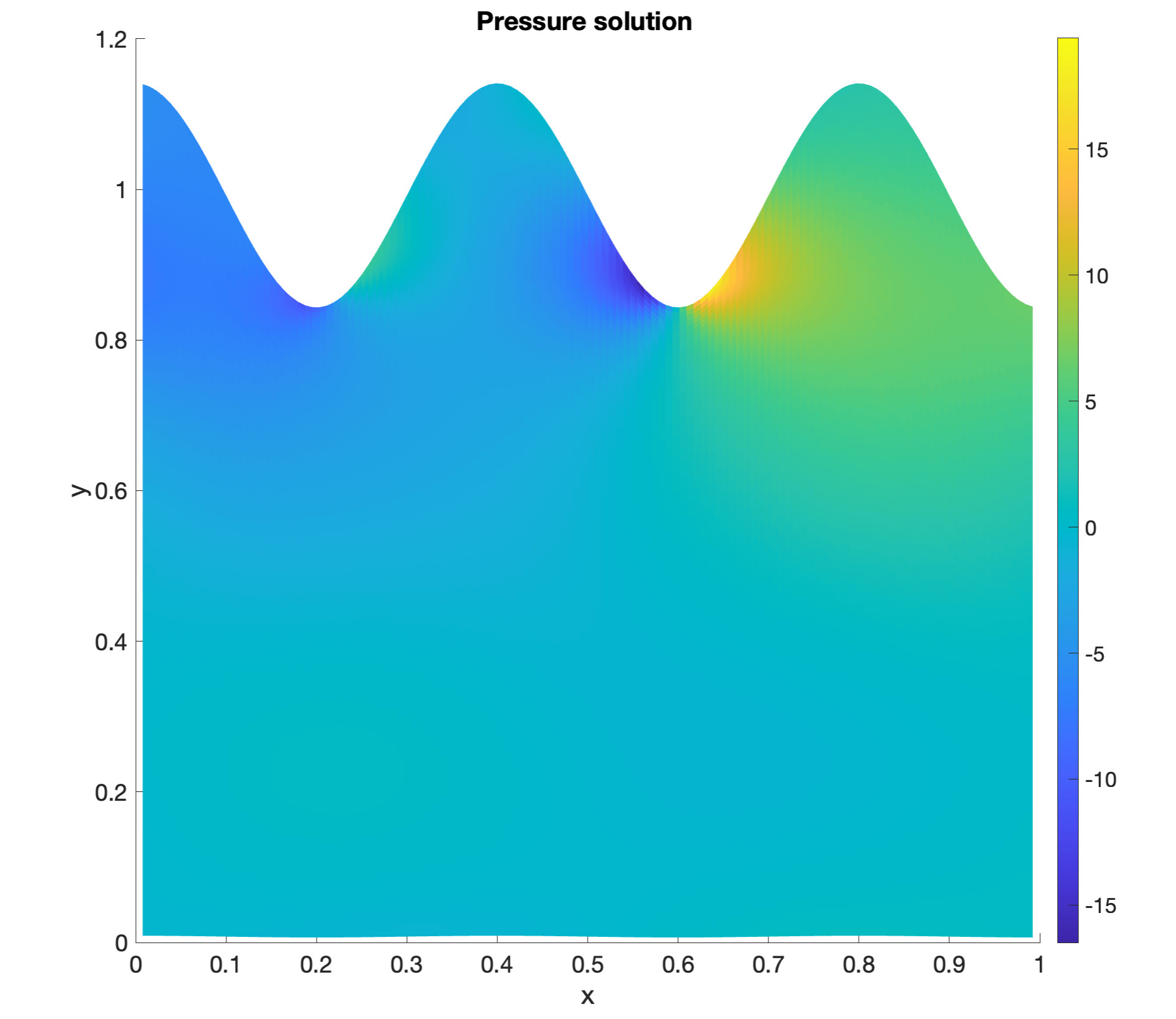}};
\node at (7,0) {\includegraphics[width=0.5\textwidth]{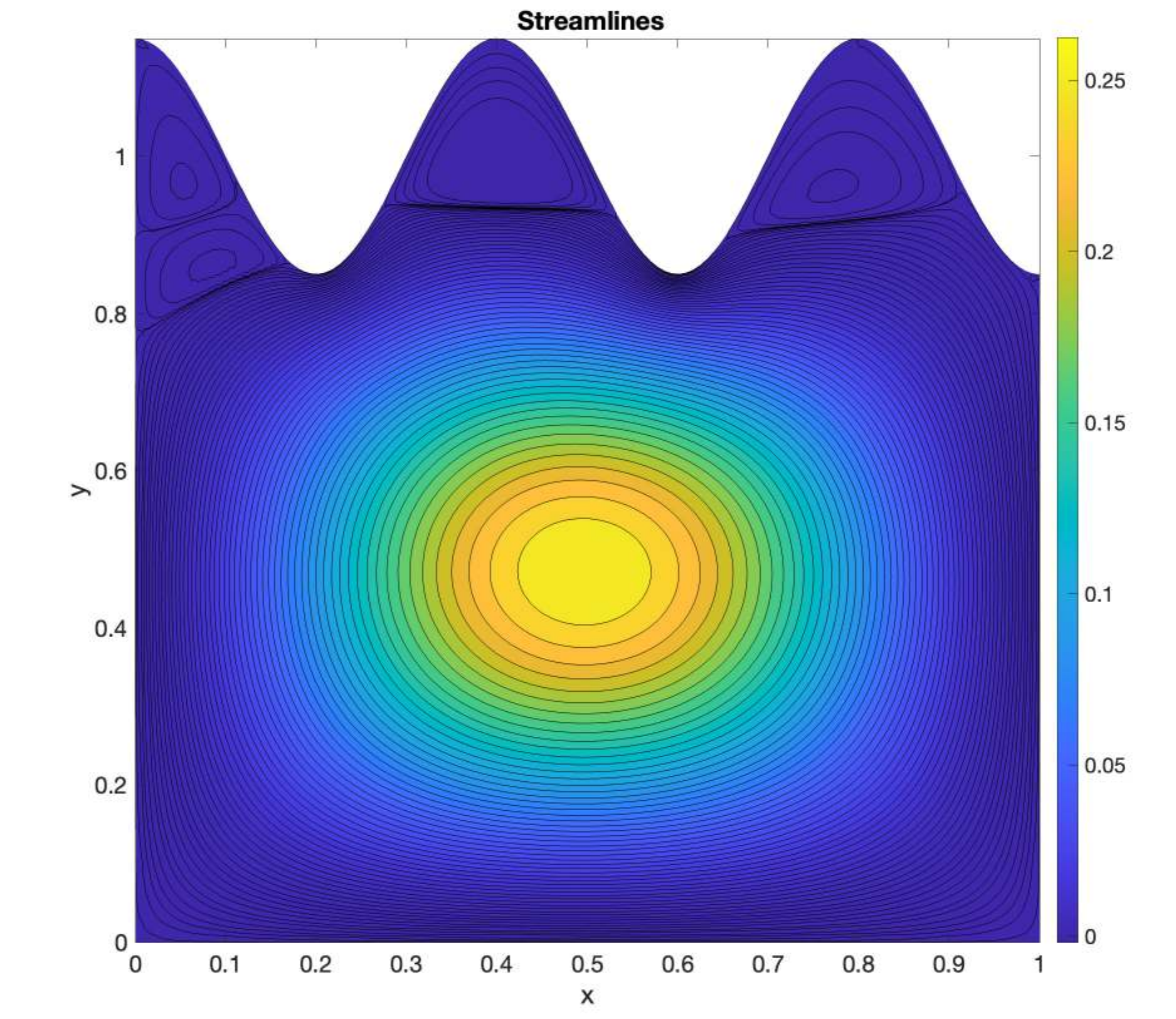}};
		\end{tikzpicture}
	}
	\caption{Pressure and streamlines of velocity field solution of Navier-Stokes equation for $T^{(1)}$ at $y=1$}{}\label{NS fig:pressure 1}
\end{figure}

\begin{figure}[h]
	\centerline{
		\begin{tikzpicture}
\node at (7,0) {\includegraphics[width=0.5\textwidth]{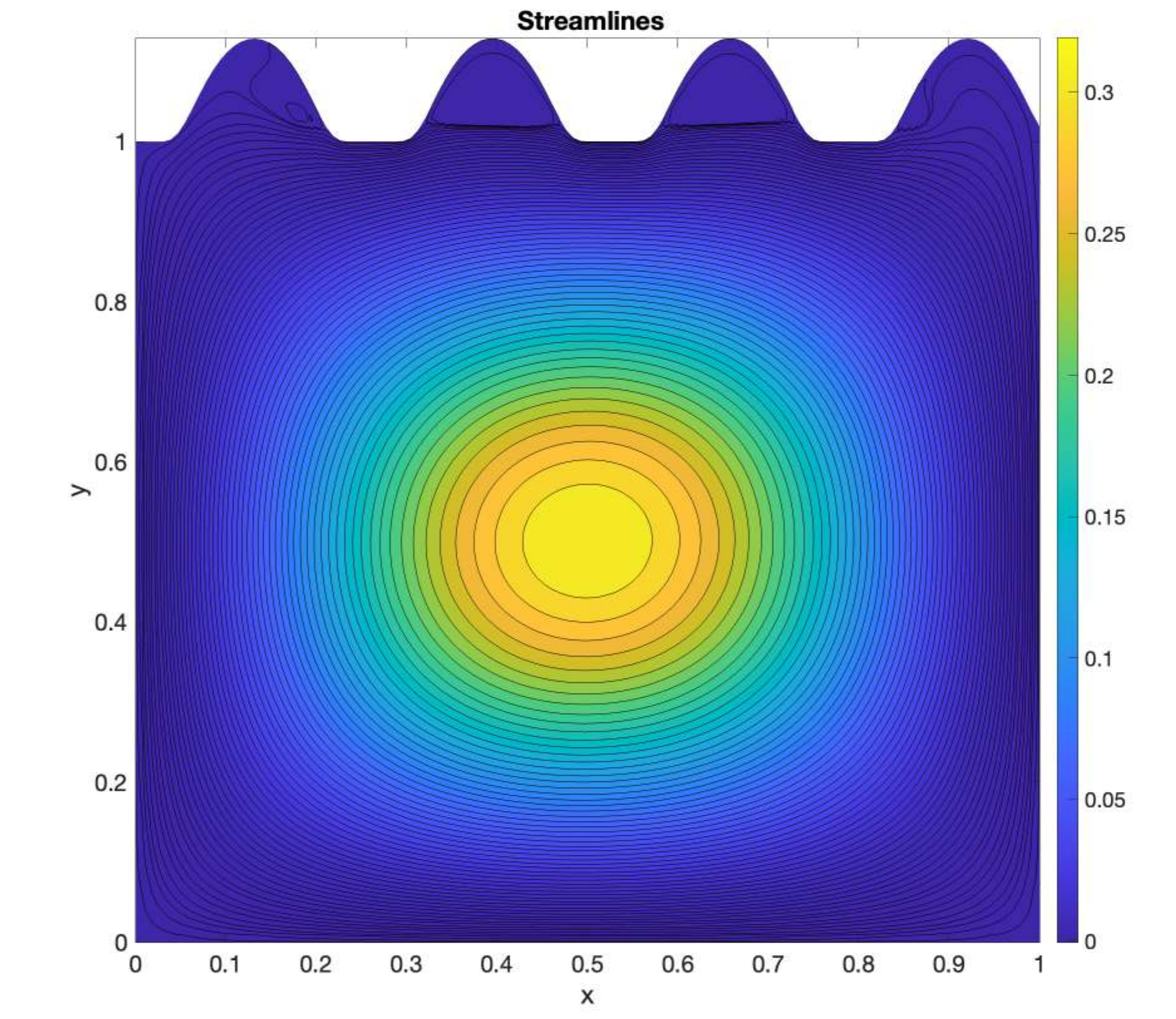}};
			\node at (0,0) {\includegraphics[width=0.5\textwidth]{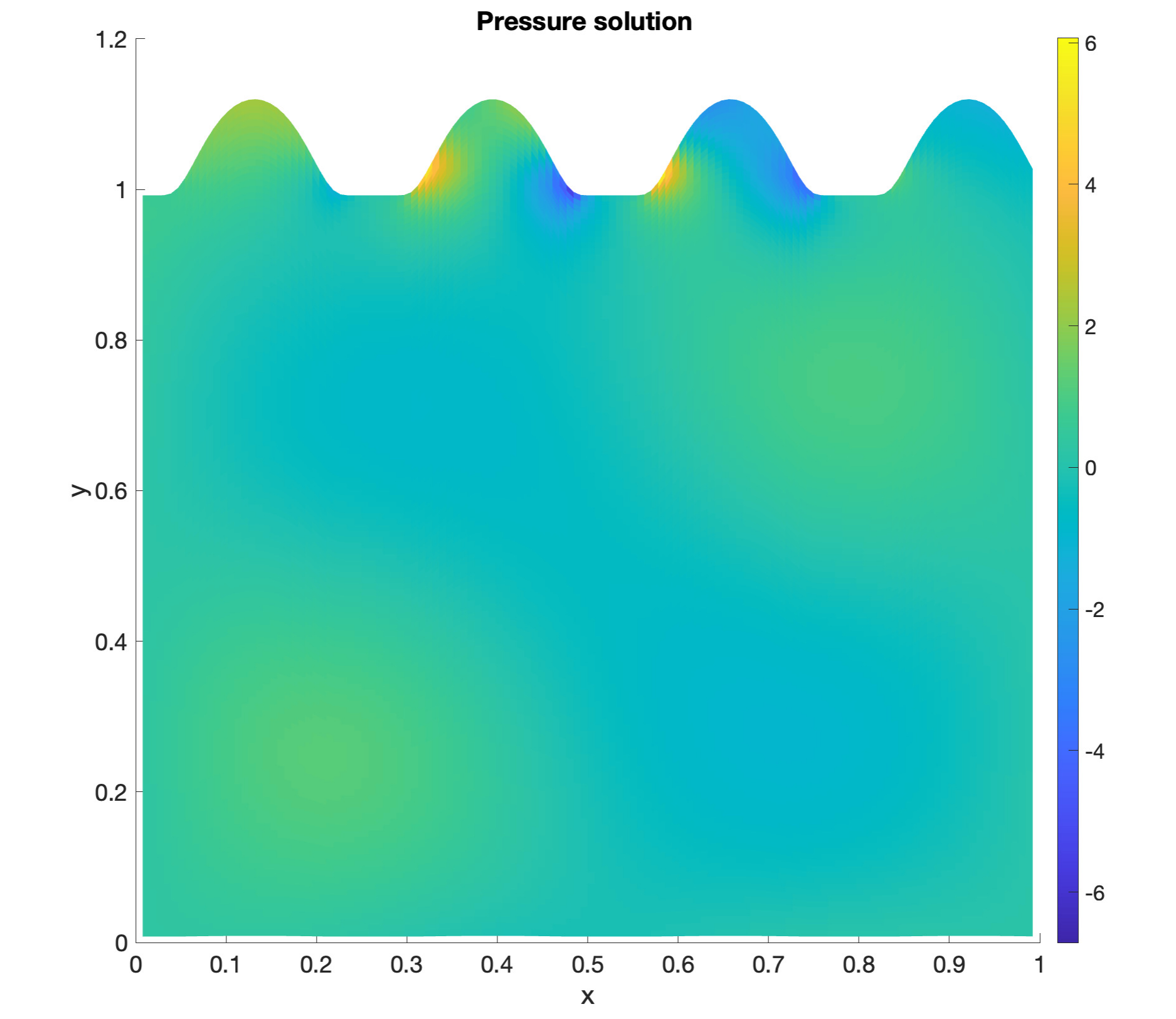}};
		\end{tikzpicture}
	}
	\caption{Pressure and streamlines of velocity field solution of Navier-Stokes equation for $T^{(2)}$ at $y=0.76$}{}\label{NS fig:pressure 2}
\end{figure}

In order to observe the behaviour predicted in \eqref{eps1} and \eqref{eps2} we solve deterministic equations~\eqref{NS variational form} in every quadrature point $y = \xi_i$ on a very fine finite element grid having $16.641$ degrees of freedom. Since $I(p)$ is not available in closed form, we approximate it by the very fine Gauss-Legendre quadrature $Q_{n^*}[p]$ with $n^* = 200$ quadrature nodes for both case $T^{(1)}$ and $T^{(2)}$. As a solver we use Finite Element Method (FEM) based on the software package IFISS, see e.g. \cite{ifiss, ers14, ers07}, with the jump-stabilized $Q_1$-$P_0$ method on a uniform grid. For solving nonlinear system, we used fixed-point iteration as in \eqref{contracting mapping} with nonlinear tolerance $\tau= 10^{-7}$. In both cases, due to the FEM and the nonlinear solver, achieving an error smaller than $10^{-10}$ is not feasible. However, the relative error is sufficient to demonstrate the convergence rate of Gauss-Legendre quadrature.

In {the left panel of} Figure \ref{NS fig:GL}, we plot the relative error $\varepsilon_n^{(1)}$ against the number of quadrature points $n$ in the semi-logarithmic scale. The reference line clearly shows the linear trend of the type $-r n + \log C$ and thereby confirms \eqref{eps1}.
 
In {the right panel of Figure \ref{NS fig:GL}} , we plot the relative error $\varepsilon_n^{(2)}$ with respect to $m := n^{2/3}$ in the semi-logarithmic scale. Here we can also observe the linear trend of the type $-r m + \log C$. This confirms \eqref{eps2} and thereby demonstrates the meaning and validity of Theorem \ref{gevrey regularity for NS original}.

\subsection{Quasi-Monte Carlo method for Gevrey functions} \label{sec: QMC}
In this section we give an example for the case of high-dimensional parametric integration and demonstrate an application of Theorem \ref{gevrey regularity for NS original}. For this, we consider the variational problem \eqref{NS variational form} with $s$ parameters $\boldsymbol{y} = (y_1,\dots,y_s)$ and the quantity of interest being the integral
\begin{align}\label{approx integral}
	I_s( {F}) =	
	\int_{U}  {F}(\boldsymbol{y}) \, d y_1 \dots dy_s,
\end{align}
where the functional $ {F}(\boldsymbol{y}):=\left\|{\boldsymbol{u}(\boldsymbol{y})}\right\|_{{\mathcal H}}^2$ is the square ${\mathcal H}$-norm of the velocity field $\boldsymbol{u}$ and the $s$-dimensional unit cube $U = [-\frac{1}{2},  \frac{1}{2}]^s$ is the parameter domain. Such problems arise when the coefficients of \eqref{NS variational form} are modelled as general random fields e.g. via their Karhunen-L\`oeve expansion \cite{HansenSchwab13,CohenDevoreSchwab2011,KuoSchwabSloan2013,KuoSchwabSloan2012} and then truncated to a certain dimension $s$ for computational reasons. For further comments on the analysis of the truncation error for Gevrey-class parametrizations we refer to \cite[Section 6.2]{ChernovLe2024a} and \cite[Section 5.2]{DjurdjevacKaarnioja2025}.

In this section we fix the dimension of the parameter domain at $s = 100$ and concentrate on estimation of $I_s(\boldsymbol{u})$ where $\boldsymbol{u}$ is a solution for the variational formulation \eqref{NS variational form} with the perturbation function $T(x_1,\boldsymbol{y})$ being one of following functions
\begin{equation}\label{QMC-a1-new}
 T^{(1)}({x_1},\boldsymbol{y})=\exp\left( \sum_{j=1}^{100} \frac{1}{4 j^3} \sin(5 j \pi {x_1}) \, y_j\right)
\end{equation}
or
\begin{equation}\label{QMC-a2-new}
 T^{(2)}({x_1},\boldsymbol{y})=1+ \sum_{j=1}^{100} \frac{\exp(1)}{ 8 j^6} \sin((j+5) \pi {x_1}) \,\left( \exp \left(-\,\frac{1}{(y_j+0.5)^2}\right)+1\right).
\end{equation}
Here $y_j$ are independent uniformly distributed in $[-\frac{1}{2},\frac{1}{2}]$ random variables for all $1 \leq j\leq s$. Notice that $T^{(1)}$ is analytic with respect to $\boldsymbol{y}$ with $\delta^{(1)} = 1$, whereas $T^{(2)}$ is of Gevrey-$\delta$ class with $\delta^{(2)} = \frac{3}{2}$,  {i.e. $T^{(1)} \in G^1(U,{\mathcal H})$ and $ T^{(2)} \in G^{\frac{3}{2}}(U,{\mathcal H})$}. Notice also that both $T^{(1)}$ and $T^{(2)}$ are $\ell^p$-summable with $p<1$. More precisely, with any $p^{(1)} > \frac{1}{2}$ and $p^{(2)} > \frac{1}{5}$. This will be important for the convergence of the QMC methods, as explained below.

From Theorem \ref{gevrey regularity for NS original}, we know that $ {F}(\boldsymbol{y})$ is Gevrey-$\delta$ regular with the same $\delta$ as for $T^{(1)}$ and $T^{(2)}$. The quantity of interest for the perturbation $T_k$ will be denoted by $I_s( {F}^{(k)})$, $k = 1,2$. The standard Monte Carlo(MC) quadrature is the sample mean approximation $Q^{MC}_{s,n}( {F}) := \frac{1}{n} \sum_{i=1}^n  {F}(\boldsymbol{y}_s^{(i)})$, where $\boldsymbol{y}_s^{(i)}$ are independent samples uniformly distributed in $U$. It satisfies the classical error estimate for the relative root mean square error
\begin{equation}\label{errorMC}
	\varepsilon^{\textrm{MC},(k)}_n := \sqrt{\mathbb E\left( \left|\frac{I_s( {F}^{(k)})-Q^{MC}_{s,n}( {F}^{(k)})}{I_s( {F}^{(k)})}\right|^2\right)}
	\lesssim {n^{-\frac{1}{2}}}.
\end{equation} 
Observe that the asymptotic error bound is insensitive to the Gevrey class regularity of $ {F}$. As an alternative to the MC quadrature, we consider a class of  {QMC} rules called randomly shifted rank-1 lattice rules defined as
\begin{align}\label{QMC quad def}
	Q^{\Delta}_{s,n}( {F}) := 
	\frac{1}{n} \sum_{i=1}^n  {F}\big(\left\{ \tfrac{i \boldsymbol{z}_{s}}{n} +\Delta  \right\} -\tfrac{1}{2}\big).
\end{align}
Here the braces in \eqref{QMC quad def} indicate the fractional part of each component of the argument vector, see for e.g. \cite{KuoNuyens2016}, \cite{KuoSchwabSloan2012}. The quadrature points are constructed using a generating vector $\boldsymbol{z}_{s} \in \mathbb N^s$, and a random shift $\Delta$, which is a random variable uniformly distributed over the cube $(0,1)^s$. As quadrature point generator, we use the software package QMC4PDE \cite{KuoNuyens2016software} to obtain the generating vector for both cases. According to~\cite{Cools2006}, the generating vector can be reused for any intermediate power of 2, forming an embedded sequence of lattice rules. This allows the quadrature rules at the lower QMC levels to be derived from the highest-level rule. We adopt this embedded lattice rule approach to reduce computational time, as demonstrated below
\begin{align*}
	\underbrace{
		\underbrace{
		\underbrace{
		\underbrace{\boldsymbol{t}_0,\boldsymbol{t}_1\,}_{\text{Level 1}}
	,\boldsymbol{t}_2,\boldsymbol{t}_3\,}_{\text{Level 2}}
	,\boldsymbol{t}_4,\boldsymbol{t}_5,\boldsymbol{t}_6,\boldsymbol{t}_7\,}_{\text{Level 3}}
	,\boldsymbol{t}_8,\boldsymbol{t}_9,\boldsymbol{t}_{10},\boldsymbol{t}_{11}
	,\boldsymbol{t}_{12},\boldsymbol{t}_{13},\boldsymbol{t}_{14},\boldsymbol{t}_{15}\,}_{\text{Level 4}}
	,\dots,
	\end{align*}
	where $\boldsymbol{t}_i=\left\{ \tfrac{i \boldsymbol{z}_{s}}{n} +\Delta  \right\} -\tfrac{1}{2}$. Since $Q^{\Delta}_{s,n}({F})$ is also a random variable, we use the relative root mean square error (RMSE) as a measure of accuracy
\begin{align}\label{error-deco}
\varepsilon^{QMC,(k)}_n := \sqrt{\mathbb{E}\left(\bigg|\frac{I_s( {F}^{(k)})-Q^{\Delta}_{s,n}( {F}^{(k)})}{I_s( {F}^{(k)})}\bigg|^2\right)},
\end{align}
where $\mathbb E$ stands for the expectation with respect to the random shifts $\Delta$. The analysis for RMSE, when the integrand is analytic (of Gevrey class with $\delta=1$) can be found  {in} \cite{KuoNuyens2016} and \cite{KuoSchwabSloan2012}. In \cite[Lemma 7.4]{ChernovLe2024a} we extended this result to integrands that may belong to Gevrey-$\delta$ class with $\delta > 1$. According to that, the RMSE of QMC quadrature for a fixed integer $s$ and $n$ being a power of 2 admits the bound
\begin{equation}\label{errorQMC}
\varepsilon^{QMC,(k)}_n
	\leq C_{s,\vartheta} n^{-\frac{1}{2\vartheta}},
\end{equation}
where $C_{s,\vartheta}$ is independent of $n$ and
\begin{align*}
		\vartheta = 
		\left\{\begin{matrix}
			\omega &\text{for any } \omega \in (\tfrac{1}{2},1) & \text{when } p\in (0,\tfrac{2}{3\delta}],
			\\
			\frac{\delta p}{2-\delta p}    &  & \text{when } p\in (\tfrac{2}{3\delta},\tfrac{1}{\delta}] .
		\end{matrix}\right. 
	\end{align*}

\begin{figure}[h]
	\centerline{\includegraphics[width=0.5\textwidth]{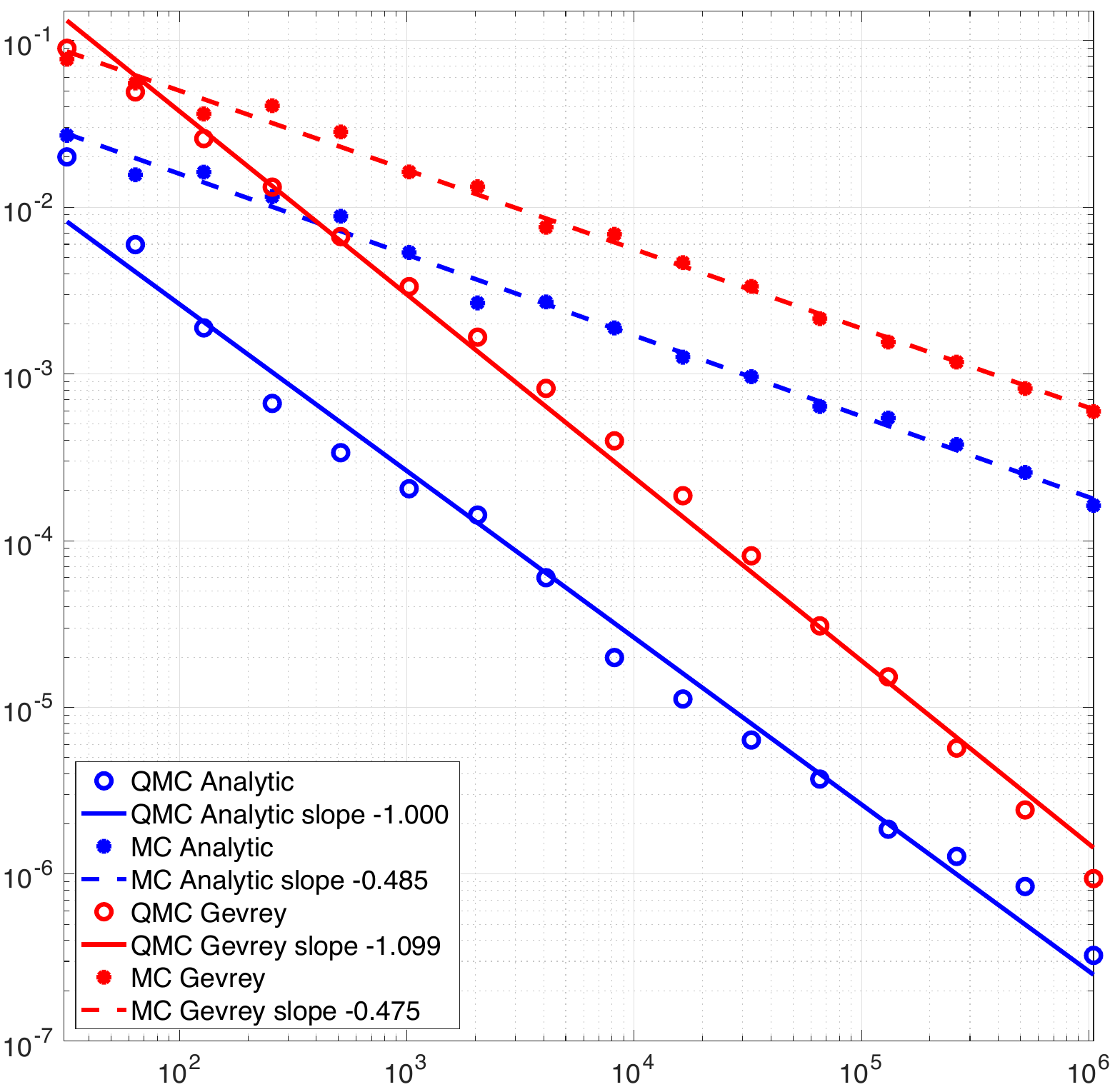}}
	\caption{Convergence of the quadrature error with respect to the number of samples $n$ for four methods: QMC analytic ($\varepsilon^{{\rm QMC},(1)}_n$), QMC Gevrey ($\varepsilon^{{\rm QMC},(2)}_n$), MC analytic ($\varepsilon^{{\rm MC},(1)}_n$), MC Gevrey ($\varepsilon^{{\rm MC},(2)}_n$).}\label{fig2}
\end{figure}

Recall that the summability parameter $p^{(k)}$ can be chosen to satisfy $p^{(k)} <\frac{2}{3 \delta^{(k)}}$ with $k= 1,2$ and therefore $\vartheta$ can be close to~$\frac{1}{2}$ and therefore $\varepsilon_n^{{\rm {QMC}},(k)}$ is close to $n^{-1}$	. This effect is very well visible in the numerical experiments in Figure \ref{fig2}. These results were obtained by solving the variational problem \eqref{NS variational form} in each quadrature point by the {FEM} on a uniform mesh with the mesh size $h=1/128$ and by means of a fixed-point iteration \eqref{contracting mapping} with the error tolerance of $10^{-7}$ with respect to the ${\mathcal H}$-norm. 

The outer expectation in \eqref{errorQMC} is approximated by the empirical mean of $R = 16$ runs, i.e., for $\Delta^{(j)}$, $1 \leq j \leq R$ being an independent sample from the uniform distribution from the unit cube $(0,1)^s$ and $Q^{(j)}_{s,n}( {F}^{(k)})$ the corresponding QMC quadrature. Thus, we approximate the relative QMC error by 
\begin{equation}
	\varepsilon_n^{{\rm QMC},(k)}\sim
	\sqrt{
		\frac{1}{R} \sum_{j=1}^R
		\left(\left|\frac{ I_s^*( {F}^{(k)})- {Q}^{(j)}_{s,n}( {F}^{(k)})}{I_s^*( {F}^{(k)})}\right|^2\right)}
		\end{equation}
		and analogously for the standard  {MC} approximation $\varepsilon_n^{{\rm MC},(k)}$.
We choose the reference value $I_s^*( {F}^{(k)})$ as QMC approximation with a very high number of points, namely $n=2^{16}$ with 16 random shifts, in both cases $k=1,2$.  {
It is worth noting that the relative error $\varepsilon_n^{{\rm QMC},(k)}$ is independent of dimension $s$, as shown in, for example, \cite[Lemma 7.4]{ChernovLe2024a}.}

Figure \ref{fig2} clearly demonstrates that predicted MC and QMC convergence rates of $n^{-\frac{1}{2}}$ and about $n^{-1}$, cf. \eqref{errorMC} and \eqref{errorQMC}, are very well reproduced in the numerical experiments.  By choosing  the QMC rule with the largest number of quadrature points as the reference value $I_s^*( {F})$ in both QMC and MC methods, and given that the MC method converges at a rate of approximately $n^{-\frac{1}{2}}$ in both cases, we ensure that the QMC method indeed converges to $I_s(F)$.

\backmatter

\bmhead{Acknowledgements}

The authors acknowledge the HPC cluster ROSA at the University of Oldenburg (Germany) for providing computing resources. ROSA was funded by the German Research Foundation (DFG) through the Major Research Instrumentation Programme (INST 184/225-1 FUGG) and by the Ministry of Science and Culture (MWK) of Lower Saxony.

\end{document}